\documentclass[leqno,12pt]{article}
\usepackage[american]{babel}
\usepackage{amsmath} 
\usepackage{latexsym}
\usepackage{amssymb}
\usepackage{theorem}
\usepackage{diagrams}
\theoremstyle{change}
\newtheorem{Thm}{Theorem}[section]
\newtheorem{Cor}[Thm]{Corollary}
\newtheorem{Prop}[Thm]{Proposition}
\newtheorem{Lem}[Thm]{Lemma}
\newtheorem{Problem}[Thm]{Problem}
{\theorembodyfont{\rmfamily}
\newtheorem{Example}[Thm]{Example}
\newtheorem{Num}[Thm]{}
\newtheorem{Rem}[Thm]{Remark}
\newtheorem{Def}[Thm]{Definition}}

\renewcommand{\phi}{\varphi}
\newcommand{\rk}{\mathrm{rk}}
\newcommand{\crk}{\mathrm{cor}}
\newcommand{\bra}[1]{\langle#1\rangle}
\newcommand{\proof}{\par\medskip\rm\emph{Proof. }}

\newcommand{\qed}{\ \hglue 0pt plus 1filll $\Box$}
\newcommand{\mapstoo}{\longmapsto}
\newcommand{\RR}{\mathbb{R}}
\newcommand{\ZZ}{\mathbb{Z}}
\newcommand{\NN}{\mathbb{N}}
\renewcommand{\SS}{\mathbb{S}}

\newcommand{\CC}{\mathbb{C}}
\newcommand{\HH}{\mathbb{H}}
\newcommand{\OO}{\mathbb{O}}
\newcommand{\bi}{\boldsymbol i}
\newcommand{\bj}{\boldsymbol j}
\newcommand{\bk}{\boldsymbol k}
\newcommand{\bl}{\boldsymbol\ell}

\newcommand{\FF}{\mathbb F}
\newcommand{\id}{\mathrm{id}}
\newcommand{\lk}{\mathrm{lk}}
\newcommand{\Hom}{\mathrm{Hom}}
\newcommand{\Aut}{\mathrm{Aut}}
\newcommand{\AutTop}{\mathrm{AutTop}}
\newcommand{\RRP}{{\RR\mathrm P}}
\newcommand{\CCP}{{\CC\mathrm P}}
\newcommand{\OOP}{{\OO\mathrm P}}
\newcommand{\Cham}{\mathrm{Cham}}
\newcommand{\SKIP}[1]{}
\newcommand{\eps}{\varepsilon}
\renewcommand{\emptyset}{\varnothing}
\renewcommand{\setminus}{-}

\newcommand{\fg}{\mathfrak{g}}
\newcommand{\fa}{\mathfrak{a}}
\newcommand{\fp}{\mathfrak{p}}
\newcommand{\fb}{\mathfrak{b}}
\newcommand{\fq}{\mathfrak{q}}
\newcommand{\fd}{\mathfrak{d}}
\newcommand{\fh}{\mathfrak{h}}
\newcommand{\fn}{\mathfrak{n}}
\newcommand{\fso}{\mathfrak{so}}
\newcommand{\fsu}{\mathfrak{su}}
\newcommand{\fsp}{\mathfrak{sp}}
\newcommand{\fu}{\mathfrak{u}}
\newcommand{\ft}{\mathfrak t}
\newcommand{\fc}{\mathfrak c}
\newcommand{\fz}{\mathfrak z}
\newcommand{\pr}{\mathrm{pr}}
\newcommand{\Cen}{\mathrm{Cen}}
\newcommand{\Lie}{\mathrm{Lie}}
\newcommand{\SO}{\mathrm{SO}}
\newcommand{\SU}{\mathrm{SU}}
\newcommand{\Sp}{\mathrm{Sp}}
\newcommand{\Spin}{\mathrm{Spin}}
\newcommand{\cG}{\mathcal G}
\newcommand{\colim}{\varinjlim}
\begin{document}
\title{{\boldmath\bf Homogeneous compact geometries}}
\author{Linus Kramer and Alexander Lytchak 
\thanks{This work was supported by the SFB 878 'Groups, Geometry and Actions'. 
The second author
was supported by a Heisenberg Fellowship}}
\maketitle

\begin{abstract}
We classify compact homogeneous geometries of irreducible spherical type
and rank at least $2$ which admit a transitive action of a compact
connected group, up to equivariant $2$-coverings. We apply our classification
to polar actions on compact symmetric spaces.
\end{abstract}
We classify compact homogeneous geometries which look locally
like compact spherical buildings. Geometries which look locally like
buildings arise naturally in various recognition problems in group theory. 
Tits' seminal paper \emph{A local approach to buildings} \cite{TitsLocal}
is devoted to them. Among other things, Tits proved there
that a geometry which looks locally like a building can be $2$-covered
by a building if and only if a local geometric obstruction vanishes.
The condition is that the links of all corank $3$ simplices of 
type $\mathsf C_3$ and $\mathsf H_3$ admit coverings by buildings.

There exists 
a famous finite geometry of type $\mathsf C_3$ which is not covered by any
building, the so-called Neumaier Geometry \cite{Neumaier}
(see also \ref{NeumaierGeometry} below).
It seems to be an open problem if there exist other (finite) examples
of non-building $\mathsf C_m$ geometries,
and if there exist geometries of type $\mathsf H_3$
(note that we assume geometries to be thick).
Assuming a transitive group action, Aschbacher classified
all finite homogeneous geometries of type $\mathsf C_3$,
see  \cite{Asch}  and \cite{Yos}. Using this result, Aschbacher
classified the finite homogeneous geometries with irreducible spherical diagrams
\cite[Thm.~3]{Asch}.
Our Theorem A below may be viewed a Lie group analog of his classification.
More results and references can be found in Pasini's book \cite{Pasini}.

We are here concerned with the classification of geometries on which compact
Lie groups act transitively. Such geometries arise in the 
classification of polar actions. 
For example Thorbergsson's classification of
isoparametric submanifolds in spheres \cite{Tho} relied heavily
on the Burns-Spatzier classification of compact connected spherical
buildings admitting a strongly transitive action \cite{BS}. 
In the last section 
we describe an application of our results to polar actions
on compact symmetric spaces.
\begin{center}
 {*}\hspace{1cm}{*}\hspace{1cm}{*}
\end{center}
Our results are as follows.
With Grundh\"ofer and Knarr, the first author obtained in the mid-90s
a complete classification of irreducible homogeneous compact connected
spherical buildings. We recall the result (which is built on earlier work of
Salzmann, L\"owen, and Burns-Spatzier \cite{BlauesBuch} \cite{BS}).

\smallskip\noindent\textbf{Theorem \cite{GKK1,GKK2,GKVW}}
{\em Let $\Delta$ be a compact building of
irreducible spherical type and rank at least $2$,
with connected panels. Assume that its topological 
automorphism group acts transitively on 
the chambers of $\Delta$. Then $\Delta$ is the spherical building 
associated to a noncompact simple Lie group.}

\smallskip
Using a combination of this and results in
Tits' local approach \cite{TitsLocal},
we prove the following main results.
The exceptional $\mathsf C_3$ geometry that appears in Theorem A
was discovered by
Podest\`a-Thorbergsson \cite{PTh}. We describe this geometry in detail in
Section \ref{TheExceptionalC3Geometry}. For the unexplained definitions we
refer to our paper below, in particular to Section~1 and 2.
A map between geometries is called a $2$-covering if it is bijective
on the links of all simplices of corank $2$, see \ref{k-CoveringDefinition} below.

\smallskip\noindent\textbf{Theorem A}
{\em Let $\Delta$ be a compact geometry of
irreducible spherical type and rank at least $2$,
with connected panels. Assume that a compact
group acts continuously and transitively on the chambers of $\Delta$.

If $\Delta$ is not of type $\mathsf C_3$, then there exists a simple noncompact
Lie group $S$, a compact chamber-transitive subgroup $K\subseteq S$ and a $K$-equivariant 
$2$-covering $\widetilde\Delta\rTo\Delta$, where $\widetilde\Delta$ is
the canonical spherical building associated to $S$.

If $\Delta$ is of type $\mathsf C_3$, then either there exists
a building $\widetilde\Delta$ and a $2$-covering $\widetilde\Delta\rTo\Delta$ as in the
previous case, 
or $\Delta$ is isomorphic to the unique exceptional homogeneous compact
$\mathsf C_3$ geometry which cannot be $2$-covered by  any building.}

\smallskip
More general results are proved in \ref{MainTheoremInTheReducibleCase},
\ref{MainTheoremInTheIrreducibleCase},
\ref{C3ClassificationTheorem}, \ref{UniquenessOfException}.
In this way we obtain a complete classification of homogeneous compact geometries
with connected panels whose irreducible factors are of spherical type and rank at least
$2$, up to equivariant $2$-coverings. 
In certain situations the conclusion of Theorem A may be strengthened. For example,
we prove the following in \ref{AutomaticallyBuilding}.

\smallskip\noindent\textbf{Proposition}
{\em Let $\Delta$ be a homogeneous compact geometry as in Theorem A.
If $\Delta$ is of type $\mathsf A_m$ or $\mathsf E_6$ or
if all panels are $2$-dimensional, then $\Delta$ is the building associated to
a noncompact simple Lie group $S$, and the compact connected group induced on 
$\Delta$ is a maximal compact subgroup of $S$.}

\smallskip
One application of this classification is the following result,
which builds heavily on results by the second author \cite{Lytchak}.
See \ref{cons} below for more details, an outline of proof, and how this relates
to independent work on the classification of polar actions in positive
curvature by Fang-Grove-Thorbergsson \cite{FGTh}.

\smallskip\noindent
\textbf{Theorem B}
\emph{Suppose that $G\times X\rTo X$ is a polar action of a compact connected Lie group $G$
on a symmetric space $X$ of compact type. 
Then, possibly after replacing $G$ by a larger orbit equivalent
group, we have splittings $G=G_1 \times\cdots\times G_m$ and $X=X_1  \times\cdots\times X_m$, 
such that the action of $G_i$ on $X_i$ is either trivial or hyperpolar or the
space $X_i$ has rank $1$, for $i=1,\ldots m$.}
\begin{center}
 {*}\hspace{1cm}{*}\hspace{1cm}{*}
\end{center}
The following problems seem to be open.

\smallskip\noindent\textbf{Problem~1}
{\em Are there geometries of type $\mathsf H_3$, $\mathsf H_4$, 
$\mathsf F_4$, or $\mathsf C_m$, $m\geq 4$, that are not $2$-covered by buildings?}
Note that we assume geometries to be thick. Possibly, in some cases affirmative answers can be 
obtained along the lines of \cite{BK}.

\smallskip\noindent\textbf{Problem~2}
{\em Is the topological automorphism group of a compact geometry of spherical type 
locally compact in the compact-open topology?} Our approach avoids this question. However,
we show that the compact groups that appear are automatically Lie groups, provided that the
panels are connected.

\smallskip\noindent\textbf{Problem~3}
{\em Does the conclusion of Theorem A still hold if we just assume that the topological
automorphism group acts transitively on the chambers?}

\smallskip\noindent\textbf{Problem~4}
{\em Are there non-homogeneous compact geometries of irreducible spherical type and
rank at least
$3$ that are not $2$-covered by buildings?} We remark that non-homogeneous compact geometries
which are $2$-covered by buildings arise naturally from polar foliations, see for example~\cite{DV}.

\smallskip\noindent\textbf{Problem~5}
{\em Is the exceptional $\mathsf C_3$ geometry from Section \ref{TheExceptionalC3Geometry}
simply connected? Is there an analogy with the Neumaier Geometry? Can this geometry be
defined over other fields?}

\begin{center}
 {*}\hspace{1cm}{*}\hspace{1cm}{*}
\end{center}
The paper is organized as follows. In Section~1 we introduce the relevant
combinatorial notions and explain Tits' results. In Section~2 we introduce
a convenient category of homogeneous compact geometries, and we show the
existence of universal objects. In Section~3 we review the known examples
of homogeneous compact geometries of type $\mathsf C_3$, and in Section~4
we prove that this list is complete. In the final Section~5 we combine our 
classification results and prove, among other things, 
Theorem A and explain the main steps for the proof of Theorem B.

\medskip\noindent\textbf{Acknowledgment}
We thank Karl Hofmann, Misha Kapovich, 
Andreas Kollross, Karl-Her\-mann Neeb, Antonio Pasini and the anonymous referee
for helpful comments.

\section{Geometries and buildings}
In order to make this paper self-contained, we first introduce some 
elementary combinatorial terminology.
For the following facts and definitions we refer to Tits \cite{TitsLNM,TitsLocal}.
Additional material can be found in \cite{Bue1,Bue2,Pasini}. The geometries
which we consider here are a standard tool in the structure theory
of the finite simple groups.
We allow ourselves a few small deviations from Tits' terminology. These will
be indicated where they appear.
\begin{Num}\textbf{Chamber complexes\ }
\label{ChamberComplexDefinition}
Let $V$ be a (nonempty) set and $\Delta$ a collection of finite subsets of $V$.
If $V$ is closed under going down 
(i.e. $\alpha\subseteq\beta\in\Delta$ implies $\alpha\in\Delta$) and if
$V=\bigcup\Delta$, then the poset $(\Delta,\subseteq)$ is called a
\emph{simplicial complex}. The elements of $V$ are called \emph{vertices}
and the elements of $\Delta$ are called \emph{simplices}. 
The \emph{rank} of a simplex $\alpha$ is the number of its vertices,
$\rk(\alpha)=\mathrm{card}(\alpha)$.
Two vertices which are contained in a simplex are called \emph{adjacent}.
If the simplex $\alpha$ is contained in the simplex $\beta$, we call
$\alpha$ a \emph{face} of $\beta$.
A simplicial complex is called a \emph{flag complex} if every finite set
of pairwise adjacent vertices is a simplex
('every non-simplex contains a non-edge').
The \emph{link} of a 
simplex $\alpha$ is the subcomplex 
\[
\lk_\Delta(\alpha)=\lk(\alpha)=\{\beta\in\Delta\mid \alpha\cap\beta=\emptyset\text{ and }
\alpha\cup\beta\in\Delta\}
 \]
and the \emph{residue} of $\alpha$ is the set
$\Delta_{\geq\alpha}=\{\beta\in\Delta\mid \beta\supseteq\alpha\}$
of all simplices having $\alpha$ as a face.
The link and the residue of $\alpha$ are poset-isomorphic. 
We note that $\Delta$ is the link of the empty simplex.

A \emph{simplicial map} between simplicial complexes is a map between their vertex
sets which maps simplices to simplices.
We call a simplicial map \emph{regular} if its restriction to every simplex is
bijective;
these are Tits' \emph{morphisms} \cite[1.1]{TitsLNM}. The \emph{geometric realization}
$|\Delta|$ of $\Delta$ consists of all functions $\xi:V{\rTo}[0,1]$
whose support $\mathrm{supp}(\xi)=\{v\in V\mid \xi(v)>0\}$
is a simplex, and with $\sum_{v\in V}\xi(v)=1$. 
We also write
\[
\xi=\sum_{v\in V}v\cdot\xi(v).
\]
The weak
topology turns $|\Delta|$ into a CW complex 
which we denote by $|\Delta|_w$.
A simplicial map $f:\Delta\rTo\Delta'$ induces (by piecewise linear continuation)
a continuous map 
$|\Delta|_w\rTo|\Delta'|_w$ which we denote by the same symbol $f$.

A simplicial complex is called \emph{pure} if every simplex is contained in a
maximal simplex and if all maximal simplices have the same rank $n$. In this case
we say that $\Delta$ has rank $n$ and we call the simplices of rank $n$ \emph{chambers}.
The set of all chambers is denoted $\Cham(\Delta)$. The \emph{corank} of a simplex $\alpha$ 
is then defined as 
\[
\crk(\alpha)=n-\rk(\alpha)
\]
(the corank coincides with the codimension of
the simplex in the geometric realization). The residue of a corank $1$ simplex
is called a \emph{panel}. Abusing notation slightly, we call the link of such a simplex
also a panel.
Given a simplex $\alpha$ and $k\geq \crk(\alpha)$, we denote by $\mathcal E_k(\Delta,\alpha)$ the union
of the links of the corank $k$ faces of $\alpha$,
\[
\mathcal E_k(\Delta,\alpha)=\bigcup\left\{\lk_\Delta(\beta)\mid \beta\subseteq\alpha\text{ and }
\crk(\beta)= k\right\}.
\]
A \emph{gallery} in a pure simplicial complex is a sequence of chambers
$(\gamma_0,\ldots,\gamma_r)$, where $\gamma_{i-1}\cap\gamma_{i}$ 
has corank at most $1$. A gallery \emph{stammers} if $\gamma_{i-1}=\gamma_i$ holds for
some $i$.
A pure simplicial complex
where any two chambers can be connected by some gallery is called a
\emph{chamber complex}. 
A gallery $(\gamma_0,\ldots,\gamma_r)$ is called
\emph{minimal} if there is no gallery from $\gamma_0$ to $\gamma_r$ with
less than $r+1$ chambers.
If every panel contains at least $3$ different chambers, the chamber
complex is called \emph{thick}.
\end{Num}
\begin{Num}\textbf{Geometries\ }
\label{GeometryDefinition}
Suppose that $\Delta$ is a thick chamber complex of rank $n$ with vertex set $V$
and that $I$ is a finite set of $n$ elements. A \emph{type function} is a
map $t:V\rTo I$ whose restriction to every simplex is injective.
We view the type function also as a regular simplicial map $t:\Delta\rTo 2^I$
and extend it to the geometric realizations, $t:|\Delta|\rTo|2^I|$.
The latter map is, for obvious reasons, sometimes called the \emph{accordion map}.
We call $(\Delta,t)$ a \emph{geometry} if $\Delta$ has the following
two properties.

\medskip
(1) $\Delta$ is a flag complex.

(2) The link of every nonmaximal simplex is a chamber complex.

\smallskip\noindent
We remark that what we call here a geometry is called a 
\emph{thick residually connected geometry} in \cite{TitsLocal}.
The \emph{type} (resp. \emph{cotype}) of a simplex $\alpha$ is $t(\alpha)$ 
(resp. $I\setminus t(\alpha)$).
If $\alpha$ is a simplex of cotype $J$, then $\lk(\alpha)$ is a geometry over $J$.
The simplicial join of two geometries is again a geometry.
The \emph{type of a nonstammering gallery} $(\gamma_0,\ldots,\gamma_r)$ is the sequence
$(j_1,\ldots,j_r)\in I^r$, where $j_k$ is the cotype of $\gamma_{k-1}\cap\gamma_k$.
Automorphisms and homomorphisms of geometries are defined in the obvious way; they are
regular simplicial maps which preserve types.

The idea behind this is that the vertices in a geometry are points, lines,
planes and so on. The type function says what kind of geometric object a
given vertex is and the simplices are the flags.
The set of all simplices
of a given type $J\subseteq I$ is the \emph{flag variety} $V_J(\Delta)$.
The chambers are thus the maximal flags, $V_I=\Cham(\Delta)$. 
A gallery shows how one maximal flag can be altered into another maximal flag
by exchanging one vertex at a time.
The type of the gallery records what types of exchanges occur.
\end{Num}
\begin{Num}\textbf{\boldmath Generalized $n$-gons\ }
Let $n\geq 2$ be an integer. A geometry of rank $2$ is a bipartite simplicial 
graph. It is called a \emph{generalized $n$-gon} if it has girth $2n$
and diameter $n$,
i.e. if it contains no circles of length less than $2n$ and if the 
combinatorial distance between two vertices is at most $n$.

A generalized digon is the same as a complete bipartite graph, i.e. the 
simplicial join of two vertex sets (of cardinalities at least $3$, because of
the thickness assumption). A generalized
triangle is the same as an abstract projective plane; one type gives the points
and the other the lines. The axioms above then say that any two distinct
lines intersect in a unique point, and that any two distinct points lie
on a unique line.
\end{Num}
\begin{Lem}
\label{LocalRecognitionOfGeometry}
Let $\Delta$ be a simplicial flag complex with a type function $t:\Delta\rTo2^I$.
Suppose that the link of every vertex $v$ 
is a thick chamber complex of rank $\mathrm{card}(I)-1$.
If $|\Delta|_w$ is connected as a topological space, then $\Delta$ is a geometry.

\proof
The simplicial complex $\Delta$ is pure (since this a local condition).
We have to show that it is gallery-connected. Since the
$1$-skeleton $\Delta^{(1)}$ is connected, it suffices to show
that any two chambers that have a vertex $v$ in common can be joined by
a gallery. But this is true since $\lk(v)$ is a chamber complex.
\qed
\end{Lem}

\begin{Num}\textbf{\boldmath Geometries of type $M$\ }
Suppose that $M:I\times I\rTo\mathbb N$ is a Coxeter matrix,
i.e. $M_{i,j}=M_{j,i}\geq 2$ for all $i\neq j$, and $M_{i,i}=1$ for all $i$.
A geometry $(\Delta,t)$ 
is \emph{of type $M$}  if the link of every simplex $\alpha$ of corank
$2$ and cotype $\{i,j\}$ is a generalized $M_{i,j}$-gon. 

We put
$M_\alpha=M_{i,j}$ for short.
The link of a simplex $\alpha$ of cotype $J$ is a geometry
of type $M'$, where $M'$ is the restriction of $M$ to $J\times J$.
The Coxeter group associated to $M$ is 
\[
W=\bra{I\mid (ij)^{M_{i,j}}=1},
\]
see \cite[5.1]{Hum}.
The Coxeter group and diagram for $M'$ will be called the Coxeter group
and diagram of the simplex $\alpha$.
A gallery is called \emph{reduced} if the word which is represented by
its type in $W$ is reduced in the sense of Coxeter groups, see \cite[5.2]{Hum}.
Recall that a Coxeter group is called \emph{spherical} if it is finite.
We will be mainly concerned with geometries of spherical type.

For the irreducible spherical Coxeter groups we use the standard names 
$\mathsf A_k$, $\mathsf C_k$, $\mathsf D_k$ and so on as in \cite{Hum}. 
By $\mathsf C_3$ and $\mathsf H_3$ we mean in particular the octahedral and the icosahedral
group.
The dihedral group of order $2n$ is denoted $\mathsf I_2(n)=\bra{a,b\mid a^2=b^2=(ab)^n=1}$.
\end{Num}
If $\Delta$ is a geometry of type $M$ whose Coxeter diagram
is not connected, then $\Delta$ is in a natural way a join of
two geometries, see \cite[6.1.3]{TitsLocal}. 
It therefore suffices in many cases to consider
geometries with connected Coxeter diagrams.
A geometry of type $\mathsf A_1$ is a set without further structure. Therefore
a geometry whose Coxeter diagram has an isolated node is a join
of a set with a geometry. For this reason we will often exclude
geometries whose Coxeter diagrams have isolated nodes.
\begin{Lem}
\label{MinimalGalleriesAreReduced}
Suppose that $\Delta$ is a geometry of type $M$. Then every minimal
gallery is reduced. In particular there is a uniform upper bound
on the length of minimal galleries if $M$ is of spherical type.

\proof
This is an easy consequence of the reduction process of words in Coxeter groups,
see 3.4.1--3.4.4 in \cite{TitsLocal}.
\qed
\end{Lem}
\begin{Num}\textbf{Homogeneous geometries\ }
If a group $G$ acts (by type preserving automorphisms) transitively on the chambers
of a geometry $\Delta$, we call the pair $(G,\Delta)$ a
\emph{homogeneous geometry}.
We denote the stabilizer of a simplex $\alpha$ by $G_\alpha$.
If $(G,\Delta)$ is homogeneous, then $(G_\alpha,\lk_\Delta(\alpha))$
is also homogeneous.
The following fact about the bounded generation of stabilizers
will be important on several occasions.
\end{Num}
\begin{Lem}
\label{FiniteStepGeneration}
Let $(G,\Delta)$ be a homogeneous geometry of type $M$.
Let $\gamma$ be a chamber and suppose that $\beta\subseteq \gamma$
is a face of corank at least $1$ whose Coxeter group is of
spherical type. Let $\alpha_1,\ldots,\alpha_t\subseteq\gamma$ be
the faces of corank $1$ which contain $\beta$.
Let $s$ be the length of the longest word in the Coxeter group
of $\beta$. Then the $st$-fold multiplication map
\[
\bigl(G_{\alpha_1}\times\ldots\times G_{\alpha_t}\bigr)^s\rTo G
\]
which sends a sequence of $st$ group elements to their product
has $G_\beta$ as its image.

\proof
It is clear that the image of the multiplication map is
contained in $G_\beta$, since each of the groups
$G_{\alpha_k}$ is contained in $G_\beta$.
Suppose that $g$ is in $G_\beta$. Then there is a
gallery $\gamma=\gamma_0,\gamma_1,\ldots,\gamma_{r-1},\gamma_r=g(\gamma)$
in $\Delta_{\geq\beta}$, with $r\leq s$. 
We show by induction on $r$ that $g$ is in the image of the multiplication map.
For $r=0,1$ we have 
$g\in G_{\alpha_1}\cup\cdots\cup G_{\alpha_t}$. 
For $r>1$ we find $h\in G_{\alpha_1}\cup\cdots\cup G_{\alpha_t}$
with $h(\gamma_0)=\gamma_1$. 
Then $h^{-1}g(\gamma)$ can be connected to $\gamma$ by a gallery
of length $r-1$ in $\Delta_{\geq\beta}$. By the induction hypothesis $h^{-1}g$ can be
written as  a product of $r-1$ elements from $G_{\alpha_1}\cup\cdots\cup G_{\alpha_t}$
and the claim follows.
\qed
\end{Lem}
\begin{Num}\textbf{Simple complexes of groups\ }
\label{SimpleComplexOfGroupsDefinition}
A \emph{simple complex of groups} $\cG$ is a cofunctor from a poset to
the category of group monomorphisms, see \cite[II.12.11]{BH}.
In other words, it assigns in a functorial way to every poset element
$\alpha$ a group $G_\alpha$ and to every inequality $\beta\leq\alpha$ a group
monomorphism $G_\beta\lTo G_\alpha$, such that all resulting triangles of
maps commute. If $H$ is a group, then
a simple homomorphism $\phi:\cG\rTo H$ 
consists of a collection of homomorphisms
$\phi_\alpha:G_\alpha\rTo H$ such that all resulting triangles commute.

Let $\gamma$ be a chamber in a homogeneous geometry $(G,\Delta)$.
The stabilizers $G_\alpha$ of the nonempty simplices
$\emptyset\neq\alpha\subseteq\gamma$ form in a natural way a simple complex
of groups $\cG$, with a simple homomorphism $\cG\rTo G$.
The pair $(G,\Delta)$ is completely determined by this datum $\cG\rTo G$.
We will see that in certain situations $(G,\Delta)$ is already determined by $\cG$.
This is for example true if $|\Delta|_w$ is simply connected as a topological space,
since then $G=\colim\cG$, see \cite{TitsAmalgam} or \cite[II.12.18]{BH}.
However, this condition is not so easy to check in our setting of compact Lie groups
and we will replace the 'abstract' colimit $\colim\cG$ in \ref{ConstructionOfGHat}
by a compact group $\widehat G$ which serves essentially the same purpose in the
category of compact groups. We remark that an analogous construction works for
finite geometries and groups.
\end{Num}
Finally, we need Tits' notion of a $k$-covering of geometries \cite{TitsLocal}.
\begin{Num}\textbf{\boldmath$k$-Coverings\ }
\label{k-CoveringDefinition}
Let $\Delta$ and $\Delta'$ be chamber complexes of rank $n$ and let
$\rho:\Delta\rTo\Delta'$ be a 
surjective regular simplicial map. We call $\rho$ a \emph{$k$-covering}
if for every simplex $\alpha\in\Delta$ of corank at most $k$, the induced
map $\lk_\Delta(\alpha)\rTo \lk_{\Delta'}(\rho(\alpha))$ is an isomorphism.
If $\rho$ is an $(n-1)$-covering, then $|\Delta|_w\rTo|\Delta'|_w$
is a covering in the topological sense. We call an $(n-1)$-covering
a \emph{covering} for short. 
As Tits remarks, one should view $k$-coverings
as 'branched coverings'. 
We note the following: if $\rho:\Delta\rTo\Delta'$ is a covering and 
if $\Delta'$ is a flag complex, then $\Delta$ is also a flag complex (since this is a local
condition). For $k$-coverings between geometries we always assume that they preserve types.
\end{Num}
\begin{Num}\textbf{Universal $k$-coverings\ }
\label{UniversalCoveringDef}
A $k$-covering $\tilde\rho:\widetilde\Delta\rTo\Delta$ is called \emph{universal} if it
has the following property: for every $k$-covering $\rho:\Delta'\rTo\Delta$ and
every pair of chambers $\gamma'\in\Delta'$ and $\tilde\gamma\in\widetilde\Delta$
with $\tilde\rho(\tilde\gamma)=\rho(\gamma')$,
there is a unique $k$-covering $\rho':\widetilde\Delta\rTo\Delta'$ with
$\rho'(\tilde\gamma)=\gamma'$ and $\tilde\rho=\rho\circ\rho'$.
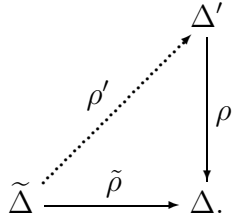
\begin{diagram}[nohug]
&&\Delta'\\
&\ruDotsto^{\rho'}&\dTo_\rho\\
\widetilde\Delta&\rTo^{\tilde\rho}&\Delta .
\end{diagram}
\end{Num}
Applying this universal property twice, we have the following.
\begin{Lem}
\label{LiftingCor}
Let $\rho:\widetilde\Delta\rTo\Delta$ be a universal $k$-covering of geometries of type $M$, 
for $k\geq2$. Suppose $g$ is an automorphism of $\Delta$. Given any two chambers 
$\gamma_1,\gamma_2\in\widetilde\Delta$ with $g(\rho(\gamma_1))=\rho(\gamma_2)$, there
exists a unique automorphism $\tilde g$ of $\widetilde\Delta$ with
$\rho\circ\tilde g=g\circ\rho$ and $\tilde g(\gamma_1)=\gamma_2$.
\qed
\end{Lem}
We call the lifts of the identity \emph{deck transformations}. The following is
an immediate consequence of the previous lemma.
\begin{Prop}
\label{LiftinProp}
Let $\rho:\widetilde\Delta\rTo\Delta$ be a universal $k$-covering of geometries of type $M$, 
for $k\geq2$. Suppose that 
$H\subseteq \Aut(\Delta)$ acts transitively on the chambers of $\Delta$.
Let $\widetilde H\subseteq \Aut(\widetilde\Delta)$ denote the collection of
all lifts of the elements of $H$ and let $F\subseteq\widetilde H$ denote the
collection of all deck transformations. Then we have the following.

(1) $\widetilde H$ is a group acting transitively on the chambers of $\widetilde\Delta$
and $F\subseteq\widetilde H$ is a normal subgroup. 

(2) The map $\rho$ is equivariant with respect to the
map $\widetilde H\rTo H\cong\widetilde H/F$.

(3) If $\alpha\in\widetilde\Delta$ is a simplex of corank at most $k$,
then $\widetilde H_\alpha\cap F=\{\id\}$ and 
$\widetilde H_\alpha$ maps isomorphically onto $H_{\rho(\alpha)}$.

(4) For a simplex $\alpha$ of corank at most $k$ in $\tilde\Delta$, the  
$\widetilde H$-stabilizer of $\rho(\alpha)$ splits as a semidirect product 
$\widetilde H_{\rho(\alpha)}=\widetilde H_\alpha F$.

\proof
From \ref{LiftingCor} we see that products and inverses of lifts are again lifts.
Thus $\widetilde H$ is a group.
The natural map $\widetilde H\rTo H$ which assigns to a lift $\tilde g$
the automorphism $g$ which was lifted is an epimorphism with kernel $F$. Therefore
we have (1) and (2).

Suppose that $\alpha\in\widetilde\Delta$ has corank at most $k$
and that $\gamma\supseteq\alpha$ is a chamber. 
A deck transformation which sends $\gamma$ to a chamber in 
$lk (\alpha)$ must fix $\gamma$ and is therefore the identity.
Thus $F\cap\widetilde H_\alpha=\{\id\}$. 
If $g\in H$ fixes $\rho(\alpha)$,
then we find a unique chamber $\gamma'\in \widetilde\Delta_{\geq\alpha}$
with $\rho(\gamma')=g(\rho(\gamma))$ and hence a lift $\tilde g\in\widetilde H_\alpha$
of $g$. This shows that $\widetilde H_\alpha\rTo H_{\rho(\alpha)}$ is
surjective, and therefore an isomorphism. Thus we have (3).

For (4) we note that $\widetilde H_\alpha F$ fixes $\rho(\alpha)$.
Conversely, suppose that $\tilde g$ in $\widetilde H$ fixes $\rho(\alpha)$
and that $\gamma$ is a chamber containing $\alpha$. Then $\rho(\tilde g(\gamma))$
contains $\rho(\alpha)$. There exists an element $f\in F$ such that
$f(\tilde g(\gamma))\in\widetilde\Delta_\alpha$, because $F$ acts
transitively on the preimage of $g(\rho(\gamma))$. This proves (4).
\qed
\end{Prop}
\begin{Rem}
The existence of a universal $2$-covering of a geometry seems in general to be an open problem.
The existence of a universal $n-1$-covering is not an issue; see also
Pasini \cite[Ch.~12]{Pasini}. For $n\geq 2$, the topological universal
covering $\widetilde{|\Delta|_w}\rTo|\Delta|_w$ is 
the universal $n-1$-covering, as one sees 
from \ref{LocalRecognitionOfGeometry}.
We remark also that an analog of the construction that we give in \ref{ConstructionOfGHat}
below gives universal homogeneous geometries in the class of finite homogeneous geometries
of spherical type. In any case, we have the following important fact.
\end{Rem}
\begin{Thm}[Tits]
\label{BuildingsAre2Connected}
Suppose that $\rho:\widetilde\Delta\rTo\Delta$ is a $2$-covering of geometries.
If $\widetilde\Delta$ is a building, then $\rho$ is universal.

\proof
This follows from Theorem 3 and 2.2 in \cite{TitsLocal}.
\qed
\end{Thm}
We close this section with the following deep result due to Tits. It says that the only
obstruction to the existence of a $2$-covering by a building lies in
the rank $3$ links. We remark that (thick) buildings of
type $\mathsf H_3$ and $\mathsf H_4$ do not exist, see \cite[Addenda]{TitsLNM}.
(Tits' result applies also to non-thick geometries.)
\begin{Thm}[Tits]
\label{TitsCoveringTheorem}
Let $\Delta$ be a geometry of type $M$. Then the following are equivalent.

(1) There exists a building $\widetilde\Delta$ and a $2$-covering
$\widetilde\Delta\rTo\Delta$.

(2) For every simplex $\alpha\in\Delta$ of corank $3$ whose Coxeter diagram
is of type $\mathsf C_3$ or $\mathsf H_3$, there exists a building $\Gamma$ and
a $2$-covering $\Gamma\rTo \lk_\Delta(\alpha)$.

\proof
Our assumptions allow us to go back and forth between chamber systems
and geometries. The result follows thus from 5.3 in \cite{TitsLocal}.
\qed
\end{Thm}
\begin{Num}
The following facts illustrate two interesting cases:

(a) Every geometry $\Delta$
of type $\mathsf A_n$ is a projective geometry and in particular
a building, see \cite[6.1.5]{TitsLocal}. Therefore  
$\id:\Delta\rTo\Delta$ is the universal $2$-covering
(and $\Delta$ admits no quotients).

(b) Suppose that $\Delta$ is a geometry of type $\mathsf C_3$ and that we
call the three types of vertices points, lines and hyperlines
as in \cite[p.~542]{TitsLocal}.
Then $\Delta$ is a building if and only if any two lines which have at
least two distinct points in common are equal, i.e.~if there are no digons,
see \cite[6.2.3]{TitsLocal}.
\end{Num}
\begin{Num}\textbf{The Neumaier Geometry\ }
\label{NeumaierGeometry}
We briefly explain the one known finite geometry of type $\mathsf C_3$ which
is not covered by a building.
Let $V_1$ be a set consisting of seven \emph{points}
and let $V_2=\binom{V_1}{3}$ denote the set of all $3$-element subsets of
$V_1$. These are the \emph{lines} of the geometry.
There are $30$ ways of making $V_1$ into a projective plane
by choosing $7$ appropriate lines in $V_2$; let $X\subseteq \binom{V_2}{7}$
be this set. Finally, let $G=\mathrm{Alt}(V_1)=\mathrm{Alt}(7)$ 
and let $V_3\subseteq X$
be one of the two 15-element $G$-orbits in $X$.
The elements of this orbit are the \emph{planes}
of the geometry. Put $V=V_1\cup V_2\cup V_3$ and define two vertices 
$v,w\in V$ to be adjacent if $v\in w$ or $w\in v$. Let $\Delta$ denote the corresponding
flag complex. Then $\Delta$ is a geometry of type $\mathsf C_3$. We note that
points (vertices of type $1$) and planes (vertices of type $3$)
are always incident. See Neumaier \cite{Neumaier} and Pasini
\cite[6.4.2]{Pasini} for more details.
\end{Num}

\section{\boldmath Compact geometries}

Now we consider actions of compact Lie groups on geometries. This leads to a different
topology on $|\Delta|$. The next definition is very much in the spirit of
Burns-Spatzier \cite{BS}; see also \cite[6.1]{GKVW}. 
Suppose that $\Delta$ is a geometry over $I$. Given a simplex $\alpha$
of type $J\subseteq I$, let $\alpha(j)$ denote its unique vertex of type $j\in J$.
In this way we can view $\alpha$ as a map $J\rTo V$ or as a $J$-tuple of vertices,
$\alpha\in V^J$.
\begin{Def}
\label{DefinitionOfCompactGeometry}
Let $\Delta$ be a geometry of type $M$ over $I$. 
Suppose that the vertex set $V$ of $\Delta$ carries a compact Hausdorff
topology and that for every $J\subseteq I$, the flag variety $V_J$ 
(viewed as a subset of the compact space $V^J$) is closed. Then we call $\Delta$
a \emph{compact geometry}. The proof of \cite[6.6]{GKVW} applies verbatim and
shows that for every simplex $\alpha\in\Delta$, the link
$\lk(\alpha)$ is again a compact geometry.
We say that $\Delta$ has \emph{connected panels} if the panels are
connected in this topology.

Examples of compact geometries arise as follows from groups.
Suppose that $(G,\Delta)$ is a homogeneous geometry of type $M$ and
that $G$ is a locally compact group. If every simplex stabilizer $G_\alpha$
is closed and cocompact (i.e. $G/G_\alpha$ is compact), then
$V$ carries a compact topology and the flag varieties
are also compact, hence closed. We then call $(G,\Delta)$ a
\emph{homogeneous compact geometry}.
The spherical buildings associated to semisimple or reductive
isotropic algebraic groups
over local fields are particular examples of homogeneous compact geometries.

The topology on $V$ can be used to define a new topology on $|\Delta|$ as follows. 
Consider the map
\begin{align*}
p:\Cham(\Delta)\times|2^I|&\rTo|\Delta|\\
\bigl(\gamma,\textstyle\sum_{i\in I} i\cdot\xi(i)\bigr)&\rMapsto 
\textstyle\sum_{i\in I}\gamma(i)\cdot\xi(i)
\end{align*}
Both $\Cham(\Delta)$ and $|2^I|_w$ are compact and we endow $|\Delta|$
with the quotient topology with respect to the map $p$.
The resulting compact space is denoted $|\Delta|_K$. 
The identity map $|\Delta|_w\rTo|\Delta|_K$
is clearly continuous, and we call the topology of $|\Delta|_K$ the
\emph{coarse topology} on $|\Delta|$.
\end{Def}
\begin{Lem}
\label{IsHausdorff}
The space $|\Delta|_K$ is compact Hausdorff. If $\Delta$ has rank at least $2$,
then $|\Delta|_K$ is path-connected.

\proof
From the continuity of the natural maps 
$\Cham(\Delta)\times|2^I|_w\rTo|\Delta|_K\rTo^t|2^I|_w$ we see that
we can separate points which have different $t$-images.

Suppose now that $x,y\in|\Delta|_K$ have the same type
$\xi=t(x)=t(y)\in|2^I|$. 
We let $J=\mathrm{supp}(\xi)=\{j\in I\mid\xi(j)>0\}$ denote the support of 
$\xi$ and we put
\[
u(\xi)=\{\zeta\in|2^I|\mid \mathrm{supp}(\zeta)\supseteq \mathrm{supp}(\xi)\}.
\]
Then $u(\xi)$ is an open neighborhood of $\xi$.
Let $U\subseteq V_J$ be open and let $U_C\subseteq\Cham(\Delta)$
denote the open set of all chambers whose face of type $J$ is in $U$.
We claim that $U_C\times u(\xi)$ is $p$-saturated.
Indeed, if $(\gamma,\zeta)\in U_C\times u(\xi)$ and if
\[
p(\gamma,\zeta)=\textstyle
\sum \gamma(i)\cdot\zeta(i)=\sum\gamma'(i)\cdot\zeta'(i)=p(\gamma',\zeta'),
\]
then $\zeta=\zeta'$ and $t(\gamma\cap\gamma')\supseteq J$.
It follows that the $p$-image of $U_C\times u(\xi)$ is open.

Now for $x,y$ as above, we choose disjoint open neighborhoods
$X,Y\subseteq V_J$ of the type $J$ simplices containing them. 
Then the $p$-images of $X_C\times u(\xi)$ and  $Y_C\times u(\xi)$ are 
disjoint open neighborhoods.

Finally, we note that $|\Delta|_w$ is path-connected if $\Delta$ has rank
at least $2$, so the same is true for $|\Delta|_K$.
\qed
\end{Lem}
\begin{Lem}
\label{FlagVarietiesConnected}
Let $\Delta$ be a compact geometry with connected panels.
Then all flag varieties $V_J$ are connected (in the coarse topology).

\proof
We show first that $\Cham(\Delta)$ is connected.
If $(\gamma_0,\ldots,\gamma_r)$ is a gallery, then
$\gamma_{k-1},\gamma_k$ are in a common panel and hence in
a connected subset. Since $\Delta$ is gallery-connected,
$\Cham(\Delta)$ is connected. For $J\subseteq I$ we have
a continuous surjective map $\Cham(\Delta)\rTo V_J$,
hence $V_J$ is also connected.
\qed
\end{Lem}
\begin{Lem}
\label{JoinDecompositionLemma}
Let $\Delta$ be a geometry of type $M$ over $I$. Suppose that
$\emptyset\subsetneq J\subsetneq I$ and that $M_{j,k}=2$ holds for all 
$j\in J$ and $k\in K=I\setminus J$. Then $\Delta$ is a join of two
geometries $\Delta_1$, $\Delta_2$ of types $M|_{J\times J}$
and $M|_{K\times K}$. If $\Delta$ is a compact geometry, then
this decomposition is compatible with the topology.

\proof
The proof in \cite[6.7]{GKVW} applies verbatim.
\qed
\end{Lem}
A compact homogeneous geometry of type $\mathsf A_1$ is just a compact space
with a transitive group action. Therefore we will often 
assume that the Coxeter diagram of the geometry has no isolated nodes.
A compact homogeneous geometry $(K,\Delta)$ of type $\mathsf A_1\times\mathsf A_1$
consists of two compact spaces $X,Y$ and a transitive $K$-action on $X\times Y$
which is equivariant with respect to the maps $X\lTo^{\pr_1} X\times Y\rTo^{\pr_2} Y$.
Suppose that $X=\SS^m$, that $Y=\SS^n$ and that $K$ is a compact connected Lie group
acting faithfully and transitively on $\SS^m\lTo\SS^m\times\SS^n\rTo\SS^n$. 
We note that such a group $K$ embeds into $\SO(m+1)\times\SO(n)$, see \cite[96.20]{BlauesBuch}.
In this case, a classification is possible. The result that we need is as follows. 
\begin{Lem}
\label{StructureOfD}
Let $K \subseteq \SO(m + 1) \times \SO(n + 1)$ be a compact connected
group acting transitively on $\SS^m\times \SS^n$. Let 
$K_1$ and $K_2$ denote the projections of $K$ to 
$\SO(m + 1)$ and $\SO(n + 1)$ respectively. Assume that $m=1,2,4,8$
and that $K_1$ = $\SO(m+1)$. If $m = 1$ assume in addition that 
$K_2 = \SO(n+1)$ or that $K_2$ is a compact simple Lie group.
Then $K = K_1 \times K_2$, unless 
$m = 2$, $n = 4k - 1$ and $K = \Sp(1) \cdot \Sp(k)$
acting on $\mathrm{Pu}(\mathbb H)\oplus\mathbb H^k$
via $(a,g)\cdot(u,v)=(au\bar a,gv\bar a)$.

\proof
We decompose the Lie algebra $\Lie(K)$ into the ideals
$\fh_1=\Lie(K)\cap(\fso(m+1)\oplus0)$,
$\fh_2=\Lie(K)\cap(0\oplus\fso(n+1))$
and a supplement $\fh_0$, such that
$\Lie(K)=\fh_1\oplus\fh_2\oplus\fh_0$
and $K=(H_1\times H_2)\cdot H_0$, where $H_i$ is the closed
connected normal subgroup with Lie algebra $\fh_i$.
Since $\Lie(K)/\fh_2\cong\fh_1\oplus\fh_0\cong\fso(m+1)$
is either $1$-dimensional or simple,
we have necessarily $\fh_1=0$ or $\fh_0=0$.
We consider these two cases separately.

(a) Assume that $\mathfrak h_0=0$.
Then we have  a product decomposition
of the Lie algebra and therefore $K=H_1\times H_2=K_1\times K_2$.

(b)
Assume that $\fh_1=0$.
Then we have $\fh_0\cong\fso(m+1)$
and thus the Lie algebra of the group induced on the second factor
$\SS^n$ is $\fso(m+1)\oplus\fh_2$. We compare this with the
classification of transitive actions of compact connected Lie groups
on spheres, see \cite[p.~227]{Oni} or \cite[96.20--23]{BlauesBuch}
or \cite[6.1]{Kramer2Trs}. 
We note also that the $K$-stabilizer of a
nonzero vector $v\oplus0\in \RR^{m+1}\oplus \RR^{n+1}$ acts
transitively on $\SS^n$.

If $m=8$ then $n=8,15$ and $\fh_2=0$.
However, a group with Lie algebra $\mathfrak{so}(8)$ cannot act 
transitively on $\SS^8$ or $\SS^{15}$, so this case cannot occur.
If $m=4$, then $n=4,7$ and $\fh_2\subseteq\fso(3)$.
Again, a group with Lie algebra $\fso(4)\oplus\fh_2$ cannot
act transitively on $\SS^4$ or $\SS^7$.
If $m=2$, then $H_0$ cannot act transitively on
$\SS^m\times\SS^n$ by a similar argument as in the case $m=8$. 
Hence $H_0$ is not transitive on $\SS^n$ and the only
remaining possibility is that 
$n=4k-1$ and $\mathfrak h_2=\mathfrak{sp}(k)$.
The case $m=1$, with $\Lie(K)=\RR\oplus\mathfrak h_2$ is excluded
by our assumptions.
\qed
\end{Lem}
A general classification of transitive action on products of spheres can be found in 
Onishchik \cite[p.~274]{Oni}, Straume \cite[Table II]{Straume}. 
The transitive $\SO(4)$-action on $\SS^2\times\SS^3\subseteq\mathrm{Pu}(\HH)\oplus\HH$
will play a role for the exceptional $\mathsf C_3$ geometry.

The following result is a main ingredient in our classification of
homogeneous compact geometries. A spherical building is \emph{Moufang}
if it has a 'large' automorphism group; see \cite[pp.~274]{TitsLNM}. 
The buildings associated to
reductive isotropic algebraic groups have this property and
conversely, the spherical Moufang buildings can be classified in terms
of certain algebraic data \cite{TitsWeiss}. A deep result due to Tits
says that all irreducible spherical buildings of rank at least
$3$ are Moufang, see \cite[4.1.2]{TitsLNM} \cite{WeissSpherical}.
Spherical buildings of rank $2$ need not be Moufang.
\begin{Thm}
\label{BasicStrucureOfRank2Links}
Let $(G,\Delta)$ be a homogeneous compact geometry of type $M$ with connected panels.
Suppose that $\alpha$ is a simplex of corank $2$ and cotype $\{i,j\}$,
and with $M_{i,j}=M_\alpha\geq 3$.
Then $M_\alpha\in\{3,4,6\}$ and $\lk(\alpha)$ is a compact connected Moufang $M_\alpha$-gon
(and explicitly known).

The panels of cotype $i$ and $j$ are homeomorphic to
spheres (in the coarse topology), of dimensions $m_i,m_j\geq 1$. 
If $G$ is compact, then the panel stabilizers act linearly (i.e. as subgroups 
of orthogonal groups) on these panels.
The space $|\lk(\alpha)|_K$ is homeomorphic to a sphere
of dimension $M_{i,j}(m_i+m_j)+1$.

If $M_{i,j}=3$, then $m_i=m_j=1,2,4,8$, 
if $M_{i,j}=6$, then $m_i=m_j=1,2$ and 
if $M_{i,j}=4$, then either $1\in\{m_i,m_j\}$ or $m_i=m_j=2$ or $m_i+m_j$ is odd
(with further number-theoretic restrictions).

In particular, there are no homogeneous compact geometries of type
$\mathsf H_3$ with connected panels.

\proof
The link $\lk(\alpha)$
is a homogeneous compact generalized $M_\alpha$-gon with connected panels.
By the main results in \cite{Szm,KnarrNE,GKK1,GKK2} we have 
$M_\alpha\in\{3,4,6\}$ and $\lk(\alpha)$ it is the compact connected
Moufang $M_\alpha$-gon associated to a simple real Lie group $S$ of $\RR$-rank $2$.
Moreover, $|\lk(\alpha)|_K$ is $G$-equivariantly homeomorphic to unit sphere 
in the tangent space of the symmetric space $S/G$, where $G\subseteq S$ is a maximal
compact subgroup. The principal $G$-orbits have dimension $M_\alpha(m_1+m_2)$ and 
codimension $1$ in this sphere. A different, purely topological proof
that $|\lk(\alpha)|_K$ is a sphere of this dimension is given in \cite{KnarrNE}. 
A complete classification of these compact geometries and their chamber-transitive 
closed connected subgroups is given in \cite{GKK2}.
\qed
\end{Thm}
\begin{Cor}
\label{CommutatorIsTransitive}
Let $(G,\Delta)$ be a homogeneous compact geometry of type $M$ with connected panels.
Suppose also that $G$ is compact.
If the Coxeter diagram of $M$ has no isolated nodes and if all panels have in
(the coarse topology) dimension at least $2$, then the commutator group $[G,G]$
of $G$ acts transitively on the chambers of $\Delta$.

\proof
Let $\alpha\in\Delta$ be a simplex of corank $1$ and let $m$ denote the dimension
of the sphere $|\lk(\alpha)|_K\cong\SS^m$. The stabilizer $G_\alpha$ induces
a transitive subgroup of $\mathrm{O}(m+1)$ on this panel.
Since $m\geq 2$, the commutator group of $G_\alpha$ still acts transitively
on $\lk(\alpha)$, see \cite[p.~94]{Oni}. 
In particular, $([G,G])_\alpha$ acts transitively on
$\lk(\alpha)$. Since any two chambers can be connected by some gallery,
$[G,G]$ acts transitively on the chambers.
\qed
\end{Cor}
\begin{Lem}
\label{EffectiveOnE_1}
Let $(G,\Delta)$ be a homogeneous compact geometry of type $M$ with connected panels.
If $G$ is compact and acts faithfully on $\Delta$, then $G_\gamma$ acts
faithfully on $\mathcal E_1(\Delta,\gamma)$, for every chamber~$\gamma$.

\proof
Suppose that $g\in G_\gamma$ fixes $\mathcal E_1(\Delta,\gamma)$ pointwise. 
Let $\alpha\subseteq\gamma$
be a face of corank $2$. If $M_\alpha\geq 3$,
then $g$ fixes $\lk(\alpha)$ pointwise by \cite[2.2]{GKK2}.
If $M_\alpha=2$, then $\lk(\alpha)$ is a join of two panels and therefore $g$
fixes $\lk(\alpha)$ pointwise. 
Thus $g$ fixes $\mathcal E_2(\Delta,\gamma)$ pointwise. If 
$(\gamma,\gamma')$ is a gallery, then 
$\mathcal E_1(\Delta,\gamma')\subseteq \mathcal E_2(\Delta,\gamma)$,
hence $g$ fixes $\mathcal E_1(\Delta,\gamma')$ pointwise. 
Since $\Delta$ is gallery-connected, we conclude that $g$ acts trivially on $\Delta$.
\qed
\end{Lem}
In order to show that compact groups acting transitively on compact geometries
are automatically Lie groups, we use the following fact.
\begin{Lem}
\label{NSSLemma}
Let $K$ be a compact group and $H\unlhd K$ a closed normal subgroup.
If $H$ and $K/H$ are Lie groups, then $K$ is a Lie group as well.

\proof
We show that $K$ has no small subgroups, see \cite[2.40]{HMCompact}. 
Let $W\subseteq K/H$ be
a neighborhood of the identity which contains no nontrivial subgroup
and let $V$ be its preimage in $K$.
Let $U\subseteq K$ be a neighborhood of the identity such that
$U\cap H$ does not contain a nontrivial subgroup of $H$. Then
$U\cap V$ contains no nontrivial subgroup of $K$.
\qed
\end{Lem}
\begin{Thm}
\label{GIsLieGroup}
Let $(G,\Delta)$ be a homogeneous compact geometry of spherical type $M$,
with connected panels.
Assume that $G$ is compact and acts effectively, and that the Coxeter diagram
of $\Delta$ has no isolated nodes. Then $G$ is a compact Lie group.

\proof
We first show that certain simplex stabilizers are compact Lie groups.
Let $\gamma$ be a chamber. By \ref{EffectiveOnE_1}, $G_\gamma$ acts
faithfully on $\mathcal E_1(\Delta,\gamma)$. From
\ref{BasicStrucureOfRank2Links} we see that $G_\gamma$ injects into a
finite product of orthogonal groups. Thus $G_\gamma$ is a Lie group.
Now let $\alpha\subseteq\gamma$ be a face of corank $1$.
Let $N\unlhd G_\alpha$ denote the kernel of the action of
$G_\alpha$ on the panel $\lk(\alpha)$.
Then $N$ is a closed subgroup of $G_\gamma$ and hence a Lie group.
The quotient $G_\alpha/N$ is by \ref{BasicStrucureOfRank2Links}
a closed subgroup of an orthogonal group and therefore also a Lie group.
By \ref{NSSLemma}, $G_\alpha$ is a Lie group.

Let now $\alpha_1,\ldots,\alpha_t$ denote the corank $1$ faces
of $\gamma$. Let $s$ be the length of the longest word in the Coxeter group
of $\Delta$. Recall from \ref{FiniteStepGeneration} that we have
a surjective multiplication map
$\bigl(G_{\alpha_1}\times\ldots\times G_{\alpha_t}\bigr)^s\rTo G$.
If we compose it with the projection $G\rTo G/\overline{[G,G]}$, it becomes
a surjective continuous homomorphism, since the target group is
abelian. Thus $G/\overline{[G,G]}$ is a compact
abelian Lie group. From the multiplication map we see also that $G$ has only finitely
many path components, and that the path components of $G$ are closed, and
therefore open.
In particular, $G^\circ$ is an open and path-connected subgroup.
Since $\Cham(\Delta)$ is connected by \ref{FlagVarietiesConnected},
the identity component
$G^\circ$ acts transitively and $(G^\circ,\Delta)$ is a homogeneous compact geometry.
It now suffices to show that $G^\circ$ is a compact Lie group,
and for this we may as well assume that $G=G^\circ$ is connected.

Then $G$ is a central quotient $\bigl(Z\times\prod_{\nu\in N}S_\nu\bigr)/D$,
where $(S_\nu)_{\nu\in N}$ is a (possibly infinite) family of simply connected 
compact almost simple Lie groups, $Z$ is a compact
connected abelian group, and $D$ is a compact totally disconnected central
subgroup of the product $Z\times\prod_{\nu\in N}S_\nu$, see \cite[9.24]{HMCompact}.
We claim that the index set $N$ of the product is finite. Otherwise, $G$
admits a homomorphism onto a semisimple Lie group $H$ of dimension 
strictly bigger
than $r=s\dim(G_{\alpha_1}\times\ldots\times G_{\alpha_t})$.
The composite
$\bigl(G_{\alpha_1}\times\ldots\times G_{\alpha_t}\bigr)^s\rTo G\rTo H$
is a smooth map between Lie groups.
Therefore its image has (by Sard's Theorem, see eg.~\cite[Ch.~3]{Mil}) dimension at most
$r$, a contradiction.
So the index set $N$ is finite, and $[G,G]=\overline{[G,G]}$ is a compact semisimple
Lie group. By \ref{NSSLemma}, the group $G$ is a compact Lie group.
\qed
\end{Thm}
The following byproduct of the proof will be useful later.
\begin{Cor}
\label{IdentityComponentIsTransitive}
Under the assumptions of \ref{GIsLieGroup}, the identity component $G^\circ$ acts
transitively on $\Cham(\Delta)$ and $(G^\circ,\Delta)$ is a homogeneous compact
geometry.
\qed
\end{Cor}
We do not know the answer to the following problem (see
Problem 2 and Problem 3 in the introduction).
For compact connected buildings,
it is in both cases affirmative.
\begin{Problem}
Suppose that $M$ has no isolated nodes.
Is the automorphism group of a compact geometry of type $M$ locally compact in the
compact-open topology?
If the geometry is homogeneous, does there necessarily exist a compact 
chamber-transitive group?
\end{Problem}

A first application of \ref{GIsLieGroup} 
is that there is an upper bound for the topological dimension of the chamber set.
\begin{Cor}
\label{DimensionBound}
Let $(G,\Delta)$ be a homogeneous compact geometry of spherical type $M$,
with connected panels. Assume that $G$ is compact and acts effectively,
and that the Coxeter diagram of $\Delta$ has no isolated nodes.
Each panel of cotype $i$ is by \ref{BasicStrucureOfRank2Links} a sphere,
of dimension $m_i$.
Let $(i_1,\ldots,i_r)$ be a representation of the
longest word in the Coxeter group $W$ of $M$ and let $\gamma$ be a chamber. Then
\[
\dim(G)-\dim(G_\gamma)\leq m_{i_1}+m_{i_2}+\cdots+m_{i_r}.
\]

\proof
Let $\alpha\subseteq\gamma$ be a face of corank $1$ and
cotype $i$. The canonical map 
\[
\Cham(\Delta)\cong G/G_\gamma\rTo G/G_\alpha\cong V_{I\setminus\{i\}}
\]
is a locally trivial $G_\alpha/G_\gamma=\SS^{m_i}$-bundle.
We fix a chamber $\gamma=\gamma_0$ and a  
sequence $(i_1,\ldots,i_r)\subseteq I^r$.
Pulling these sphere bundles back several times, we see that
the space of stammering galleries (the 'Bott-Samelson cycles')
\[
\{(\gamma_0,\ldots,\gamma_r)\in \Cham(\Delta)^{r+1}\mid
I\setminus\{i_k\}\subseteq t(\gamma_{k-1}\cap\gamma_k),\,k=1,\ldots r\}
\]
is a smooth manifold of dimension $m_{i_1}+\cdots+m_{i_r}$, see also \cite[7.9]{KramerLoop}.
The map sending such a stammering gallery
$(\gamma_0,\ldots\gamma_r)$ to $\gamma_r$ is smooth,
hence its image has (by Sard's Theorem) dimension at most
$m_{i_1}+\cdots+m_{i_r}$. By \ref{MinimalGalleriesAreReduced},
every chamber can be reached
from $\gamma_0$ by a gallery whose type is reduced.
Since there are only
finitely many reduced words in a spherical Coxeter group
we obtain an upper bound for the dimension. 
From the Bruhat order on the Coxeter group
we see that this upper bound is of the form that we claim,
(and does not depend on the chosen representation of the
longest word), see \cite[5.10]{Hum} and \cite[7.9]{KramerLoop}.
\qed
\end{Cor}
If $\Delta$ is a homogeneous compact building with connected panels
and if the Coxeter diagram is spherical and has no isolated nodes,
then $|\Delta|_K$ is homeomorphic to a sphere.
This is the 'Topological Solomon-Tits-Theorem', which has been proved in various
degrees of generality, see \cite{Mitchell,KnarrNE,KramerDiss}.
For a homogeneous compact geometry, $|\Delta|_K$ need not be a manifold.
However, we have the following result for geometries of rank $3$.
\begin{Prop}
\label{ManifoldProp}
Suppose that $(G,\Delta)$ is a homogeneous compact geometry of irreducible spherical
type $M$ with connected panels and that $G$ is compact. 
If $\Delta$ has rank $3$, then
$|\Delta|_K$ is a closed connected 
topological manifold of dimension $\dim(G/G_\gamma)+2$.

\proof
Replacing $G$ by $G/N$, where $N$ is the kernel of the action,
we may by \ref{GIsLieGroup} assume that $G$ is a compact Lie group
acting faithfully and transitively on the chambers.
We put $I=\{1,2,3\}$ and we fix a chamber $\gamma\in\Delta$.
Recall from \ref{DefinitionOfCompactGeometry} the closed quotient maps 
\begin{diagram}
\Cham(\Delta)\times|2^I|_w & \rTo^p & |\Delta|_K & \rTo^t & |2^I|_w\\
(\gamma',\zeta) &\rMapsto & \textstyle\sum\gamma'(i)\cdot\zeta(i) &\rMapsto&\zeta.
\end{diagram}
We use the same notation as in the proof of \ref{IsHausdorff}.
Suppose that $x\in|\Delta|_K$ has type $t(x)=\xi\in|2^I|$. 
Let $J=\mathrm{supp}(\xi)=\{j\in I\mid \xi(j)\neq0\}$ and 
$u(\xi)=\{\zeta\in|2^I|\mid \mathrm{supp}(\zeta)\supseteq\mathrm{supp}(\xi)\}$.
Let $W=\Cham(\Delta)\times u(\xi)$. This set is $p$-saturated, 
hence its $p$-image is open (compare the proof of \ref{IsHausdorff})
and a neighborhood of $x$. We claim that $p(W)$ is a tube around the orbit
$G(x)\subseteq|\Delta|_K$, see Bredon \cite[II.4]{Bredon}. 
Let $\alpha\subseteq\gamma$ denote the face of type $J$.
We put 
\[\textstyle
r\bigl(\sum\gamma'(i)\cdot\zeta(i)\bigr)=\sum\gamma'(i)\cdot\xi(i).
\]
From the commutative diagram
\begin{diagram}
W & \rTo & G/G_\gamma\\
\dTo^p&&\dTo\\
p(W)&\rTo^r &G(x).
\end{diagram}
we see that $r$ is continuous, since the preimage in $W$ of an open set in
$G(x)\cong G/G_\alpha$ is $p$-saturated (by the same arguments as in the proof
of \ref{IsHausdorff}). Thus $r$ is an equivariant retraction.
From the chamber-transitivity of $G_\alpha$ on $\lk(\alpha)$ we have that
$G_\alpha(p(\{\gamma\}\times u(\xi))=r^{-1}(x)$.
By \cite[II.4.2]{Bredon}, the set $p(W)$ is a tube with slice $S=r^{-1}(x)$ and
\[
p(W)\cong G\times_{G_\alpha}S.
\]
By construction, the slice $S$ is $G_\alpha$-equivariantly 
homeomorphic to a product of $\RR^{\mathrm{card}(J)-1}$ and 
the open cone over $|\lk(\alpha)|_K$.
If $I=J$, then the slice $S$ is thus an open $2$-disk. If $J=\{1,2\}$, then $S$
is an open $m+2$-disk, since $|\lk(\alpha)|_K$ is an $m$-sphere.
If $J=\{1\}$, then $|\lk(\alpha)|_K$ is $G_\alpha$-equivariantly homeomorphic 
to unit sphere in a polar representation of $G_\alpha$ of cohomogeneity $2$
(here we use the classification \ref{BasicStrucureOfRank2Links}).
Thus $p(W)$ is in each case equivariantly homeomorphic to an open disk bundle
over $G/G_\alpha$ and therefore a manifold. 
\qed
\end{Prop}
The previous proof works also for irreducible spherical types of higher rank 
if we assume that all proper links arise from polar representations. 

We collect a few more elementary facts about homogeneous compact geometries.
\begin{Lem}
\label{CoveringInjectiveOnStabilizers}
Suppose that $\rho:(G',\Delta')\rTo(G,\Delta)$ is a continuous equivariant $k$-covering
between homogeneous compact geometries of type $M$ with connected panels,
for $k\geq 2$, and that $G'$ is compact and acts faithfully on $\Delta'$.  
If $\alpha\in\Delta'$ is a simplex of corank at most $k$, then
$G'_\alpha\rTo G_{\rho(\alpha)}$ is injective.

\proof
Let $\gamma$ be a chamber containing $\alpha$. Then
$\mathcal E_1(\Delta',\gamma)\rTo^\rho \mathcal E_1(\Delta,\rho(\gamma))$
is a $G_\gamma'$-equivariant bijection.
By \ref{EffectiveOnE_1}, the group $G_\gamma'$ acts faithfully on $\mathcal E_1(\Delta',\gamma)$,
hence $G_\gamma'\rTo G_{\rho(\gamma)}$ is injective. 
By assumption, $\rho$ maps  $\lk_{\Delta'}(\alpha)$ bijectively onto $\lk_{\Delta}(\rho(\alpha))$.
So if $g\in G_\alpha'$ is in the kernel of $G_\alpha'\rTo G_{\rho(\alpha)}$, then $g\in G_\gamma'$,
and therefore $g=1$.
\qed
\end{Lem}
\begin{Def}
We call a homogeneous compact geometry $(G,\Delta)$ \emph{minimal} if $G$ has no closed
normal chamber-transitive subgroup $N\subseteq G$.
\end{Def}
Such minimal actions are called \emph{irreducible} in Onishchik \cite{Oni}, but this
terminology would conflict with our terminology for geometries.
Since compact Lie groups satisfy the descending chain condition, we have the following fact.
\begin{Lem}
\label{MinimalExists}
Suppose that the Coxeter diagram of $M$ is spherical and has no isolated nodes and
that $(G,\Delta)$ is a homogeneous compact geometry of type $M$ with connected panels.
If $G$ is compact and acts faithfully, 
then there exists a closed connected normal subgroup $K\unlhd G^\circ$
such that $(K,\Delta)$ is minimal. 

\proof
The group $G^\circ$ acts transitively on the chambers by \ref{IdentityComponentIsTransitive}.
Among all closed normal connected chamber transitive subgroups of $G^\circ$,
let $K\subseteq G^\circ$ be a smallest one. Every closed connected normal subgroup of $K$
is also normal in $G^\circ$, hence $(K,\Delta)$ is minimal.
\qed
\end{Lem}
Under the assumptions of \ref{MinimalExists}, the group $K$ is necessarily connected 
(by \ref{IdentityComponentIsTransitive})
and if all panels have dimension at least $2$, then $K$ is semisimple (by \ref{CommutatorIsTransitive}).
\begin{Num}
\label{RecoverWholeGroup}
In the setting of \ref{MinimalExists},
the group $G^\circ$ can be recovered from $K$ as follows. Let $\alpha$ be a simplex and
put $N=\mathrm{Nor}_K(K_\alpha)$. Then $H=N/K_\alpha$ acts from the right on $K/K_\alpha$.
It is not difficult to see that in this action, $H$ is isomorphic to the
centralizer of $K$ in the symmetric group of the set $K/K_\alpha$. Now
let $L\unlhd G^\circ$ be a connected normal complement of $K$, i.e. $G=K\cdot L$ is a
central product with $K\cap L$ finite. The group $L$ is therefore a 
closed connected subgroup of $H$.
See \cite[3.5 and 3.6]{KramerHabil} or Onishchik \cite[p.~75]{Oni} for more details.
Note that this applies to every nonempty simplex $\alpha$. In particular, we have
$K=G^\circ$ if one of the $K$-stabilizers is self-normalizing in $K$.
\end{Num}
\begin{Num}\textbf{\boldmath The category $\mathbf{HCG}(M)$\ }
Our aim is the classification of compact homogeneous geometries
of a given spherical type $M$. To this end, we consider the following
category $\mathbf{HCG}(M)$. Its objects are
homogeneous compact geometries $(G,\Delta)$ of spherical
type $M$ with connected panels, 
where $G$ is a compact group acting transitively and faithfully on $\Cham(\Delta)$.
The morphisms are equivariant $2$-coverings which are continuous with respect
to the coarse topologies on the respective geometries. We note that the
continuity condition can be also phrased as follows: the homomorphisms between the
groups are continuous.

In what follows, we have sometimes to compare 'abstract' homomorphisms
in the sense of \ref{GeometryDefinition} with homomorphisms which are in
addition continuous in the coarse topology. The group of all continuous
automorphisms of a compact geometry $\Delta$ will be denoted $\AutTop(\Delta)$,
in contrast to the group $\Aut(\Delta)$ of all abstract automorphisms
of the underlying combinatorial structure.
\end{Num}
There is a Moufang spherical building $\Delta$ associated to every noncompact 
simple Lie group $S$, which can be defined in various ways. For example, 
there is a Riemannian symmetric space $X=S/K$ of noncompact type
whose connected isometry group is $S$, and whose Tits boundary $\partial_\infty X$
is the (metric) realization $|\Delta|$. The cone topology on $\partial_\infty X$
coincides with the coarse topology on $|\Delta|$. See Eberlein
\cite{Eberlein} and Bridson-Haefliger \cite{BH} for more details. 
The building $\Delta$ can also
be defined in group-theoretic terms: the Lie group $S$ has a canonical
Tits system (or BN-pair) whose building is $\Delta$, see Warner \cite{Warner}.
The latter approach is used in the next result.
\begin{Thm}
\label{AutomorphismStructure}
Let $\Delta$ denote the Moufang building associated to a centerless
simple real Lie group $S$ of real rank $k\geq 2$ and let $K\subseteq S$ be a 
maximal compact subgroup. Then $(K,\Delta)$ is in  $\mathbf{HCG}(M)$,
where $M$ is the relative diagram of $S$, see \cite[Ch.~X Table VI]{He}.

If $S$ is absolutely simple (i.e. if $\Lie(S)\otimes_\RR\CC$ is simple), 
then every automorphism of $\Delta$ is
continuous in the coarse topology, $\Aut(\Delta)=\AutTop(\Delta)$.
Moreover, $\Aut(\Delta)\subseteq \Aut_\RR(\Lie(S))$ is a second
countable Lie group.

If $S$ is a complex Lie group, then $\Aut(\Delta)$ is
a semidirect product of the group $\Aut_\CC(\Delta)$ of $\CC$-algebraic
automorphisms of $\Delta$ and the (uncountable) automorphism group of
the field $\CC$. The group $\AutTop(\Delta)$ is a semidirect
product of $\Aut_\CC(\Delta)$ and $\mathrm{Gal}(\CC/\RR)$ (this is again a 
second countable Lie group). An automorphism of $\Delta$ which is
(in the coarse topology) continuous on at least one panel is continuous
everywhere.

\proof
The simple Lie group $S$ is simple as an abstract group, see
eg. \cite[94.21]{BlauesBuch}. Therefore it coincides with the group $S^\dagger$
generated by the roots groups of $\Delta$.
Thus we have $\Aut(\Delta)\subseteq \Aut(S)$.
From the Iwasawa decomposition $S=K\!AU$ and the fact that the
$S$-stabilizer of a chamber is the Borel subgroup $B=M\!AU$, with $M=\Cen_K(A)$,
we see that $K$ acts transitively on the chambers.

If $S$ is absolutely simple, then its abstract automorphism group $\Aut(S)$
coincides with $\Aut(\Lie(S))$ and is itself a Lie group
by Freudenthal's Continuity Theorem \cite{Freudenthal}, see also
\cite{BorelTits} or \cite{KramerLie}. If $S$ is a complex Lie group,
then its abstract automorphism group
$\Aut(S)$ is a semidirect product of $\Aut_\CC(\Lie(S))$ and 
$\Aut(\CC)$, see \cite{BorelTits} or \cite{KramerLie}.
From the description of the building through the flag varieties of
$S$ it is clear that the group $\Aut(\CC)$ acts indeed on $\Delta$,
and that the action is continuous if and only if the field automorphism
is continuous.
See also Ch.~5 in \cite{TitsLNM} for more details about the automorphism group
of a spherical building over an arbitrary field.

For the last claim, suppose that the abstract automorphism $g$ is
continuous on some panel. Since $S$ acts transitively on the chambers,
we may assume that $g$ fixes a chamber $\gamma$, and that
$g$ is continuous on a panel of cotype $i$ containing $\gamma$. 
Let $i\neq j$ be another 
cotype, and let $\alpha\subseteq\gamma$ be of cotype $\{i,j\}$.
If $M_\alpha=3,4,6$, then there is a continuous bijection between
the two panels which commutes with $g$. 
This is a special property of the complex algebraic generalized
polygons: there exit so-called \emph{projective points}, 
see \cite[2.10]{KramerHolo}. Therefore $g$ is also
continuous on the panel of cotype $j$.
Since the Coxeter diagram of $M$ is connected, we see that
$g$ is continuous on $\mathcal E_1(\Delta,\gamma)$, and hence
everywhere, see \cite[6.16]{GKVW} (or the arguments in \cite[5.1]{BS}).
\qed
\end{Thm}
The previous theorem says in particular that we have a good class of objects in
our category $\mathbf{HCG}(M)$. The next result shows that the continuity of
$2$-coverings is almost automatic if the covering geometry is a building.
In the proof we require the following lemma.
\begin{Lem}
\label{DetectComplexGroups}
Let $S$ be a noncompact simple centerless Lie group. Then $S$ is absolutely simple
if and only if $\Lie(S)$ has a simple rank $1$ Levi factor which is not of type
$\mathfrak{sl}_2\CC$.

\proof
This is Lemma 10 in \cite{KramerLie}. It follows also from the classification
of the real simple Lie algebras and a case-by-case inspection of their root groups, 
see Ch.~X, Table VI in \cite{He}.
\qed
\end{Lem}
We note that $\mathrm{PSL}_2\CC$ is the only connected Lie group acting
$2$-transitively on $\SS^2$, see \cite{Kramer2Trs}. 
Therefore a simple centerless Lie group is complex
if and only if all its root groups are of real dimension $2$.
\begin{Thm}
\label{BuildingCoveringsAreContinuous}
Suppose that $(G,\Delta)$ is a homogeneous compact geometry in $\mathbf{HCG}(M)$
and that the diagram $M$ is spherical and without isolated nodes.
Assume that $\widetilde\Delta$ is a building and that $\rho:\widetilde\Delta\rTo\Delta$
is an abstract $2$-covering. Then $\widetilde\Delta$ is the Moufang building
associated to a semisimple Lie group $S$ of noncompact type. Moreover, there
exists a compact chamber-transitive 
subgroup $K\subseteq S$ and an abstract automorphism $\phi$
of $\widetilde\Delta$ such that $\rho\circ\phi:(K,\widetilde\Delta)\rTo(G,\Delta)$
is a morphism in $\mathbf{HCG}(M)$, i.e. an equivariant continuous $2$-covering.

\proof
Suppose that $\beta\in\widetilde\Delta$ is a simplex of corank $2$, with 
$M_\beta>2$. Then $\lk_{\widetilde\Delta}(\beta)\cong\lk_\Delta(\rho(\beta))$
is by \ref{BasicStrucureOfRank2Links} a Moufang generalized $M_\beta$-gon
associated to a simple noncompact Lie group. Since we excluded factors of type
$\mathsf A_1$, the irreducible factors  of the building $\widetilde\Delta$
are Moufang buildings associated to real simple Lie groups. This holds
because the panels encode the defining field(s) of an
irreducible  spherical Moufang building, see Tits-Weiss \cite[40.22]{TitsWeiss}.

We now fix a chamber $\gamma\in\widetilde\Delta$, with corank $1$ faces
$\alpha_1,\ldots,\alpha_t$. For $i\neq j$ we put $\alpha_{i,j}=\alpha_i\cap\alpha_j$.

\medskip
\emph{Claim: There exists an automorphism $\phi$ of $\widetilde\Delta$ 
fixing $\gamma$ such that
$\mathcal E_2(\widetilde\Delta,\gamma)\rTo^{\rho\circ\phi}\mathcal E_2(\Delta,\rho(\phi(\gamma)))$ is 
a homeomorphism in the coarse topology.}\\
Suppose first that $S$ is absolutely simple. 
If $M_{\alpha_{i,j}}>2$, then 
$\lk_{\widetilde\Delta}(\alpha_{i,j})\rTo^\rho \lk_\Delta(\rho(\alpha_{i,j}))$ is a
homeomorphism by \ref{AutomorphismStructure}. Since the Coxeter diagram is
irreducible and of rank at least $2$, 
$\mathcal E_2(\widetilde\Delta,\gamma)\rTo^\rho\mathcal E_2(\Delta,\rho(\gamma))$ is
a homeomorphism.

Suppose next that $S$ is a complex simple Lie group and that $M_{\alpha_{i,j}}>2$. 
By \ref{AutomorphismStructure} we find
a field automorphism $\phi$ of $\CC$ such that $\phi$ acts on $\widetilde\Delta$, fixes
$\gamma$, and such that 
$\lk_{\widetilde\Delta}(\alpha_{i,j})\rTo^{\rho\circ\phi}\lk_\Delta(\rho(\alpha_{i,j}))$
is a homeomorphism. It follows from \ref{AutomorphismStructure} that
$\lk_{\widetilde\Delta}(\alpha_{i,k})\rTo^{\rho\circ\phi}\lk_\Delta(\rho(\alpha_{i,k}))$
is a homeomorphism whenever $M_{\alpha_{i,k}}>2$.
An easy induction shows now that
$\mathcal E_2(\widetilde\Delta,\gamma)\rTo^\rho\mathcal E_2(\Delta,\rho(\gamma))$ is
a homeomorphism.

Finally, suppose that $S$ has several simple factors. Then the Coxeter diagram of $M$ has
several components and both $\widetilde\Delta$ and $\Delta$ factor as joins. This
factorization is compatible with $\rho$ and we may apply the previous arguments
to the irreducible factors. 
This finishes the proof of the claim. 

\medskip
Replacing $\rho$ by $\rho\circ\phi$, we assume from
now on that $\mathcal E_2(\widetilde\Delta,\gamma)\rTo^{\rho}\mathcal E_2(\Delta,\rho(\gamma))$
is a homeomorphism. We let $K\subseteq \Aut(\widetilde\Delta)$ denote the collection all lifts
of the elements of $G$, see \ref{LiftinProp}.
It remains to prove that $K$ acts continuously and is compact.
To this end we now consider an arbitrary corank~$1$ face
$\alpha=\alpha_i\subseteq\gamma$.

\medskip
\emph{Claim: The stabilizer $K_\alpha$ acts faithfully and continuously
on $\mathcal E_2(\widetilde\Delta,\alpha)$.}\\
We have 
\[
\mathcal E_1(\widetilde\Delta,\gamma)\subseteq 
\mathcal E_2(\widetilde\Delta,\alpha)\subseteq
\mathcal E_2(\widetilde\Delta,\gamma).
\]
Suppose that $g\in K_\alpha$ acts trivially on $\mathcal E_2(\widetilde\Delta,\alpha)$.
Then $g$ fixes $\gamma$ and acts trivially on $\mathcal E_1(\Delta,\rho(\gamma))$,
hence $g$ is a lift of the identity fixing a chamber. By \ref{LiftingCor},
the deck transformation $g$ is the identity.
From the $\rho$-equivariance we see that $K_\alpha$ acts continuously on
$\mathcal E_2(\widetilde\Delta,\alpha)\cong \mathcal E_2(\Delta,\rho(\alpha))$.

\medskip
\emph{Claim: The stabilizer $K_\alpha$ fixes a simplex
$\alpha'$ opposite $\alpha$.}\\
Let $\beta\subseteq\alpha$ be a corank~$2$ simplex. Then $K_\alpha$ acts on
the generalized polygon $\Gamma=\lk_{\widetilde\Delta}(\beta)$. In this action,
it centralizes a Cartan involution of $\Aut(\Gamma)$,
because it acts in the same way as the compact group $G_{\rho(\alpha)}$ on $\Gamma$.
Therefore it fixes a vertex opposite $\alpha\setminus\beta$
in $\Gamma$. Thus $K_\alpha$ fixes a corank $1$ face in $\widetilde\Delta$
having a corank $2$ face in common with $\alpha$. Continuing in this
way, we obtain a geodesic gallery-like sequence of corank $1$ faces fixed
by $K_\alpha$. Eventually, this sequence reaches a corank $1$ face opposite
$\alpha$.

\medskip
\emph{Claim: The group $K_\alpha$ acts continuously on $\widetilde\Delta$ and is compact.}\\
We noticed already that $K_\alpha$ acts continuously on $\mathcal E_2(\widetilde\Delta,\alpha)$.
Since $\mathcal E_1(\widetilde\Delta,\gamma)\subseteq \mathcal E_2(\widetilde\Delta,\alpha)$,
this implies by \cite[6.16]{GKVW} that $K_\alpha$ acts continuously on $\widetilde\Delta$.
We noted above that $K_\alpha$ fixes a simplex $\alpha'$ opposite $\alpha$.
Let $L=\AutTop(\widetilde\Delta)_{\alpha,\alpha'}$ denote the stabilizer of $\alpha,\alpha'$.
The group $L$ acts faithfully on the set $B=\mathcal E_2(\widetilde\Delta,\alpha)$,
and $B$ is compact in the	 coarse topology. The identity map from 
$L$ with the Lie topology to $L$ with the compact-open topology with respect to the 
$L$-action on $B$ is continuous,
and $K_\alpha\subseteq L$ has a compact image in the latter. Thus $K_\alpha\subseteq L$ is
closed in the Lie topology and therefore a second countable Lie group.
It follows from the open mapping theorem that $K_\alpha$ is compact in the Lie topology.

\medskip
The claim of the theorem follows now. 
From \ref{FiniteStepGeneration} we see that $K$ and all stabilizers in $K$ are compact.
Let $s$ denote the length of the longest word in
the Coxeter group of $M$. We have by \ref{FiniteStepGeneration} a commutative diagram
\begin{diagram}
\bigl(K_{\alpha_1}\times\ldots\times K_{\alpha_t}\bigr)^s&&\rTo^{\text{\small closed}}&&K\\
\dTo^{\text{\small homeomorphism}}&&&&\dDotsto\\
\bigl(G_{\rho(\alpha_1)}\times\ldots\times G_{\rho(\alpha_t)}\bigr)^s&&
\rTo^{\text{\small continuous}}&&G.
\end{diagram}
Therefore the dotted homomorphism is continuous.
\qed
\end{Thm}
\begin{Rem}
The proof of \ref{BuildingCoveringsAreContinuous} above relies on
properties of Moufang buildings and Lie groups. 
There is a completely different proof
which constructs the topology on the abstract building $\widetilde\Delta$
from the topology of $\Delta$, without using the group,
see Lytchak \cite{Lytchak} and Fang-Grove-Thorbergsson \cite{FGTh}.
\end{Rem}
Under the assumptions of the previous Theorem \ref{BuildingCoveringsAreContinuous}, 
$G=K/F$ where $F\subseteq K$ is, by \ref{LiftinProp}, a
closed normal subgroup which intersects the stabilizers of corank $k$ 
simplices trivially (where $k\geq 2$ is the largest integer such that $\rho$ is
a $k$-covering). Since we know the possibilities for the compact group $K$
(at least for the irreducible case) from \cite{EH}, a great deal can be said about
the possibilities for $F$. We indicate for a few examples how such a classification
works.
\begin{Prop}
\label{AutomaticallyBuilding}
Assume that $(G,\Delta)$ is a homogeneous compact geometry in $\mathbf{HCG}(M)$
and that the Coxeter diagram of $M$ is irreducible. In the following 
three situations,
$\Delta$ is necessarily the building associated to a simple Lie group $S$, and
$G^\circ$ is a maximal compact subgroup of $S$.

(1) The diagram $M$ is of type $\mathsf A_n$.

(2) All panels have dimension $2$.

(3) The diagram $M$ is of type $\mathsf E_6$.

\proof
A geometry of type $\mathsf A_n$ is always a building by \cite[6.1.5]{TitsLocal}.
By the previous theorem, $\Delta$ is the building associated to a simple
Lie group $S$ (for $n\geq 3$ this is due to  Kolmogoroff \cite{Kolmogoroff}).
From \cite{EH} we see that $G^\circ$ is a maximal compact 
subgroup of $S$. Thus we have the result (1).

Assume now that all panels have dimension $2$. By \ref{C3ClassificationTheorem} 
below, a  $\mathsf C_3$ geometry with $2$-dimensional panels is $2$-covered by a
building.
From \ref{TitsCoveringTheorem} we see that $\Delta$ is $2$-covered by
a building $\widetilde\Delta$. By \ref{BuildingCoveringsAreContinuous}, 
the building corresponds to a simple centerless Lie group $S$. By 
the remark following \ref{DetectComplexGroups}, the Lie group $S$ is complex.
Thus a maximal compact subgroup $K\subseteq S$ is centerless simple.
By \cite{EH}, $K$ has no chamber-transitive proper closed subgroups,
hence $K=G^\circ$.

For (3) we note that all panels are either $1$- or $2$-dimensional, and
the $2$-dimensional case is covered by (2). In the $1$-dimensional case,
we have $G^\circ=\mathrm{PSp(4)}$ by \cite{EH}, and this group is simple.
\qed
\end{Prop}
For the buildings of type  $\mathsf E_7$ and $\mathsf E_8$ with $1$-dimensional panels,
the Lie algebra of $G^\circ$ is simple, but $G^\circ$ has nontrivial finite center.

In order to proceed with the classification of homogeneous compact geometries,
we need a substitute for the building, i.e. a good universal object in
the class of homogeneous compact geometries. The remainder of this section
will be devoted to the construction of this compact universal geometry.
\begin{Num}\textbf{\boldmath Simple complexes of groups in $\mathbf{HCG}(M)$\ }
Let $M$ be a Coxeter matrix of spherical type over the index set $I$.
Let $\cG$ be a simple complex of compact groups and continuous homomorphisms,
indexed by the poset of nonempty subsets of $I$, i.e. 
$\cG=\{G_J\mid \emptyset\neq J\subseteq I\}$.

We consider the following category $\mathbf{HCG}_\cG(M)$.
Its objects are quadruples $(G,\Delta,\gamma,\psi)$,
where $(G,\Delta)$ is a geometry in $\mathbf{HCG}(M)$ and
$\gamma$ is a chamber of $\Delta$, and $\psi$ is an isomorphism between
$\cG$ and the simple complex of groups 
$\{G_\alpha\mid \emptyset\neq\alpha\subseteq\gamma\}$.
We assume that for each
group $G_J\in\cG$ we have $\psi(G_J)=G_\alpha$, where
$\alpha$ is the unique face of type $J$ of $\gamma$.

A morphism in $\mathbf{HCG}_\cG(M)$ is an equivariant morphism between the
geometries in $\mathbf{HCG}(M)$ which preserves the preferred chambers and
which commutes with the isomorphisms between $\cG$ and the stabilizer complex.
We remark that such a morphism is unique.
\end{Num}
Our aim is to show that there is a universal object in this category. The main
ingredient is the following construction.
\begin{Num}\textbf{The basic coset construction\ }
\label{TheBasicCosetConstruction}
Let $(G,\Delta)$ be a homogeneous compact geometry in $\mathbf{HCG}(M)$.
Let $\gamma\in\Delta$ be a chamber and let $\cG$ denote the
simple
complex of groups formed by the stabilizers $G_\alpha$, 
for $\emptyset\neq\alpha\subseteq\gamma$, see \ref{SimpleComplexOfGroupsDefinition}.
Suppose that $H$ is a topological group and that $\psi:\cG\rTo H$ is
a continuous simple homomorphism (i.e. that each homomorphism 
$\psi:G_\alpha\rTo H$ is continuous and that all triangles commute). 
In this situation we construct a new homogeneous compact geometry
$(G',\Delta')$ in $\mathbf{HCG}(M)$ and a covering 
\[
\rho:(G',\Delta')\rTo(G,\Delta)
\]
as follows.

For $g\in G_\alpha$ put $\psi'(g)=(\psi(g),g)\in H\times G$. 
This defines a continuous and injective simple homomorphism $\psi':\cG\rTo H\times G$. 
We put $G'_{\alpha'}=\psi'(G_\alpha)\subseteq H\times G$ and we let
$G'\subseteq H\times G$ denote the group which is 
algebraically generated by the $G'_{\alpha'}$.
In order to construct $\Delta'$, we use the following standard method,
see Tits \cite[1.4]{TitsLNM} and Bridson-Haefliger \cite[II.12.18--22]{BH}.

Let $v_1,\ldots,v_t$ denote the vertices of the chamber $\gamma$.
The set of cosets $G'/G'_{v_1'}\cup\cdots\cup G'/G'_{v_t'}$
covers $G'$. We let $\Delta'$ denote the nerve of this cover. 
It is easy to see that the simplices of $\Delta'$ correspond
bijectively to the cosets $gG'_{\alpha'}$, for $\emptyset\neq\alpha\subseteq\gamma$ and $g\in G'$.
The inclusion of simplices corresponds to the reversed inclusion of cosets. 
In particular we see that $\Delta'$ is a pure simplicial complex.
The residue $\Delta'_{\geq G'_{\alpha'}}$ of the simplex 
$G'_{\alpha'}$ consists of all cosets $gG'_{\beta'}$ with $g\in G'_{\alpha'}$ and
$\beta\supseteq\alpha$.
Moreover, there is a well-defined type function on $\Delta'$ which maps
$gG'_{v_i'}$ to the type $t(v_i)$. We note also that the projection
$pr_2:H\times G\rTo G$ induces a continuous surjective homomorphism
$p:G'\rTo G$, and a regular simplicial map $p:\Delta'\rTo\Delta$
which maps $gG'_{\alpha'}$ to $pr_2(g)(\alpha)$.

\medskip
{\em Claim: $\Delta'$ is a thick chamber complex.}\\
Every element $g\in G'$ can be written as a product $g=g_1\cdots g_r$,
where $g_k$ is in the stabilizer of a corank~$1$ face of $\gamma$.
This gives a gallery from $G'_{\gamma'}$ to $gG'_{\gamma'}$.
The panels of $\Delta'$ have the same cardinalities as the panels of $\Delta$,
hence $\Delta'$ is thick.

\medskip
{\em Claim: $\Delta'$ is a geometry over $I$, of the same type $M$, and 
$p:\Delta'\rTo\Delta$ is a covering.}\\
From the description of the residues above we see that the
link of a nonempty simplex in $\Delta'$ maps
isomorphically onto a link in $\Delta$. Thus $p$ is a covering.
In particular, $\Delta'$ is a flag complex.

\medskip
{\em Claim: $(G',\Delta')$ is a homogeneous compact geometry. The group
$G'$ is compact and acts faithfully.}\\
Obviously, $\Delta'$ is a homogeneous geometry. The groups $G'_{\alpha'}$
are by construction compact. From the bounded generation \ref{FiniteStepGeneration} 
we see that $G'$ is also compact. Suppose that $(h,g)\in G'_\gamma$ acts trivially on
$\Delta'$. Then $g=\id_\Delta$. Since $(h,g)\in G'_{\gamma'}$, we have $h=1$.

\medskip
We record a few more useful facts about $(G',\Delta')$.

\medskip
{\em Fact: The subgroup of $H$ generated by the $\psi(G_\alpha)$ is compact.}\\
This groups is the image of the compact group $G'$ under $pr_1:H\times G\rTo H$.

\medskip
{\em Fact: Let $F\subseteq G'$ denote the kernel of $G'\rTo^p G$. Then
$F$ intersects every simplex stabilizer trivially, i.e. $F$ acts freely on $\Delta'$.
The $G'$-stabilizer of $\alpha\in\Delta$ is a semidirect product $G'_\alpha=G'_{\alpha'}F$.}\\
Consider an element $(h,id_\Delta)\in F\cap G'_{\alpha'}$. Since $G_\alpha\rTo^{\psi'} G'_{\alpha'}$
is bijective, we have $h=1$. Suppose now that the group element
$(h,g)\in G'$ fixes the simplex $\alpha$.
Then we have $g\in G_\alpha$. Let $h_1=\psi(g)$. Then we have 
$\psi'(g^{-1})=(h_1^{-1},g^{-1})\in G'_{\alpha'}$ and 
$(h,g)(h_1^{-1},g^{-1})=(hh_1^{-1},id_\Delta)\in F$.
\end{Num}
We now use the Basic Coset Construction \ref{TheBasicCosetConstruction} 
in order to construct a universal object in  $\mathbf{HCG}_\cG(M)$.
\begin{Thm}
\label{ConstructionOfGHat}
Suppose that $M$ is spherical without isolated nodes over the index set $I$, 
that $\cG$ is a simple complex of compact groups and continuous homomorphisms over 
the collection of the nonempty subsets $J\subseteq I$ and that
$\mathbf{HCG}_\cG(M)$ is not empty.
Then there exists a homogeneous compact geometry $(\widehat G,\widehat\Delta,\hat\gamma,\hat\psi)$
in $\mathbf{HCG}_\cG(M)$ which has a unique morphism $\rho$ to every $(G,\Delta,\gamma,\psi)$ in
$\mathbf{HCG}_\cG(M)$. 

\proof
We choose a 'transversal' in $\mathbf{HCG}_\cG(M)$, i.e. a family 
$(G_\nu,\Delta_\nu,\gamma_\nu,\psi_\nu)_{\nu\in N}$
of objects in $\mathbf{HCG}_\cG(M)$ which contains one member of each isomorphism class.
Such a family exists since there are only countably many isomorphism classes
of compact Lie groups (every compact Lie group can be realized as an
algebraic matrix group). The cardinality of the index set $N$ is not important here;
we need only the fact that such a set exists.
The $\psi_\nu$ fit together to a continuous simple
homomorphism $\psi:\cG\rTo \prod_{\nu\in N}G_\nu$. 
Let $\widehat G\subseteq\prod_{\nu\in N}G_\nu$ denote the group generated algebraically by the
groups $\psi(G_J)$. 
The basic coset construction \ref{TheBasicCosetConstruction}
gives us a homogeneous compact geometry $(\widehat G,\widehat\Delta)$
and for each $\nu$ a continuous equivariant covering 
$\rho_\nu:(\widehat G,\widehat\Delta)\rTo(G_\nu,\Delta_\nu)$.
This morphism is unique, as we remarked above.
\qed
\end{Thm}
\begin{Def}
\label{UniversalHomogeneous}
We call the pair $(\widehat G,\widehat\Delta)$ constructed
in \ref{ConstructionOfGHat} a
\emph{universal homogeneous compact geometry} for the pair $(\cG,M)$
(obviously this homogeneous compact geometry is unique up to isomorphism).
If an element of the class $\mathbf{HCG}_\cG(M)$ can be covered by a 
building $\widetilde\Delta$, 
then this building is the universal homogeneous compact geometry by Tits' result
\ref{BuildingsAre2Connected} and by \ref{BuildingCoveringsAreContinuous}.
\end{Def}
The group $\widehat G$ has the following universal property.
\begin{Prop}
\label{UniversalProperty}
Suppose that $(\widehat G,\widehat\Delta)$ is a
universal homogeneous compact geometry for the pair $(\cG,M)$.
Suppose that $H$ is a topological group and that $\phi:\cG\rTo H$ is a
simple continuous homomorphism. Then there is a unique
continuous homomorphism $\psi:\widehat G\rTo H$ such that
the diagram 
\begin{diagram}[nohug]
\cG & \rTo & \widehat G\\
& \rdTo^\phi &\dTo_\psi \\
&& H
\end{diagram}
commutes.

\proof
The Basic Coset Construction \ref{TheBasicCosetConstruction} applied to 
$\cG\rTo H\times \widehat G$
gives us a geometry $(G',\Delta')$ and a map $G'\rTo H$. From the universal
property of $(\widehat G,\widehat\Delta)$ we have a homomorphism
$\widehat G\rTo G'\rTo H$. The uniqueness is clear.
\qed
\end{Prop}
Finally, we note that we can pass to a minimal universal homogeneous compact geometry.
\begin{Prop}
\label{UniversalAndMinimal}
Suppose that $M$ is of spherical type and without isolated nodes.
Suppose that $\cG$ is a simple complex of compact groups and that
$\mathbf{HCG}_\cG(M)$ is nonempty. Then there exists a simple
complex of compact groups $\mathcal K$ formed by subgroups 
$K_J\subseteq G_J$ such that $\mathbf{HCG}_{\mathcal K}(M)$ is
nonempty, with the following properties.

(1) The universal homogeneous compact geometry $(\widehat K,\widehat\Delta)$
in  $\mathbf{HCG}_{\mathcal K}(M)$ is minimal.

(2) For every geometry $(G,\Delta)$ in $\mathbf{HCG}_\cG(M)$
there is an equivariant covering
\[
(\widehat K,\widehat\Delta)\rTo(G,\Delta).
\]

\proof
We construct a sequence of equivariant coverings 
\[
\cdots\rTo(G_{k+1},\Delta_{k+1})\rTo(G_k,\Delta_k)\rTo\cdots
\rTo (G_0,\Delta_0)
\]
and simple complexes of compact groups $\cG_k$
as follows. Let $\cG_0=\cG$ and let $(G_0,\Delta_0)$ denote
the corresponding universal homogeneous compact geometry in $\mathbf{HCG}_{\cG_0}(M)$.
Given $\cG_k$ and $(G_k,\Delta_k)$, we choose a closed chamber transitive 
subgroup $H\subseteq G_k$ such that $(H,\Delta_k)$ is minimal.
Let $\cG_{k+1}$ denote the simple complex of groups formed by the
simplex stabilizers of $H$, and let $(G_{k+1},\Delta_{k+1})$
denote the corresponding universal homogeneous compact geometry in
$\mathbf{HCG}_{\cG_{k+1}}(M)$.
If $\cG_{k+1}\neq\cG_k$, then the stabilizers have become strictly smaller.
Since there are no infinite descending sequences of
closed compact Lie groups, this process becomes stationary in finite time $k$, and
we may put $\mathcal K=\cG_k$.
\qed
\end{Prop}
We remark that a completely analogous construction works for the class of
finite homogeneous geometries.

\section{\boldmath Homogeneous compact geometries of type~$\mathsf C_3$}
\label{Section3}
In this section we review the known examples of universal homogeneous compact geometries
$(G,\Delta)$ of type $\mathsf C_3$. In Section~\ref{Section4} we will show that this
list of examples is complete: such a geometry is either a building (a polar space of
rank $3$), or the exceptional geometry discovered by Podest\`a-Thorbergsson
\cite{PTh}, which we describe in Section~\ref{TheExceptionalC3Geometry}. 
Some of the results in the present section will be used in the classification.
We begin with the classical geometries, the buildings of type $\mathsf C_3$.

\subsection{Projective and polar spaces and their Veronese representations}
Almost all buildings of type $\mathsf C_3$ arise from hermitian forms. 
We review the relevant linear algebra, since we will use it in our
classification in Section~\ref{Section4}. Buildings of type
$\mathsf C_n$ 
are also called \emph{polar spaces} of rank $n$.
Buildings of type $\mathsf C_2$ are called \emph{generalized quadrangles}.
\begin{Num}\textbf{Polar spaces\ }
\label{PolarSpaceDefinition}
Let $\FF=\RR,\CC,\HH$ and let $\sigma$ be an involution on $\FF$, i.e. an
additive map with $a^{\sigma^2}=a$ and $(ab)^\sigma=b^\sigma a^\sigma$, for all $a,b\in \FF$. 
The involution $\sigma$ extends in a natural way to matrices, acting by matrix transposition
combined with entry-wise application of $\sigma$.
Let $V$ be a finite dimensional right $\FF$-module and let $\eps=\pm1$. 
A nondegenerate \emph{$(\eps,\sigma)$-hermitian form} is a biadditive map
$h:V\times V\rTo\FF$ with the properties
\[
h(v,w)=\eps{h(w,v)^\sigma}\qquad h(va,wb)=a^\sigma h(v,w)b\qquad h(V,w)=0\ \Rightarrow\ w=0.
\]
The relevant examples are 

\smallskip
\emph{symmetric bilinear forms}  with $(\eps,\sigma)=(1,\id)$ and $\FF=\RR,\CC$,

\smallskip
\emph{symplectic forms} with $(\eps,\sigma)=(-1,\id)$ and $\FF=\RR,\CC$, and

\smallskip
\emph{$(\eps,\sigma)$-hermitian forms} with $a^\sigma=\bar a$ and $\FF=\CC,\HH$. 

\medskip\noindent
A nonzero subspace $W\subseteq V$ is called \emph{totally isotropic} if $W\subseteq W^{\perp_h}$.
The form $h$ is called \emph{isotropic} if there exist totally isotropic subspaces. The maximal
dimension $k$ of a totally isotropic subspace is the \emph{Witt index} of $h$.
The corresponding
geometry $\Delta$ has as its vertices the collection of all isotropic subspaces.
The simplices in $\Delta$ are the ascending chains of isotropic subspaces.
This simplicial complex is a building of type $\mathsf C_k$,
unless $\dim(V)=2k$ and $(\eps,\sigma)=(1,\id)$.
In the latter case, a slightly modified simplicial complex is a building of type $\mathsf D_{k}$,
see \cite[7.12]{TitsLNM}. We refer to \cite[Ch.~7,~8]{TitsLNM} for more details.

The automorphism group of this building is (an extension of) the projective unitary group
of the form $h$. Its identity component $S$ is a noncompact simple Lie group
of classical type and $(S,\Delta)$ is a homogeneous compact geometry.
If $G\subseteq S$ is a maximal compact subgroup, then also $(G,\Delta)$ is a compact
homogeneous geometry.
\end{Num}
For each of these polar spaces mentioned above, it is possible to describe the 
associated polar representation in terms of certain tensors and geometric algebra.
We first recall the definition of a polar representation.
\begin{Num}\textbf{Polar representations\ }
\label{DefinitionOfPolarRepresentation}
An orthogonal representation of a compact Lie group is called \emph{polar}
if there exists a linear subspace that meets every orbit orthogonally.
Polar representations were classified up to orbit equivalence by
Hsiang-Lawson \cite{HsiangLawson} and Dadok \cite{Dadok}. 
See Eschenburg-Heintze \cite{EH1,EH} for a modern account.
The result is that every polar representation is orbit equivalent to an
\emph{$s$-representation}. An $s$-representation is defined as follows.
Let $S$ be a semisimple centerless Lie group of noncompact type,
let $G\subseteq S$ be a maximal compact subgroup and let $\Lie(S)=\Lie(G)\oplus\mathfrak P$
be the corresponding Cartan decomposition. The adjoint representation of $G$ on $\mathfrak P$
is the associated $s$-representation. It is polar, and if $\Delta$ is the associated
building, then $|\Delta|_K$ is $G$-equivariantly homeomorphic to the unit sphere 
$\SS(\mathfrak P)\subseteq\mathfrak P$.
\end{Num}
\begin{Num}\textbf{Polar representations for certain polar spaces\ }
\label{Veronesean}
Suppose that $\eps=1$ and $a^\sigma=\bar a$, for $\FF=\RR,\CC,\HH$.
Let $f_k$ denote the standard positive definite hermitian form on 
$\FF^k$, i.e. 
\[
f_k(v,w)=\sum_{j=1}^k \bar v_jw_j.
\]
Let $\mathrm U_k\FF$ denote the corresponding unitary group,
\[
\mathrm U_k\FF=\{g\in\FF^{k\times k}\mid g^\sigma g=1\}
\]
(recall that $(x_{i,j})^\sigma=(\bar x_{j,i})$).
Consider the hermitian form 
\[
h=(-f_k)\oplus f_\ell
\]
on $V=\FF^{k+\ell}$,  with $k\leq\ell$ (resp. $k<\ell$ for $\FF=\RR$). 
The Witt index of $h$ is $k$ and
$\mathrm U(k)\times\mathrm  U(\ell)$ is a maximal compact subgroup of the unitary group 
\[
\mathrm U(h)=\mathrm U_{k,\ell}\FF=\{g\in \mathrm{GL}(V)\mid h(-,-)=h(g(-),g(-))\}.
\]
We identify the tensor product $\FF^k\otimes_\FF(\FF^{\ell})^\sigma$
with the $\RR$-module 
$\FF^{k\times\ell}$ and we note that $\mathrm U_k\FF\times\mathrm U_\ell\FF$
acts in a natural way on $\FF^{k\times\ell}$, via
\[
(g_1,g_2,X)\mapstoo g_1Xg_2^\sigma.
\]
This is the polar representation we are interested in. 

There are natural projections $\FF^k\lTo^{\pr_1}\FF^{k+\ell}\rTo^{\pr_2}\FF^\ell$.
For every $t$-dimensional totally isotropic subspace $W\subseteq \FF^{k+\ell}$ there exists a basis
$w_1,\ldots,w_t$ such that $\{u_1=\pr_1(w_1),\ldots,u_t=\pr_1(w_t)\}\subseteq\FF^k$ and
$\{v_1=\pr_2(w_1),\ldots,v_t=\pr_2(w_t)\}\subseteq\FF^\ell$ are orthonormal. The map which sends
the subspace $W$ to 
$\frac{1}{\sqrt{t}}\bigl(u_1\otimes v_1^\sigma+\cdots+u_t\otimes v_t^\sigma\bigr)\in \FF^k\otimes_\FF(\FF^{\ell})^\sigma$
is well-defined and $\mathrm U_k\FF\times \mathrm U_\ell\FF$-equivariant. 
This map extends to a mapping
\[
|\Delta|\rTo\SS(\FF^k\otimes_\FF(\FF^{\ell})^\sigma)=\SS^{k\cdot\ell\cdot\dim_\RR\FF-1}
\]
which is a homeomorphism in the coarse topology.
This map is called the \emph{Veronese representation} of $\Delta$.
The Veronese representation of $\Delta$ lends itself to computations of vertex stabilizers
in $\mathrm U_k\FF\times \mathrm U_\ell\FF$. 

Finally, we note that for $\FF=\RR,\CC$ and $k<\ell$ the groups 
$\SO(k)\times\SO(\ell)$ and $\SU(k)\times\SU(\ell)$
act transitively on the chambers. According to Eschenburg-Heintze \cite{EH} 
and \cite{GKK2}  these are the smallest
compact chamber-transitive groups $K$, unless we are in one of the following exceptional cases.

$(k,\ell)=(2,7)$ and $K=\SO(2)\cdot\mathrm{G_2}$

$(k,\ell)=(2,8)$ and $K=\SO(2)\cdot\Spin(7)$

$(k,\ell)=(3,8)$ and $K=\SO(3)\cdot\Spin(7)$
\end{Num}
We remark that for the other types of $(\eps,\sigma)$-hermitian forms, similar models
for Veronese representations can be worked out in terms of tensor products and
exterior products. These will not be needed here. However, we need also polar
representations for the classical projective geometries over $\RR$, $\CC$, $\HH$ and
the Cayley algebra $\mathbb{O}$.
\begin{Num}\textbf{Polar representations for projective geometries\ }
\label{PolarRepresentationOfProjectiveGeometry}
Let $\FF=\RR,\CC,\HH$ and $V=\FF^{n+1}$, endowed 
with the standard hermitian form $f_{n+1}$. 
The projective geometry over $\FF$ (of type $\mathsf A_n$)
is the simplicial complex $\Delta$
whose vertices are the proper nonzero subspaces of $V$. The simplices
are the partial flags. The noncompact Lie group 
$\mathrm{GL}(V)=\mathrm{GL}_{n+1}\FF$ acts transitively on
the chambers of $\Delta$.
A maximal compact subgroup is the unitary group $\mathrm U_{n+1}\FF$.
Suppose that $W\subseteq V$ is a $t$-dimensional subspace, with an orthonormal basis
$w_1,\ldots,w_t$. The map which sends $W$ to the traceless hermitian matrix 
$\bigl(w_1\otimes w_1^\sigma+\cdots +w_t\otimes w_t^\sigma\bigr)-t\boldsymbol1$
is well-defined and $\mathrm U_{n+1}\FF$-equivariant. It extends to a map
\[
|\Delta|\rTo\SS^\ell,
\] 
where $\ell=\frac{n(n+1)}{2}\dim_\RR(\FF)+n-1$.
One can check that this map induces a homeomorphism $|\Delta|_K\rTo \SS^\ell$.
This is the \emph{Veronese representation} of $\Delta$. Again, the computation of vertex stabilizers
in $\mathrm U_{n+1}\FF$ is easily done in this representation. We note that for $n=2$, the minimal
transitive faithful compact groups are $\SO(3)$, $\mathrm{PSU}(3)$ and $\mathrm{PSp}(3)$.

The projective Cayley plane has no simple description in terms of $\mathbb O^3$, since this
is not a module over the Cayley algebra $\mathbb O$. Nevertheless, the right-hand
side of the Veronese representation in terms of traceless hermitian $3\times 3$-matrices
over $\mathbb O$ makes sense and leads to a Veronese representation of this geometry,
The compact group in question is the centerless simple compact Lie group
$\mathrm F_4$, with vertex stabilizers 
$\Spin(9)$ and chamber stabilizer $\Spin(8)$. We refer to Freudenthal \cite{FreudenthalOktaven},
Salzmann et al.~\cite[Ch.~1]{BlauesBuch}, and to Section~\ref{TheExceptionalC3Geometry} below.
\end{Num}
\begin{Num}\textbf{The nonembeddable polar space\ }
\label{TheNonEmbeddablePolarSpace}
Over the reals, there is one polar space $\Delta$ of type $\mathsf C_3$ which is not associated to a
hermitian form, see \cite[Ch.~9]{TitsLNM}. Instead, it is related to the Cayley algebra.
The corresponding simple noncompact
Lie group is of type $\mathrm E_{7(-25)}$ and its maximal compact subgroup is
$G=\mathrm E_6\cdot\SO(2)$. (In Cartan's classification, this is the noncompact
case E\,VII, see Helgason \cite[Ch.~X, Table V]{He}.)
Its Veronese representation $|\Delta|\rTo\SS^{53}$
arises from the corresponding $s$-representation as in \ref{DefinitionOfPolarRepresentation}.
The panels have dimensions $8$ and $1$ and the links of the vertices are
projective 
Cayley planes and generalized quadrangles belonging to the symmetric bilinear form
$h=(-f_2)\oplus f_{10}$ on $\RR^{2+10}$.
Apparently, no 'simple model' for $\Delta$ and its Veronese 
representation is known; the abstract
construction is purely Lie-theoretic.
\end{Num}

\subsection{\boldmath The exceptional $\mathsf C_3$ geometry}
\label{TheExceptionalC3Geometry}
We now construct the exceptional geometry of type $\mathsf C_3$ that
was found by Podest\`a-Thorbergsson \cite[2B.3]{PTh}. We use the Veronese
representation of the Cayley plane as a focal manifold of an isoparametric
foliation in $\SS^{25}$, corresponding to the $s$-representation of 
the symmetric space for $(\mathfrak e_{6(-26)},\mathfrak f_4)$.
For the description of the Cayley plane which we use, see also 
Cartan \cite{CartanIsop}, Console-Olmos \cite{ConsoleOlmos},
Freudenthal \cite{FreudenthalOktaven}, Karcher \cite{Karcher}, Knarr-Kramer \cite{KnarrKramer}
and Salzmann et al. \cite[Ch.~1]{BlauesBuch}. 

\begin{Num}\textbf{\boldmath The action of $\SU(3)$ on $\OO$\ }
\label{SU(3)Action}
We first recall some algebraic facts.
The real Cayley division algebra $\OO$ is
bi-associative (any two elements generate an
associative subalgebra) and therefore in a natural way a (right) complex vector space,
see \cite[11.13]{BlauesBuch}.
The norm of $\OO$ is a quadratic form which induces
a positive definite complex hermitian form on $\OO$.
As a unitary $\CC$-basis of $\OO$ we fix the elements
$1,\bj,\bl,\bj\bl\in\OO$, see \cite[11.34]{BlauesBuch}. The $\Aut(\OO)$-stabilizer
of $\bi\in\CC$ acts $\CC$-linearly and can be
identified with the matrix group 
\[
\Aut_\CC(\OO)\cong\SU(3).
\]
In order to specify such an isomorphism with the matrix group,
we use the ordered $\CC$-basis $(\bj,\bl,\bj\bl)$ of $\CC^\perp\subseteq\OO$.
Every element of $u\in\OO$ has a unique expression as
\[
u=u_0+\bj u_1+\bl u_2 + (\bj\bl)u_3, \qquad\text{ with }u_0,u_1,u_2,u_3\in\CC.
\]
We define an isomorphism $\CC^\perp\rTo^\cong\CC^3$ via 
\[
 \bj u_1+\bl u_2 + (\bj\bl)u_3,\mapstoo (u_1,u_2,u_3).
\]
For $z\in\CC$ and $u\in\CC^\perp$ we have $zu=u\bar z$. In what follows,
we will use this identity frequently. The $\CC$-component of a product 
\[
 uv=(u_0+\bj u_1+\bl u_2 + (\bj\bl)u_3)(v_0+\bj v_1+\bl v_2 + (\bj\bl)v_3)
\]
is therefore $u_0v_0-(\bar u_1v_1+\bar u_2v_2+\bar u_3v_3)$. The group $\Aut_\CC(\OO)$
preserves this quantity and fixes $u_0v_0$. Hence $\Aut_\CC(\OO)$ preserves the hermitian form
$\bar u_1v_1+\bar u_2v_2+\bar u_3v_3$ on $\CC^\perp\cong\CC^3$.
It follows that $\Aut_\CC(\OO)\cong\SU(3)$ acts via standard matrix multiplication from
the left on $\CC^3\cong\CC^\perp$.
\end{Num}
\begin{Num}\textbf{\boldmath The model of $\OOP^2$\ }
\label{OP2Definition}
We view the Cayley plane $\OOP^2$ as the set of all idempotent hermitian
$3\times 3$-matrices over $\OO$ with trace $1$ (the \emph{rank~$1$ projectors}).
This is slightly different from \ref{PolarRepresentationOfProjectiveGeometry}
above, where we considered traceless hermitian matrices. The change of the trace
simplifies matrices without changing the stabilizers. 
The euclidean inner product of two $\OO$-hermitian $3\times 3$-matrices is
defined as
\[
 \bra{X,Y}=\mathrm{trace}(XY)=
 \sum_{i=1}^3X_{i,i}Y_{i,i}+2\sum_{1\leq i<j\leq 3}\mathrm{Re}(X_{i,j}\bar Y_{i,j}).
\]
The euclidean distance between two elements $\xi,\xi'\in\OOP^2$ is given by
\[
 ||\xi-\xi'||^2=2-2\bra{\xi,\xi'},
\]
because $||\xi||^2=||\xi'||^2=1$.
The Cayley plane $\OOP^2$ is in particular a Riemannian submanifold of the
$26$-dimensional euclidean space of $\OO$-hermitian $3\times 3$-matrices
with trace $1$, and the compact group $\mathrm F_4$ acts transitively and isometrically
on $\OOP^2$.
A point with
affine coordinates $(x,y)\in\OO\times\OO$ is identified with the projector
\[
\frac1{x\bar x+y\bar y+1}
\begin{pmatrix}
x \\ y \\1
\end{pmatrix}
\begin{pmatrix}
\bar x &\bar y & 1
\end{pmatrix}
=
\frac1{x\bar x+y\bar y+1}
\begin{pmatrix}
x\bar x & x\bar y & x\\
y\bar x & y\bar y & y\\
\bar x & \bar y &1
\end{pmatrix},
\]
see \cite[p.~84]{BlauesBuch}. By means of this coordinate chart we view the 
\emph{affine Cayley plane} $\OO\times\OO$ as an open dense subset of $\OOP^2$.
We note that under this chart the image of a real line in 
$\OO\times\OO$ passing through the origin is a geodesic in $\OOP^2$. Also,  
the chart is conformal at the origin $(x,y)=(0,0)$, as is easily seen by
differentiating.
The complement of the range of the chart is the cut locus $L$ of the point $(0,0)$ in $\OOP^2$,
or, in terms of projective geometry, the projective line at infinity
of the affine Cayley plane $\OO\times\OO$, an $8$-sphere.
\end{Num}
\begin{Num}\textbf{\boldmath The action of $\SU(3)\times\SU(3)$ on $\OOP^2$\ }
The group $\SU(3)$ acts in the standard way isometrically on the set of all hermitian
$3\times 3$-matrices over $\OO$ with trace $1$, preserving $\OOP^2$, and
with $\CCP^2$ as one orbit, via
\[
 g(X)=gXg^{-1}.
\]
Due to the bi-associativity of $\OO$, this matrix product is well-defined.
The action is faithful, since $\OO$ is not commutative.
On the other hand, $\SU(3)=\Aut_\CC(\OO)$ acts as in \ref{SU(3)Action}
entry-wise on the $\OO$-hermitian matrices.
In this way, the compact group
\[
K=\Aut_\CC(\OO)\times\SU(3)=\SU(3)\times\SU(3)
\]
acts isometrically on $\OOP^2$. 
Our aim is to understand the orbit structure of this action.
We begin with the point 
\[
q=(0,0)\in\OO\times\OO.
\]
\end{Num}
\begin{Num}\textbf{\boldmath The $\CCP^2$-orbit and its normal isotropy representation\ }
The affine coordinates $(0,0)$ are complex and therefore the $K$-orbit of $q$ is 
$K(q)=\CCP^2\subseteq\OOP^2$. Since $\Aut_\CC(\OO)$ acts trivially on
$\CCP^2$, the $K$-stabilizer of $q$ is isomorphic to $\Aut_\CC(\OO)\times\mathrm U(2)$.
The projector corresponding to $q=(0,0)$ is 
$\left(\begin{smallmatrix}
0 & 0 & 0\\
0 & 0 & 0\\
0 & 0 & 1
\end{smallmatrix}\right)$ and thus $K_q$ consists of the block matrices of the form
\[
 Y_1
\times
\begin{pmatrix}
 Y_2 & 0\\
 0 & y
\end{pmatrix}
\in\SU(3)\times\SU(3),
\]
with $Y_1\in\SU(3)$ and $Y_2\in\mathrm U(2)$, and $y=\overline{\det(Y_2)}$.
The group $K_q$ stabilizes the polar line (the cut locus) $L$ of
$q$ in $\OOP^2$ and acts on the affine Cayley plane
$\OO\times\OO$. In this way we are reduced to a linear action.
The representation of $K_q$ on $\OO\times\OO$ splits off a representation
of (real) dimension $4$ on $\CC\times\CC$, with $\Aut_\CC(\OO)$ acting trivially and
$\mathrm U(2)$ acting via matrix multiplication from the left by $\det(Y_2)Y_2$ on $\CC^2$.

On the complement $\CC^\perp\times\CC^\perp\cong\CC^{2\times 3}$ we have the
following representation. We represent a point with affine coordinates
\[
(u,v)=(\bj u_1+\bl u_2 + (\bj\bl)u_3, \bj v_1+\bl v_2 + (\bj\bl)v_3)\in\CC^\perp\times\CC^\perp
\]
as
\[
\begin{pmatrix}
 u_1 & v_1\\
 u_2 & v_2\\
 u_3 & v_3
\end{pmatrix}
\in\CC^{2\times 3}
\]
Then $\Aut_\CC(\OO)\cong\SU(3)$ acts in the standard way from the left on $\CC^{2\times 3}$.
We put
\[
Y_2=\begin{pmatrix}
               c & -a\bar s\\
               s & \phantom-a\bar c
              \end{pmatrix},
\]
where $c,s,a$ are complex numbers with $c\bar c+s\bar s=a\bar a=1$, and $y=\bar a$.
Then $\left(\begin{smallmatrix}
 Y_2 & 0\\
 0 & y
\end{smallmatrix}\right)$ maps the point $(u,v)\in\CC^\perp\times\CC^\perp$ to 
\[
(cua-a\bar sva,sua+a\bar cva)=(u\bar ca-v s,u\bar sa+vc)=
(u,v)\begin{pmatrix}
               \bar ca & \bar s a\\
               -s & c
              \end{pmatrix}
              =
(u,v)\begin{pmatrix}
               \bar c & \bar s \\
               -s\bar a & c\bar a
              \end{pmatrix}a              
              .
\]
Hence $Y_2$ acts on $\CC^{2\times3}$ through matrix multiplication from the right by $Y_2^*\bar y$,
where $Y_2^*=\bar Y_2^T$ and $\bar y=\det(Y_2)$. Summing this up, we have the $K_q$-action
\[
 Y_1
\times
\left(
\begin{smallmatrix}
 Y_2 & 0\\
 0 & y
\end{smallmatrix}
\right)
:
\left(
\begin{smallmatrix}
 u_1 & v_1\\
 u_2 & v_2\\
 u_3 & v_3
\end{smallmatrix}\right)
\mapstoo
Y_1 \left(\begin{smallmatrix}
 u_1 & v_1\\
 u_2 & v_2\\
 u_3 & v_3
\end{smallmatrix}\right)Y_2^*\bar y.
\]
In particular, the action of $K_q$ on $\CC^{2\times 3}$ is orbit equivalent with the polar action
described in \ref{Veronesean}.
Since the tangent space $T_q\OOP^2$ splits also $K_q$-equivariantly
as $T_q\OOP^2=T_q\CCP^2\oplus\bot_q\CCP^2$, this gives us at the same time the
normal isotropy representation of $K_q$. In particular, the normal
isotropy representation of $K_q$ on $\bot_q\CCP^2\subseteq T_q\OOP^2$ is polar. 

We put 
\[
d=(0,\bj)\qquad\text{and}\qquad p=(-\bl,\bj),
\]
Suppose that $\lambda,\mu$ are nonnegative reals and consider the point 
\[
o=o_{\lambda,\mu}=p\lambda+d\mu=(-\lambda\bl,(\lambda+\mu)\bj)\in\CC^\perp\times\CC^\perp.
\]
The corresponding point in $\CC^{2\times 3}$ is $\left(\begin{smallmatrix}
                                                        \phantom-0 & \lambda+\mu\\-\lambda&0\\\phantom-0&0
                                                       \end{smallmatrix}\right)$.
Since the action of $K_q$ on $\bot_q\CCP^2$ is polar,
every $K_q$-orbit in $\bot_q\CCP^2$ contains exactly one
such point $o_{\lambda,\mu}$. 
The $K_q$-stabilizer of $o_{\lambda,\mu}$ can
easily be computed. 
For $\lambda,\mu>0$, it consists of the matrices 
\[\left(\begin{smallmatrix}
 \bar z & \\
 & zy \\
  && \bar y
\end{smallmatrix}\right)
\times
\left(\begin{smallmatrix}
  z & \\
 & \overline{zy} \\
  &&  y
\end{smallmatrix}\right)\in\SU(3)\times\SU(3),
\]
with $y,z\in\mathrm U(1)$. For $\lambda=0<\mu$, it consists of the matrices
\[
\bigl(\begin{smallmatrix}
  \bar z & \\
  &  Z_2
\end{smallmatrix}\bigr)
\times
\left(\begin{smallmatrix}
 z & \\
 & \overline{yz} \\
  && y
\end{smallmatrix}\right)
\in\SU(3)\times\SU(3),
\]
with $Z_2\in\mathrm U(2)$.
For $\mu=0<\lambda$ it consists of matrices of the form
\[
\bigl(\begin{smallmatrix}
 \bar Y_2 & \\
  & \bar y
\end{smallmatrix}\bigr)
\times
\bigl(\begin{smallmatrix}
  Y_2 & \\
  &  y
\end{smallmatrix}\bigr)\in\SU(3)\times\SU(3),
\]
with $Y_2\in\mathrm U(2)$.
\end{Num}
\begin{Num}\textbf{Euclidean distances between orbits\ }
The projector $\xi_{\lambda,\mu}\in\OOP^2$ corresponding to $o_{\lambda,\mu}$ is
\[\xi=\xi_{\lambda,\mu}=
\textstyle
\frac1{\lambda^2+(\lambda+\mu)^2+1}
\begin{pmatrix}
\lambda^2 & -\lambda(\lambda+\mu)\bj\bl&-\lambda\bl\\
\lambda(\lambda+\mu)\bj\bl&(\lambda+\mu)^2&(\lambda+\mu)\bj\\
\lambda\bl & -(\lambda+\mu)\bj& 1.
\end{pmatrix}
\]
We note that the off-diagonal entries of this matrix
are all Cayley numbers which are perpendicular to $\CC$. 
We denote the euclidean distance between $\xi$ and $\CCP^2$
by 
\[
\delta(\xi)=\min\{||\zeta-\xi||\mid\zeta\in\CCP^2\}. 
\]
In order to compute this distance, we note that every point in $\CCP^2\subseteq\OOP^2$
is of the form 
\[
\zeta=\begin{pmatrix}
|u|^2 &u\bar v &u\bar w\\
v\bar u&|v|^2&v\bar w\\
w\bar u&w\bar v&|w|^2
\end{pmatrix},
\]
where $u,v,w$ are complex numbers with $|u|^2+|v|^2+|w|^2=1$.
The point $q$ corresponds to $(u,v,w)=(0,0,1)$.
The euclidean inner product between $\xi=\xi_{\lambda,\mu}$ and $\zeta$
is given by
\[
\bra{\xi,\zeta}=\frac{\lambda^2|u|^2+(\lambda+\mu)^2|v|^2+|w|^2}{\lambda^2+(\lambda+\mu)^2+1},
\]
because the off-diagonal entries of $\xi$ are perpendicular to $\CC$.
From this formula we see the following. 
We have 
\[
\delta(\xi)=||q-\xi||\quad\text{ if and only if }\quad\lambda\leq 1\text{ and }\lambda+\mu\leq 1.
\]
This condition defines a linear simplex (recall that $\lambda,\mu\geq 0$).
From the formula for $\bra{\xi,\zeta}$, the following is immediate.

\medskip\noindent
(1) 
If  $\lambda+\mu<1$, then $q$ is the unique point in $\CCP^2$ at distance
$\delta(\xi_{\lambda,\mu})$ from $\xi_{\lambda,\mu}$.
In particular, $K_{\xi_{\lambda,\mu}}\subseteq K_q$.

\medskip\noindent
(2)
If $\lambda+\mu=1\neq\lambda$, then every point with complex coordinates
$(u,v,w)\in\SS^5$ and $u=0$ realizes the distance $\delta$. 
This condition defines a complex projective line in $\CCP^2$.
Also, the point $\tilde q$ with complex coordinates $(u,v,w)=(1,0,0)$ 
is in this case the unique
point in $\CCP^2$ at maximal distance from $\xi_{\lambda,\mu}$, hence
$K_{\xi_{\lambda,\mu}}\subseteq K_{\tilde q}$.

\medskip\noindent
(3)
If $\lambda+\mu=1=\lambda$, then every point $\zeta$ in $\CCP^2$ has
distance $\delta(\xi_{\lambda,\mu})$ from $\xi_{\lambda,\mu}$.
\end{Num}
\begin{Lem}
\label{FundamentalDomain}
Every $K$-orbit contains a unique point $\xi_{\lambda,\mu}$ with
$0\leq\lambda,\mu$ and $\lambda+\mu\leq 1$.

\proof
Let $\eta\in\OOP^2$ and let $\zeta\in\CCP^2$ be a point that has
minimal euclidean distance from $\eta$. There exists $g\in K$ with
$g(\zeta)=q$. Then $g(\eta)$ is not in the cut locus $L$ of $q$,
since $L\cap\CCP^2$ contains points which are strictly closer to 
any given point in $L$ than $q$ (we omit this short calculation). 
Hence $g(\eta)\not\in L$. If one of the off-diagonal entries of the projector
$g(\eta)$ is not perpendicular to $\CC$, then the
inner product shows that $q$ is not the closest point 
to $g(\eta)$ on $\CCP^2$. Thus
$g(\eta)$ is, as a point in $\OO\times\OO$, perpendicular to $\CC\times\CC$.
Since we have a polar action on the normal space of $q$,
there exists $h\in K_q$ such that $hg(\eta)=\xi_{\lambda,\mu}$,
for some $\lambda,\mu\geq 0$. By the observations above, we have $\lambda+\mu\leq 1$.

It remains to show the uniqueness.
Let $\xi_{\lambda,\mu}$ be a point in the simplex.
If $\lambda+\mu<1$ and if $g(\xi_{\lambda,\mu})=\xi_{\lambda',\mu'}$
is in the simplex,
then $g(q)=q$, because $q$ is the unique nearest point to $\xi_{\lambda,\mu}$.
Therefore $(\lambda,\mu)=(\lambda',\mu')$,
because the action of $K_q$ on the normal space is polar.
If $\lambda+\mu=1\neq\lambda$ and if $g(\xi_{\lambda,\mu})=\xi_{\lambda',\mu'}$,
then we see from the geometric description above that $\lambda'+\mu'=1\neq\lambda'$.
The number $\lambda$ is determined by the distance of $\xi_{\lambda,\mu}$
from $\CCP^2$, hence $(\lambda,\mu)=(\lambda',\mu')$. Finally, $p$ is the unique 
point in the simplex that has constant distance from $\CCP^2$.
\qed
\end{Lem}
The uniqueness statement of the previous lemma follows also from the fact
that the $K$-action is polar, which we prove below. 
Also, we have worked with the euclidean distance,
rather than with the inner metric of the Riemannian manifold $\OOP^2$. 
We will come back to this. But first we determine the stabilizers of the
$\xi_{\lambda,\mu}$, where $\lambda+\mu=1$.
\begin{Num}\textbf{The remaining orbit types\ }
Suppose that $\lambda+\mu=1\neq\lambda$. Then the point $\tilde q$ with complex coordinates
$(u,v,w)=(1,0,0)$ uniquely maximizes the euclidean distance from $\xi_{\lambda,\mu}$,
as we noticed above. Thus $K_{\xi_{\lambda,\mu}}\subseteq K_{\tilde q}$. The involution 
$h=\left(\begin{smallmatrix} && 1\\ &-1\\1
\end{smallmatrix}\right)\times
\left(\begin{smallmatrix} && 1\\ &-1\\1
\end{smallmatrix}\right)   \in\SU(3)\times\SU(3)$ 
interchanges $q$ and $\tilde q$. 

For $\lambda\neq0$, it maps $p\lambda +d\mu=(-\lambda\bl,\bj)$ to 
$(-\frac1\lambda\bl,-\frac1\lambda\bj)$. 
The $K$-stabilizer of 
$p\lambda +d\mu$ consists then of the matrices
\[
\bigl(\begin{smallmatrix}
  \bar z& \\
  &  \bar Z_1
\end{smallmatrix}\bigr)
\times
\bigl(\begin{smallmatrix}
  z & \\
  & Z_1
\end{smallmatrix}\bigr)\in\SU(3)\times\SU(3),
\]
with $Z_1\in\mathrm U(2)$.
By continuity, these matrices also fix $p$ and $d$.

Suppose that $\lambda=0$. The involution maps the projector coresponding to 
$d=(0,\bj)$ to 
\[
\theta ={\textstyle\frac12}\left(\begin{smallmatrix}
               \phantom-1 \\-\bj\bl\\\phantom-0
              \end{smallmatrix}
\right)(\begin{smallmatrix}1,\bj\bl,0\end{smallmatrix})=
{\textstyle\frac12}\left(\begin{smallmatrix}
                                  1 & \bj\bl &0 \\
                                  \overline{\bj\bl} &1 & 0\\
                                  0&0&0
                                 \end{smallmatrix}
\right).
\]
The element $\id\times\left(\begin{smallmatrix}
                                     c & -a\bar s \\
                                     s & \phantom-a\bar c\\
                                     && \bar a
                                    \end{smallmatrix}\right)
\in\Aut_\CC(\OO)\times\SU(3)$ maps $\theta$ to
$\frac12\left(\begin{smallmatrix}
                                  1 & \bj\bl a &0 \\
                                  \overline{\bj\bl a} &1 & 0\\
                                  0&0&0
                                 \end{smallmatrix}
\right)$ and therefore
$K_d$ consists of the matrices of the form
\[
\bigl(\begin{smallmatrix}
  \bar z& \\
  &  Z_2
\end{smallmatrix}\bigr)
\times
\bigl(\begin{smallmatrix}
  z & \\
  & Z_1
\end{smallmatrix}\bigr)\in\SU(3)\times\SU(3),
\]
with $Z_1,Z_2\in\mathrm U(2)$.

\end{Num}
\begin{Lem}
We have $K_p=\left\{\left. \bar X\times X\right| X\in\SU(3)\right\}$.

\proof
Let $H=\left\{\left.\bar X\times X\right| X\in\SU(3)\right\}$.
The block matrices 
$\left(\begin{smallmatrix} \bar Y & \\ & 1\end{smallmatrix}\right)\times
\left(\begin{smallmatrix} Y & \\ & 1\end{smallmatrix}\right)$
and 
$\left(\begin{smallmatrix} 1 & \\ & \bar Y\end{smallmatrix}\right)\times
\left(\begin{smallmatrix} 1 & \\ & Y\end{smallmatrix}\right)$
with $Y\in\SU(2)$ fix $p$ and generate $H$, 
hence $H$ fixes $p$. The Lie algebra of $H$ is maximal in
$\fsu(3)\oplus\fsu(3)$, because $\fsu(3)$ is simple. 
Thus $H=(K_p)^\circ$. If $1\times A$ is in the kernel of $\pr_1:K_p\rTo\SU(3)$,
then $A\in\Cen(\SU(3))$. Such an element fixes $p$ only if $A=1$.
Similarly, if $A\times 1$ fixes $p$, then $A=1$. It follows that $H=K_p$.
\qed
\end{Lem}
\begin{Num}\textbf{The corresponding complex of groups\ }
\label{DefinitionOfExceptionalGroupComplex}
The kernel of the $K$-action is the group $Z=K_{p,d,q}\cap\Cen(\SU(3)\times\SU(3))\cong\ZZ/3$
and we put $G=K/Z$.
The simple complex of groups in $K$ formed by the seven types of stabilizers,
corresponding to the faces of the simplex
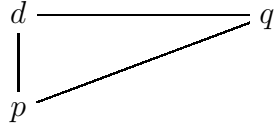
\begin{diagram}[height=1.5em,width=4em] 
d & \rLine & q\\
\dLine & \ruLine\\
p
\end{diagram}
looks as follows. 
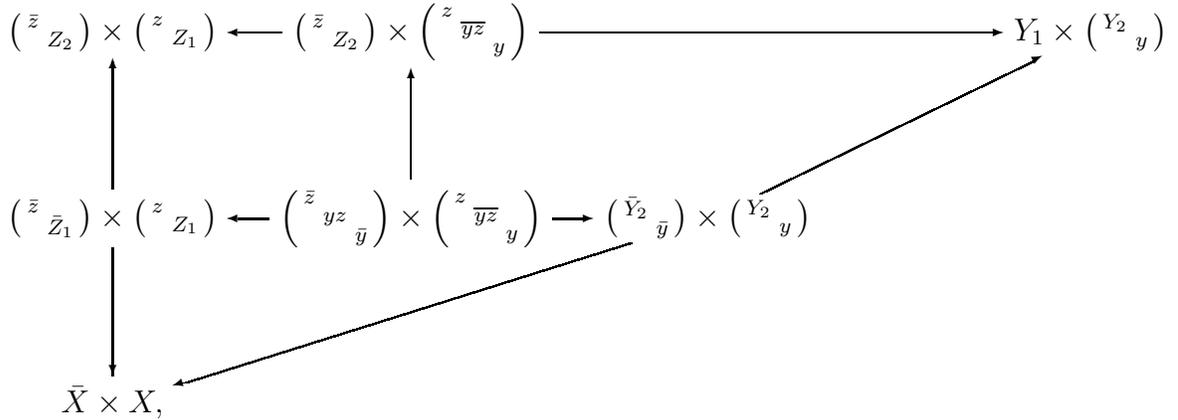
\begin{diagram}
\bigl(\begin{smallmatrix}
  \bar z& \\
  &  Z_2
\end{smallmatrix}\bigr)
\times
\bigl(\begin{smallmatrix}
  z & \\
  & Z_1
\end{smallmatrix}\bigr)
&\lTo& 
\bigl(\begin{smallmatrix}
  \bar z & \\
  & Z_2
\end{smallmatrix}\bigr)
\times
\left(\begin{smallmatrix}
 z & \\
 & \overline{yz} \\
  && y
\end{smallmatrix}\right)
&\rTo &&&&&
Y_1
\times
\bigl(\begin{smallmatrix}
 Y_2 & \\
  & y
\end{smallmatrix}\bigr)
\\
\uTo && \uTo &&&&& \ruTo(4,2) \\
\bigl(\begin{smallmatrix}
  \bar z& \\
  &  \bar Z_1
\end{smallmatrix}\bigr)
\times
\bigl(\begin{smallmatrix}
  z & \\
  &  Z_1
\end{smallmatrix}\bigr)
&\lTo &
\left(\begin{smallmatrix}
 \bar z & \\
 & yz \\
  && \bar y
\end{smallmatrix}\right)
\times
\left(\begin{smallmatrix}
  z & \\
 & \overline{yz} \\
  &&  y
\end{smallmatrix}\right)
&\rTo&
\bigl(\begin{smallmatrix}
 \bar Y_2 & \\
  & \bar y
\end{smallmatrix}\bigr)
\times
\bigl(\begin{smallmatrix}
  Y_2 & \\
  & y
\end{smallmatrix}\bigr)
\\
\dTo&&&\ldTo(4,2)\\
 \bar X
\times
 X,
\end{diagram}
with $Y_2,Z_1,Z_2\in\mathrm U(2)$ and $Y_1,X\in\SU(3)$.
The corresponding simple complex of groups in $G$ is obtained by taking matrices mod $Z$.
\end{Num}
An isometric action of a Lie group $G$ 
on a complete Riemannian manifold $M$ is called \emph{polar}
if there exists a complete submanifold $\Sigma\subseteq M$ which meets every orbit
orthogonally, i.e. 
\[
G(\Sigma)=M\quad\text{ and }\quad T_\sigma\Sigma\perp T_\sigma G(\sigma)
\text{ holds for every }\sigma\in\Sigma. 
\]
This is the case for our action.
We define an immersion $\sigma:\SS^2\rTo\OOP^2$ by putting
\[
 \sigma(x,y,z)=
\left(\begin{smallmatrix} x^2&-\bj\bl xy&-\bl zx\\\bj\bl xy&y^2&\bj yz\\\bl zx&-\bj yz&z^2\end{smallmatrix}\right)\in\OOP^2
\]
and we put $\Sigma=\sigma(\SS^2)$. The surface $\Sigma$ is isometric to $\RRP^2$.
The following is proved in \cite{PTh}.
\begin{Thm}[Podest\`a-Thorbergsson]
The action of $G=(\SU(3)\times\SU(3))/Z$ on $\OOP^2$ is polar and $\Sigma$ is a section.

\proof
The simplex which we considered above is contained in $\Sigma$, hence 
$G(\Sigma)=\OOP^2$ by \ref{FundamentalDomain}.
Let $\xi=\sigma(x,y,z)\in\Sigma$. We claim that $T_\xi\Sigma\perp T_\xi G(\xi)$.

Let $\dot g\in\fsu(3)$. We view $\dot g$ as an element of $\fsu(3)\oplus0\subseteq\Lie(\SU(3)\times\SU(3))$. 
Then $\dot g$ acts via ordinary $3\times 3$-matrix multiplication as 
$\xi\mapstoo \dot g\xi-\xi \dot g\in T_\xi G(\xi)$. 
Now let $\dot\xi\in T_\xi\Sigma$.
A short and elementary calculation shows that the matrix product $\xi\dot\xi$ is a
matrix whose off-diagonal entries are Cayley numbers 
perpendicular to $\CC$, while the entries
on the diagonal are real. The same holds for $\dot\xi\xi$.
We denote the real part of a Cayley number $a$ by $\mathrm{Re}(a)$
and extend this entry-wise to matrices. Then we have 
\begin{align*}
\bra{\dot g\xi-\xi\dot g,\dot\xi}&=\mathrm{trace}((\dot g\xi-\xi\dot g)\dot\xi)\\
&=\mathrm{Re}(\mathrm{trace}(\dot g\xi\dot\xi)-\mathrm{trace}(\xi\dot g\dot\xi))\\
&=\mathrm{Re}(\mathrm{trace}(\dot g\xi\dot\xi))-\mathrm{Re}(\mathrm{trace}(\xi\dot g\dot\xi))\\
&=\mathrm{trace}(\mathrm{Re}(\dot g\xi\dot\xi))-\mathrm{trace}(\mathrm{Re}(\dot g\dot\xi\xi))\\
&=0.
\end{align*}
Now let $\dot h$ be an element of $\Lie(\Aut_\CC(\OO))=\fsu(3)$.
Because $\dot h$ has imaginary entries
on its diagonal, we have $\bra{\dot h(\bj),\bj}=\bra{\dot h(\bl),\bl}=\bra{\dot h(\bj\bl),\bj\bl}=0$.
On the three real diagonal entries of $\xi$, the infinitesimal automorphism $\dot h$ acts as multiplication by $0$.
Therefore $\bra{\dot\xi,\dot h(\xi)}=0$. 

This shows that $T_\xi\Sigma\perp T_\xi G(\xi)$.
\qed
\end{Thm}
\begin{Num}\textbf{The Riemannian metric\ }
The linear simplex which we considered in $\OO\times\OO$ is contained in $\Sigma$
and has geodesic edges (and constant curvature) in $\OOP^2$. 
The quotient $G\backslash\OOP^2$ is isometric to a spherical simplex of shape $\mathsf C_3$.
\end{Num}
The previous results give us a geometric description of the orbits
$G(d)$ and $G(p)$. The orbit $G(p)$ consists of all points in $\OOP^2$ having
maximal (inner or euclidean) distance from $\CCP^2$. The orbit $G(d)$ consists of all points which
have the property that a (euclidean or inner-metric) ball around them touches $\CCP^2$ in a
$2$-sphere, and which have maximal distance from $\CCP^2$ with respect to this
property. We remark that the embedding of $\CCP^2$ is \emph{tight}:
every euclidean ball that touches $\CCP^2$ does this either in a unique
point, along a $2$-sphere, or everywhere.
\begin{Prop}
\label{ExceptionalGeometryProp}
Let $\Delta$ denote the simplicial complex whose nerve is the
covering of $G$ by the cosets of $G_p$, $G_d$ and $G_q$, as defined in
\ref{DefinitionOfExceptionalGroupComplex}.
Then $(G,\Delta)$ is a homogeneous compact geometry of type
$\mathsf C_3$ which is not a building. We have $|\Delta|_K=\OOP^2$.

\proof
We can identify the nonempty simplices with the cosets of the various
$G_\alpha$, for $\emptyset\neq\alpha\subseteq\{p,d,q\}$.
From the diagram in \ref{DefinitionOfExceptionalGroupComplex}
above it is clear that the link of $G_p$ is
isomorphic to the $2$-dimensional complex projective geometry.
From \ref{Veronesean} we see that the link of $G_q$ is isomorphic to the generalized
quadrangle corresponding to the hermitian form
$h=(-f_2)\oplus f_3$ on $\CC^{2+3}$. The link of $G_d$ is isomorphic to
the generalized digon $\SS^2\lTo\SS^2\times\SS^3\rTo\SS^3$. 
In particular, $\lk(d)$ is a complete bipartite graph. 
It follows that every triangle in the $1$-skeleton $\Delta^{(1)}$
which contains $d$ is filled by a $2$-simplex. From the transitive action of
$G$ we conclude that $\Delta$ is a flag complex.
From the diagram \ref{DefinitionOfExceptionalGroupComplex} we see that $G=G_pG_q$.
Thus $\Delta$ is (gallery-) connected and by 
\ref{LocalRecognitionOfGeometry} a geometry of type $\mathsf C_3$. 

Since $G=G_pG_q$, the
plane stabilizer $G_p$ acts transitively on the set of points $G/G_q$.
In other words, a point and a plane in $\Delta$ are always incident
(such geometries are called \emph{flat} in \cite{Pasini}).
This cannot hold in a polar space. 

Finally we note that that we have a $G$-equivariant bijective map 
$|\Delta|_K\rTo\OOP^2$ which sends 
$g(G_p\cdot\lambda+G_d\cdot\mu+G_q\cdot\nu)\in|\Delta|$
to $g(\xi_{\lambda,\mu,\nu})\in\OOP^2$.
\qed
\end{Prop}
\begin{Thm}
\label{UniquenessOfException}
Let $\cG$ denote the simple complex of groups from \ref{DefinitionOfExceptionalGroupComplex}.
Up to isomorphism, there is exactly one homogeneous
compact geometry $(G,\Delta)$ of type $\mathsf C_3$ belonging
to this complex of groups.

\proof
Let $(\widehat G,\widehat\Delta)$ denote the universal homogeneous
compact geometry for $(G,\Delta)$, as in \ref{ConstructionOfGHat}. 
We denote the vertex stabilizers corresponding to $\cG\rTo\widehat G$ by 
$\widehat G_\alpha$, for $\alpha\subseteq\{p,d,q\}$.
We have by \ref{TheBasicCosetConstruction} a surjective equivariant map
$(\widehat G,\widehat\Delta)\rTo(G,\Delta)$. Let $F\subseteq\widehat G$ denote
its kernel. Thus $\Lie(\widehat G)\cong\Lie(G)\oplus\Lie(F)$. 
Let
$\pr_2:\Lie(\widehat G)\rTo\Lie(F)$ denote the projection onto the
second summand and suppose that $\Lie(F)\neq 0$.
Since $\widehat G$ is generated by $\widehat G_p\cup\widehat G_q$,
see \ref{FiniteStepGeneration}, 
either $\pr_2(\Lie(\widehat G_p))\neq 0$ or $\pr_2(\Lie(\widehat G_q))\neq 0$.
Moreover, we have $\dim(F)\leq 7$ by \ref{DimensionBound}. Since $\fsu(3)$ is
simple and $8$-dimensional, we have
$\pr_2(\Lie(\widehat G_p))=0$ and $\pr_2:\Lie(\widehat G_q)\rTo\Lie(F)$
annihilates the $\fsu(3)$-summand. From the diagram above and
the fact that the $\pr_2$-image of $\Lie(\widehat G_q)$ is nontrivial, we see
that $\pr_2$ is not trivial on $\Lie(\widehat G_{p,q})$. This is a
contradiction, since $\Lie(\widehat G_{p,q})\subseteq\Lie(\widehat G_p)$.
Thus $\Lie(F)=0$ and $F$ is finite. Since $F$ acts freely by
\ref{TheBasicCosetConstruction},
$F\subseteq\pi_1(|\Delta|_K)=\pi_1(\OOP^2)=1$. This shows that $(G,\Delta)$
is universal. 

Finally, we note that the $\ZZ/2$-Lefschetz number of every self-homeomorphism $\phi$ of
$\OOP^2$ is $1$, hence $\phi$ has a fixed point.
Therefore $|\Delta|_K$ admits no continuous free action and in
particular no quotients, see eg.~Brown \cite[p.~42]{Brown} or \cite[55.19]{BlauesBuch}.
\qed
\end{Thm}
The following result is a consequence of our classification in Section~\ref{Section4}.
\begin{Prop}
Suppose that $(G,\Delta)$ is the exceptional compact homogeneous 
$\mathsf C_3$ geometry from \ref{ExceptionalGeometryProp}
and suppose that a compact connected group $H$ acts continuously, faithfully
and transitively on the chambers. Then $H$ is conjugate to the group $G$
in the group of topological automorphisms of $\Delta$.

\proof
The group $H$ is a compact connected Lie group by \ref{GIsLieGroup}. 
We consider the chamber $\gamma=\{p,d,q\}$.
The fundamental group of the set of chambers $G/G_\gamma$ is finite.
Therefore the semisimple commutator group $K=[H,H]$ acts transitively on the
chambers, see \cite[p.~94]{Oni}. From the long exact homotopy sequence
for the transitive action of $K$ on $K/K_p\cong\SU(3)$ we see that
$K_p$ is connected and semisimple. Similarly, we see that from the transitive
action of $K$ on $K/K_q\cong\CCP^2$ that $K_q$ has a $1$-dimensional center.
This is all that is needed in \ref{2*1+1} in order to determine the 
possibilities for the simple
complex of groups $\mathcal K$ formed by the stabilizers.
Thus there are at most two possibilities for $\mathcal K$,
and the corresponding universal homogeneous compact geometries are 
by \ref{C3ClassificationTheorem} either a polar
space or the exceptional geometry. Since $\Delta$ is not covered
by a building, the compact universal covering is $(G,\Delta)$.
Therefore we have a continuous isomorphism $(G,\Delta)\rTo(K,\Delta)$.
Finally, we have $H=K$ because the connected $K$-normalizer of $K_p$ is
$K_p$, see \ref{RecoverWholeGroup}.
\qed
\end{Prop}
\begin{Rem}
There is another way to approach the exceptional geometry.
Starting from the fact that the $G$-action on $\OOP^2$ is polar and has a
spherical simplex as its metric orbit space, one may consider the simple
complex of groups formed by the stabilizers corresponding to the faces of the simplex.
The horizontal simplicial complex corresponding to the action
can be shown to be a compact geometry
of type $\mathsf C_3$, whose coarse realization is homeomorphic to $\OOP^2$, see
Lytchak \cite{Lytchak} and Fang-Grove-Thorbergsson \cite{FGTh}. The proof of the main theorem
in  \cite{Lytchak} shows that this geometry cannot be covered by a building.   
The classification in the following sections shows that there
is at most one candidate for such a simple complex of groups.
Therefore, this candidate must describe our geometry.  This, very implicit way, can be used to
obtain the stabilizer of  our  action without   explicit computations.
The proof that this action is polar with a nice quotient, however, seems to require some
calculations, as in \cite[p.~151--154]{PTh}.
\end{Rem}

\section{\boldmath The classification of the universal homogeneous compact geometries
of type $\mathsf C_3$}
\label{Section4}

Our aim in this section is the classification of the 
universal homogeneous compact geometries of type $\mathsf C_3$.
The main result of this section is as follows.
\begin{Thm}
\label{C3ClassificationTheorem}
Let $(G,\Delta)$ be a homogeneous compact geometry of type $\mathsf C_3$ with connected
panels. Assume that $G$ is compact and acts faithfully, and let 
$(\widehat G,\widehat\Delta)$ denote the corresponding universal compact 
homogeneous geometry, as in \ref{ConstructionOfGHat}.
Then $\widehat\Delta$ is either a building or the exceptional geometry described in
Section~\ref{TheExceptionalC3Geometry}.
\end{Thm}
In order to prove this theorem, we classify the possibilities for the
simple complex of groups $\cG$. In view of \ref{UniversalAndMinimal}
we assume that the homogeneous geometry $(G,\Delta)$ is both minimal
and universal. The proof will be given at the end of Section \ref{Section4}.
\begin{Num}\textbf{Notation\ }
\label{TerminologyForC3}
We fix some notation that will be used throughout Section~\ref{Section4}.
We assume that $(G,\Delta)$ is a homogeneous compact geometry of type
$\mathsf C_3$ with connected panels. The Lie group $G$ is compact and
connected and acts faithfully.
In the geometry $\Delta$ we have three types of vertices,
called \emph{points}, \emph{lines} and \emph{planes}.
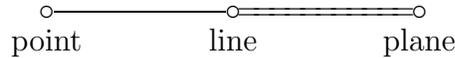
\begin{diagram}[height=1em,abut]
\circ & \rLine & \circ & \rEq & \circ\\
\text{point}&&\text{line}&&\text{plane}
\end{diagram}
We fix a chamber $\gamma=\{p,d,q\}$, where $p$ is a plane, $d$ is a line and $q$ is a point
(so $\lk(p)$ is a projective \textbf{p}lane, $\lk(d)$ is a generalized \textbf{d}igon and
$\lk(q)$ is a generalized \textbf{q}uadrangle). 
The simple complex of groups $\cG$ looks as follows.
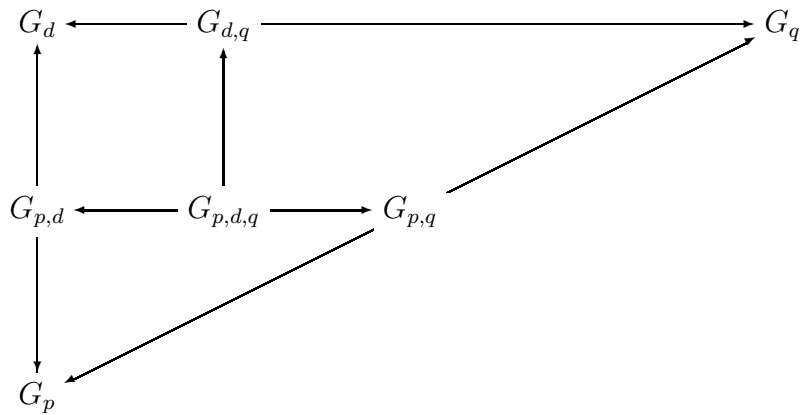
\[
\begin{diagram}
G_d &\lTo& G_{d,q}&\rTo &&&&& G_q\\
\uTo && \uTo &&&&& \ruTo(4,2) \\
G_{p,d}&\lTo &G_{p,d,q}&\rTo& G_{p,q}\\
\dTo&&&\ldTo(4,2)\\
G_p
\end{diagram}
\]
We note also that 
\[
G=\bra{G_p\cup G_q}=\bra{G_p\cup G_d}=\bra{G_q\cup G_d},
\]
because $G=\bra{G_{p,d}\cup G_{d,q}\cup G_{q,p}}$.
The link $\lk(p)$ is one of the four compact Moufang
planes over $\RR,\CC,\HH$ or $\OO$.
Accordingly, the panels $\lk(\{p,d\})$ and $\lk(\{p,q\})$
are, in the coarse topology, spheres of dimension $m=1,2,4,8$.
The link $\lk(q)$ is a compact connected Moufang quadrangle,
the panel $\lk(\{d,q\})$ is an $n$-sphere, and
the link $\lk(d)$ is a generalized digon.
It is given by the two $G_d$-equivariant maps 
\[
\SS^m\lTo\SS^m\times\SS^n\rTo\SS^n
\]
as in \ref{StructureOfD}.
We have 
\[
\dim(G/G_\gamma)\leq 6m+3n
\]
by \ref{DimensionBound}.
If $m,n\geq 2$, then $G$ is semisimple by \ref{CommutatorIsTransitive}.
\end{Num}
\begin{Num}\textbf{\boldmath Homotopy properties of $\cG$\ }
\label{HomotopyPropertiesOfcG}
Recall that a continuous map is called a \emph{$k$-equivalence} if it induces an isomorphism on
the homotopy groups in degrees less than $k$ and an epimorphism in degree $k$.
The following diagram shows the low-dimensional homotopy properties of the maps in $\cG$.
\label{ConnDiagram}
\begin{diagram}[nohug]
G_d &\lTo_{\scriptstyle (m-1)\text{-equiv.}}& G_{d,q}&\rTo_{\scriptstyle (m-1)\text{-equiv.}} &&&&& G_q\\
\uTo_{\scriptstyle (n-1)\text{-equiv.}} && \uTo_{\scriptstyle (n-1)\text{-equiv.}} &&&&& \ruTo(4,2)_{\scriptstyle (n-1)\text{-equiv.}} \\
G_{p,d}&\lTo^{\scriptstyle (m-1)\text{-equiv.}} &G_{p,d,q}&\rTo^{\scriptstyle (m-1)\text{-equiv.}}& G_{p,q}\\
\dTo_{\scriptstyle (m-1)\text{-equiv.}}&&&\ldTo(4,2)_{\scriptstyle (m-1)\text{-equiv.}}\\
G_p
\end{diagram}
They follow from the fact that the quotients of the various isotropy
groups are spheres, products of spheres or compact generalized polygons.
For example, $G_q/G_{d,q}$ is the point space of a compact generalized quadrangle
with topological parameters $(m,n)$ and admits therefore a CW decomposition
$G_q/G_{d,q}=e^0\cup e^m\cup e^{m+n}\cup e^{m+n+m}$, see \cite[3.4]{KramerDiss}.
The long exact homotopy sequence of the fibration
$G_{d,q}\rTo G_q\rTo G_q/G_{d,q}$ yields the $(m-1)$-connectivity of the 
homomorphism $G_{d,q}\rTo G_q$. The reasoning for the other homomorphisms is
similar, using the results in \emph{loc.cit.}
\end{Num}
We now consider the homotopy groups $\pi_0$ and $\pi_1$.
A compact connected Lie group is divisible (because tori are divisible and
every element is contained in some torus, see eg.~\cite[6.30 or 9.35]{HMCompact}).
This implies the following. If $H$ is a compact Lie group and if
$\phi:H\rTo F$ is a homomorphism to a finite group $F$, then
$\phi$ factors through $\pi_0(H)=H/H^\circ$. In particular, $\phi$ is automatically continuous.
\begin{Lem}
\label{IsotropyConnectedLemma}
If $m,n>1$, then all seven isotropy groups appearing in $\cG$ are connected.

If $m>n=1$, then $\pi_0(G_d)=\pi_0(G_{d,q})=\pi_0(G_q)=1$.

If $n>m=1$, then $\pi_0(G_p)=1$.

\proof
If $m,n>1$, then all maps in $\cG$ are $1$-equivalences and induce
therefore isomorphisms on $\pi_0$.
By the universal property \ref{UniversalProperty} of $G$,
there is a homomorphism $G\rTo\pi_0(G_{p,d,q})$ which
is surjective, because the natural map
$G_{p,q,z}\rTo\pi_0(G_{p,d,q})$ is surjective.
The group $G$ is connected and therefore $\pi_0(G_{p,d,q})=1$.
If $m>n=1$, then \ref{ConnDiagram} shows similarly 
that $\pi_0(G_d)=\pi_0(G_{d,q})=\pi_0(G_q)=1$. The case $n>m=1$ is analogous.
\qed
\end{Lem}
We need a similar result for the fundamental groups in order to control
the torus factors of the stabilizers. This requires some low-dimensional
cohomology.
\begin{Lem}
\label{CofunctorLemma}
Let $K$ be a compact connected Lie group. There are natural isomorphisms
\begin{diagram}
\Hom(K,\SS^1)&\rTo^\cong& H^1(K)&\rTo^\cong&  \Hom(\pi_1(K),\ZZ)\\
\uTo^\cong&&\uTo^\cong&&\uTo^\cong\\
\Hom(K/[K,K],\SS^1)&\rTo^\cong& H^1(K/[K,K])&\rTo^\cong  &\Hom(\pi_1(K/[K,K]),\ZZ),
\end{diagram}
where $H^1$ denotes $1$-dimensional
singular cohomology and $\Hom(K,\SS^1)$ denotes the group of
continuous homomorphisms $K\rTo\SS^1$.

\proof
For every path-connected space $X$ we have by the Universal Coefficient Theorem
\[
H^1(X)\cong \Hom(H_1(X),\ZZ)\cong \Hom(\pi_1(X),\ZZ),
\]
see eg.~\cite[XII.4.6 and VIII.7.1]{Massey}.
We note that $K/[K,K]$ is a torus, hence
$\pi_1(K/[K,K])$ is free abelian.
The fundamental group of the semisimple group
$[K,K]$ is finite, see \cite[5.76]{HMCompact}. 
From the split short exact sequence
\[
1\rTo\pi_1([K,K])\rTo\pi_1(K)\rTo\pi_1(K/[K,K])\rTo1
\]
we have therefore an isomorphism 
\[
\Hom(\pi_1(K),\ZZ)\lTo^\cong \Hom(\pi_1(K/[K,K]),\ZZ).
\]
Since $\SS^1$ is abelian, we also have a natural isomorphism
$\Hom(K,\SS^1)\lTo^\cong\Hom(K/[K,K],\SS^1)$. 
Since $\SS^1\simeq K(\ZZ,1)$ is an Eilenberg-MacLane space representing $1$-dimensional
cohomology with integral coefficients, see eg.~\cite[V.7.5 and 7.14]{Whitehead}, we have for every
connected Lie group $H$ a natural map 
\[
\Hom(H,\SS^1){\rTo}[H,\SS^1]_0\cong H^1(H).
\]
For the torus $H=K/[K,K]$ this map is an isomorphism, see \cite[8.57(ii)]{HMCompact}.
\qed
\end{Lem}
\begin{Cor}
\label{pi_1_vs_abelianization}
Let $\phi:H\rTo K$ be a continuous homomorphism between compact connected Lie
groups. If $\phi_*:\pi_1(H)\rTo\pi_1(K)$ is bijective (surjective), then
the abelianization $H/[H,H]\rTo K/[K,K]$ is bijective (surjective).

\proof
If $\phi_*$ is bijective/surjective on the fundamental groups,
then the map in $1$-dimensional cohomology is bijective/injective
by duality and \ref{CofunctorLemma}. Thus
\[
\Hom(H/[H,H],\SS^1)\lTo \Hom(K/[K,K],\SS^1)
\]
is bijective/injective.
Dualizing again, we obtain the claim by (Pontrjagin) duality.
\qed
\end{Cor}
We now apply this result to $\cG$ in order to control the torus factors. 
Recall that $G$ is semisimple if $m,n\geq 2$.
\begin{Prop}
\label{StabilizersAreSemisimple}
If $m=2<n$, then $G_p$ is semisimple and
$G_q$ and $G_d$ have $1$-dimensional centers,
and $G_{p,d,q}$ has a $2$-dimensional center.
If $m,n\geq 3$, then all groups appearing in
$\cG$ are semisimple.

\proof
A compact connected Lie group is perfect if and only if
it is semisimple, see \cite[6.16]{HMCompact}.
By \ref{IsotropyConnectedLemma} all groups $G_\alpha$ are connected.
We consider the abelianizations $H_\alpha=G_\alpha/[G_\alpha,G_\alpha]$.
Let $\mathcal H$ denote the diagram formed by these seven
abelian groups $H_\alpha$. Suppose that this diagram
has a continuous homomorphism to some abelian topological group $H$. By the
universal property of $G$, there is a unique homomorphism
$G\rTo H$ commuting with the maps $G_\alpha\rTo H_\alpha\rTo H$.
Since $G$ is perfect, each composite map $G_\alpha\rTo G\rTo H$ is constant.
It follows that the seven maps $H_\alpha\rTo H$ are also
constant.

If $m,n\geq 3$, then all maps in $\mathcal H$ are isomorphisms
by \ref{pi_1_vs_abelianization}. From the previous paragraph
we conclude that all groups in $\mathcal H$ are trivial.
If $m=2<n$, then all groups in $\mathcal H$ surject naturally
onto $H_p$. Again by the previous paragraph, $H_p=1$.
For $\alpha=\{p,d\},\{p,q\}$ we have 
$G_p/G_\alpha\cong\CCP^2$ and $G_\alpha/G_{p,d,q}\cong\SS^2$ 
and therefore short exact sequences
\begin{diagram}
0&\rTo&\ZZ&\rTo&\pi_1(G_\alpha)&\rTo&\pi_1(G_p)&\rTo&0\\
0&\rTo&\ZZ&\rTo&\pi_1(G_{p,d,q})&\rTo&\pi_1(G_\alpha)&\rTo&0 .
\end{diagram}
Thus $\dim H_\alpha=1$ and $\dim H_{p,d,q}=2$ by \ref{CofunctorLemma}.
\qed
\end{Prop}
\begin{Num}\textbf{\boldmath The Lie algebra diagram $\Lie(\cG)$\ }
\label{TheLieAlgebraDiagram}
Passing to the Lie algebras of the groups in $\cG$, we obtain
a commutative diagram of Lie algebra inclusions which we denote
by $\Lie(\cG)$.
The next proposition reduces in many cases the classification of
the possible complexes $\cG$ to the much simpler classification of
the complexes of Lie algebras $\Lie(\cG)$. 
For $\emptyset\neq\alpha\subseteq\gamma$, we denote by 
$\widetilde{G_\alpha}$ the simply connected group with Lie algebra
$\Lie(G_\alpha)$. 
In this way we obtain from $\Lie(\cG)$ a commutative diagram
of simply connected Lie groups which we denote by $\widetilde\cG$.
We note that $\Lie(\cG)$ and $\widetilde\cG$ encode exactly
the same information, see eg.~\cite[5.42 and A2.26]{HMCompact}.
\end{Num}
The group $\widetilde{G_\alpha}$ is the universal
covering of $(G_\alpha)^\circ$ and we have a central extension
\[
1\rTo\pi_1(G_\alpha)\rTo\widetilde{G_\alpha}\rTo (G_\alpha)^\circ\rTo 1.
\]
The identification of $\pi_1(G_\alpha)$ with the kernel of this
map is compatible with the maps on the fundamental groups in $\cG$.
\begin{Prop}
\label{LieSuffices}
If $m,n\geq 2$, then $\cG$ is uniquely determined by $\Lie(\cG)$.

\proof
We begin with a small observation. Let $z\in G_{p,d,q}$.
If $z$ is central in $G_p$ and in $G_q$,
then $z\in \Cen(G)$, because $G=\bra{G_p\cup G_q}$.
It follows that $z=1$, since $G$ acts faithfully.

By Lemma \ref{IsotropyConnectedLemma}, all groups $G_\alpha$ in
$\cG$ are connected. Therefore $\widetilde\cG$ consists of
the universal coverings of the $G_\alpha$. 
We let $\pi_1\cong\pi_1(G_{p,d,q})$ denote the kernel of
the map $\widetilde{G_{p,d,q}}\rTo G_{p,d,q}$.
From \ref{ConnDiagram} we see that for each 
$\emptyset\neq\alpha\subseteq\gamma$, the group $\pi_1$
maps onto the kernel of $\widetilde{G_\alpha}\rTo G_\alpha$.

The group $\pi_1$ can now be characterized as follows. It consists
of all elements $z\in \widetilde{G_{p,d,q}}$ whose images are
central in each $\widetilde{G_\alpha}$. Indeed, every $z\in\pi_1$
has this property. Conversely, if $z\in \widetilde{G_{p,d,q}}$ 
has this property,
then its image in every $G_\alpha$ is central and thus 
its image in $G_{p,d,q}$ is trivial by the
small observation above. Thus $\pi_1$ is determined by
$\widetilde\cG$. It follows that $\cG$ is determined
by $\Lie(\cG)$.
\qed
\end{Prop}
\begin{Num}\textbf{Kernels\ }
\label{KernelsForC3}
We introduce some more notation. 
We denote by $A$, $B$ and $C$ the kernels
of the actions of $G_p$, $G_q$ and $G_d$ on
$\lk(p)$, $\lk(q)$ and $\lk(d)$, respectively.
Their respective Lie algebras are denoted by
$\fa$, $\fb$ and $\fc$. We choose supplements
$\fp$, $\fd$ and $\fq$, such that
\[
\begin{array}{c@{\hspace{1cm}}c@{\hspace{1cm}}c}
\Lie(G_p)=\fg_p=\fp\oplus\fa &
\Lie(G_q)=\fg_q=\fq\oplus\fb &
\Lie(G_d)=\fg_d=\fd\oplus\fc \\ \\
\Lie(G_p/A)\cong\fp          &
\Lie(G_q/B)\cong\fq          &
\Lie(G_d/C)\cong\fd   ,       
\end{array}
\]
see eg.~\cite[5.78]{HMCompact}.
By \ref{EffectiveOnE_1} we have
\[
 A\cap B=B\cap C=C\cap A=1\qquad\text{and}\qquad
 \fa\cap\fb=\fb\cap\fc=\fc\cap\fa=0.
\]
Moreover, $A$ is contained in $G_{p,d,q}$ and acts 
by \ref{EffectiveOnE_1} 
faithfully and as a subgroup of $\mathrm{O}(n)$ on 
$|\lk(\{d,q\})|_K\cong\SS^n$. Similarly, 
$B$ acts faithfully as a subgroup of $\mathrm{O}(m)$ on
$|\lk(\{p,d\})|_K\cong\SS^m$, and $C$ acts faithfully as a 
subgroup of $\mathrm{O}(m)$ on $|\lk(\{p,q\})|_K\cong\SS^m$.

\end{Num}
\begin{Lem}
\label{AIsTrivial}
If $m>n=1$, then $A=1$.

\proof
By \ref{IsotropyConnectedLemma}, the group $G_{d,q}$ is connected.
Therefore it induces the group $\SO(2)$ on the $1$-sphere $\lk(\{d,q\})$.
The group $G_{p,d,q}$ acts thus trivially on $\lk(\{d,q\})$.
In particular, $A$ acts trivially on $\mathcal E_1(\{p,d,q\})$,
hence $A=1$ by \ref{EffectiveOnE_1}.
\qed
\end{Lem}
\begin{Lem}
\label{LemmaForm>n=1}
Suppose that $m>n=1$ and that $G_p$ is connected.
Then $\cG$ is determined by $\Lie(\cG)$ and the subdiagram
\begin{diagram}[height=1.5em]
G_{p,d} & \lTo & G_{p,d,q} & \rTo & G_{p,q}\\
\dTo &&& \ldTo(4,2)\\
G_p
\end{diagram}

\proof
The proof is similar to the proof of \ref{LieSuffices} above.
From our assumptions and \ref{IsotropyConnectedLemma} we see that
all seven groups in $\cG$ are connected. 
We define the diagram $\widetilde\cG$ as in \ref{TheLieAlgebraDiagram}.
The groups $\widetilde{G_\alpha}$ are thus the universal covering groups
of the $G_\alpha$. In $\widetilde\cG$ we consider the two maps
$\widetilde{G_d}\lTo^\phi\widetilde{G_{d,q}}\rTo^\psi\widetilde{G_q}$.
Since both $G_q/G_{d,q}$ and $G_d/G_{d,q}$ are $1$-connected, 
we have 
\[
\widetilde{G_d}/\phi(\widetilde{G_{d,q}})=G_d/G_{d,q}\qquad\text{and}\qquad
\widetilde{G_q}/\psi(\widetilde{G_{d,q}})=G_q/G_{d,q}.
\]
An element $z\in\widetilde{G_{d,q}}$ which acts trivially
on $G_q/G_{d,q}$ acts trivially on $\lk(q)$. If it acts in addition trivially
on $G_d/G_{d,q}$, then it acts trivially on $G_{p,d}/G_{p,d,q}$ and hence on
$\mathcal E_1(\{p,d,q\})$. By \ref{EffectiveOnE_1}, it acts then trivially
on $\Delta$.
Let $\pi_1\subseteq \widetilde{G_{d,q}}$
be the subgroup consisting of these elements. 
Then
$\pi_1$ is precisely the kernel of the map $\widetilde{G_{d,q}}\rTo G_{d,q}$,
i.e.~$\pi_1=\pi_1(G_{d,q})$. By \ref{ConnDiagram}, the group $\pi_1$
maps onto $\pi_1(G_d)$ and onto $\pi_1(G_q)$. Therefore $\widetilde\cG$
determines the diagram $G_d\lTo G_{d,q}\rTo G_q$ completely.
Since all groups in $\cG$ are connected, 
the maps in $\widetilde\cG$ determine also the maps in $\cG$.
\qed
\end{Lem}
For $m=1$ we have to deal with stabilizers that are not connected.
The identity components of the stabilizers form a subdiagram
of $\cG$ which we denote by $\cG^\circ$.
\begin{Lem}
\label{LemmaForn>m=1}
Suppose that $n>m=1$.
Then $\cG^\circ$ is determined by $\Lie(\cG)$ and the subdiagram
\begin{diagram}[height=1.5em]
G_{p,d} & \lTo & G_{p,d,q} & \rTo & G_{p,q}\\
\dTo &&& \ldTo(4,2)\\
G_p
\end{diagram}

\proof
We argue similarly as in the proof of \ref{LemmaForm>n=1}.
For the four simplices $\alpha$ with 
$p\in\alpha\subseteq\{p,d,q\}$, we know already the 
kernels $\pi_1(G_\alpha)$ of the central extensions
\[
 \pi_1(G_\alpha)\rTo\widetilde{G_\alpha}
 \rTo (G_\alpha)^\circ,
\]
since we know the groups $G_\alpha$.
For $\emptyset\neq \beta=\alpha\setminus\{p\}$ the homomorphism 
$\pi_1(G_\alpha)\rTo\pi_1(G_\beta)$ is onto by \ref{ConnDiagram},
and this homomorphism, which is the restriction of
$\widetilde{G_\alpha}\rTo\widetilde{G_\beta}$,
is in turn determined by the homomorphism
$\Lie(G_\alpha)\rTo\Lie(G_\beta)$.
Therefore we know also the groups $(G_\beta)^\circ$, and, since they
are connected, the maps between them.
\qed
\end{Lem}
The problem is then to pass from $\cG^\circ$ to $\cG$. 
This requires the following homological fact, which allows us to determine 
$G_\alpha$ once we know $(G_\alpha)^\circ$ and $G_\alpha/N$, where 
$N$ is the kernel of $G_\alpha$ on $\lk(\alpha)$.
See also Hilgert-Neeb \cite[18.2]{HilgertNeeb} 
for a slightly more special result.
\begin{Lem}
\label{HomologyLemma}
Suppose $F$ and $H$ are Lie groups, that $F$ is connected and
that $F\rTo^p H$ is an open and continuous homomorphism with discrete
kernel $D$. Consider the category $\mathcal C$ of all Lie group homomorphisms
$F\rTo^i E\rTo^q H$, where $q$ is a surjective covering map and $i$
is an open inclusion, such that $q\circ i=p$.
\begin{diagram}[height=1.5em]
F & \rDotsto^i &E\\
\dTo^p&&\dDotsto_q\\
H^\circ&\rInto& H
\end{diagram}
If the category $\mathcal C$ is nonempty, then its isomorphism classes
are parametrized by the cohomology group $H^2(\pi_0(H),D)$.
The group $\pi_0(H)$ acts on $D$ by conjugation, and the cohomology is
taken with respect to this action.

\proof
Suppose that $F\rTo^i E\rTo^q H$ is in $\mathcal C$.
We view $\pi_0(H)$ as the group of path components of $H$.
If $X$ is a path component of $H$, then $E_X\rTo X$ is, as a bundle map,
isomorphic to the bundle map $F\rTo H^\circ$. Hence every
other Lie group solution $F\rTo E'\rTo H$ is isomorphic to one 
living on the same covering space $E$, but
possibly with a different group multiplication. 

We denote the given 
multiplication on $E$ by a dot $\cdot$ and we assume that
$*:E\times E\rTo E$ is another Lie group multiplication on $E$,
compatible with $q$. Suppose that $X,Y,Z\in\pi_0(H)$ are path
components with $XY=Z$. Then $E_X*E_Y=E_Z=E_X\cdot E_Y$,
and for all $x\in E_X$, $y\in E_Y$ we have $q(x*y)=q(x\cdot y)$.
The map 
\[
c:E\times E\rTo D,\qquad
(x,y)\mapstoo (x*y)\cdot(x\cdot y)^{-1}
\]
is locally
constant and factors to a map 
\[
c:\pi_0(H)\times\pi_0(H)\rTo D.
\]
The associativity of $*$ implies the cocycle condition. More precisely, we
have the identities
\[
 c(U,VW)\cdot Uc(V,W)U^{-1}=c(UV,W)\cdot c(U,V) \qquad\text{and}\qquad c(H^\circ,U)=c(U,H^\circ)=1
\]
which say that $c$ is a normalized $2$-cocycle.
Conversely, if $c$ is a locally constant map that satisfies these two
conditions, then $x*y=c(x,y)\cdot x\cdot y$ defines a new Lie group
multiplication on $E$, as is easily checked. Finally, we have the group
of deck transformations acting on these multiplications.
These maps yield precisely the coboundaries, and the claim follows.
\qed
\end{Lem}
Note that the proof gives a method to construct all other multiplications from
a given one.
The case which is interesting for us is when $\pi_0(H)\cong\ZZ/2\cong D$.
Then the action of $\pi_0(H)$ on $D$ is trivial and we have
\[
H^2(\pi_0(H);D)=H^2(\RRP^\infty;\ZZ/2)\cong\ZZ/2.
\]
Hence there are two multiplications in this case.
\begin{Example}
\label{FakeO2}
Let $F=\SO(2)$, $D=\{\pm1\}$ and $H=\mathrm{O}(2)/D$.
There are two Lie groups $E$ which fit into the diagram
\begin{diagram}[height=2em]
\SO(2)&\rTo&E\\
\dTo&&\dTo\\
\SO(2)/D&\rTo&\mathrm{O}(2)/D.
\end{diagram}
One group is $E=\mathrm{O}(2)$, the other group is the 'fake $\mathrm{O}(2)$',
which is 
$E'=\mathrm  U(1)\cup\bj\mathrm U(1)\subseteq\Sp(1)$.
For all $g\in \mathrm{O}(2)\setminus\SO(2)$ we have $g^2=1$, whereas
$g^2=-1$ holds for all $g\in E'\setminus\mathrm U(1)$.
The group $E'$ is formally obtained from $\mathrm O(2)$ by putting
\[
 u*v=\begin{cases}
     \phantom{-} uv &\text{ if }(\det(u),\det(v))\neq (-1,-1)\\
      -uv &\text{ if }(\det(u),\det(v))=(-1,-1).
     \end{cases}
\]
\end{Example}
Now we start the actual classification. As we noted above, we can identify
$\fq$ with $\Lie(G_q/B)$ etc. In this way we obtain three diagrams for
the Lie algebras of the groups acting faithfully on the links.
\begin{diagram}[size=2em]
\fd &&& \lTo && \fd_q && \fq_d &&& \rTo &&&&&& \fq \\
\uTo &&&&& \uTo &&  \uTo &&&&&&&& \ruTo(4,2) \\
\fd_p &&& \lTo && \fd_{p,q} && \fq_{p,d} &&&\rTo && \fq_p\\ 
\fp_d &&& \lTo &&& \fp_{d,q} &&&\rTo &&& \fp_q\\
\dTo &&&&&&&&&&& 
\ldTo(12,2)\\
\fp
\end{diagram}
The groups belonging to these Lie algebras are known by
\cite{GKK1,GKK2}.
From this, we determine $\Lie(\cG)$ in the following way.
We first determine the possible isomorphisms
\begin{diagram}[size=2em]
\fq_{p,d}\oplus\fb&&&\rTo&&&\fq_p\oplus\fb\\
&&& \uTo^\cong_\iota \\
\fp_{d,q}\oplus\fa&&&\rTo&&&\fp_q\oplus\fa .
\end{diagram}
Once this is done, it turns out in each case that there is
just one possibility for the structure of $\fd$, and
one possibility to fill in $\fg_d$.
These data determine $\Lie(\cG)$. If $m,n\geq2$, this
determines $\cG$. In the cases where $1\in\{m,n\}$, further work is needed.

We now consider the cases $m=1,2,4,8$ separately.
Accordingly, $(G_p/A)^\circ$ is one of the groups
\[
\SO(3),\quad\mathrm{PSU}(3),\quad\Sp(3),\ \text{ or }\ \mathrm{F}_4,
\]
see \cite[63.8]{BlauesBuch} and \cite{GKK2}.
Also, $G_p/A$ is necessarily connected, unless we are in the
case $m=2$, where $G_p/A$ may have two components.
We begin with the case $m=8$.

\subsection{\boldmath The classification of $\cG$ for $m=8$}

By \ref{BasicStrucureOfRank2Links},
$\lk(p)$ is the projective plane over the Cayley algebra.
The subalgebras $\fp_d,\fp_{d,q},\fp_q\subseteq\fp$
form the following diagram, with the standard inclusions
corresponding to the Cayley plane.
\begin{diagram}
\fso(9)&\lTo &\fso(8)
&\rTo & \fso(9)\\
\dTo&&&\ldTo(4,2)\\
\mathfrak{f}_4
\end{diagram}
According to the Main Theorem in \cite{GKK2} there is just one possibility for
$\fq$, with $n=1$.
The corresponding part of the diagram for the subalgebras
$\fq_d,\fq_p,\fq_{p,d}\subseteq\fq$ is as follows.
\begin{diagram}
 \fso(8)\oplus\RR &\rTo &&&&& \fso(10)\oplus\RR\\
 \uTo &&&&& \ruTo(4,2) \\
\fso(8)&\rTo& \fso(9)
\end{diagram}
From \ref{KernelsForC3} we see that $\fa=\fb=0$.
Up to inner automorphism, there is a unique possibility for the
isomorphism $\iota$.
By \ref{IsotropyConnectedLemma}, the group $G_d$ is connected, and by
\ref{StructureOfD} it induces $\SO(9)\times\SO(2)$ on $\lk(d)$, in its
standard action on $\SS^8\times\SS^1$. Up to inner automorphism,
the isomorphisms which identify $\fq_{p,d}\rTo\fq_d$
and $\fd_{p,q}\rTo\fd_q$ are parametrized by a nonzero real number
$r$ (acting on the $\RR$-summand). This determines the diagram
$\Lie(\cG)$. The diagrams for different values of $r$ are
isomorphic. Thus, there is a unique possibility for
$\Lie(\cG)$.

We have $A=1$ by \ref{AIsTrivial}. All automorphisms of $\mathfrak f_4$
are inner and the corresponding compact Lie group $\mathrm F_4$ is both centerless and simply connected,
see eg.~\cite[94.33]{BlauesBuch}. In particular, $G_p\cong\mathrm F_4$ is connected and we have
determined the subdiagram
\begin{diagram}[height=2em]
G_{p,d} & \lTo & G_{p,d,q} & \rTo & G_{d,q}\\
\dTo && & \ldTo(4,2)\\
G_p 
\end{diagram}
By \ref{LemmaForm>n=1}, this diagram together with $\Lie(\cG)$ determines $\cG$ uniquely.
\begin{Prop}
If $m=8$, then there is up to isomorphism a unique possibility for the diagram $\cG$.
\qed
\end{Prop}
This unique possibility for $\cG$ is realized by the nonembeddable polar space, see
\ref{TheNonEmbeddablePolarSpace}. The minimal universal group is
$\widehat G=\mathrm E_6\cdot\SS^1$, see \cite[Ch.~X Table V]{He} and \cite{EH}.

\subsection{\boldmath The classification of $\cG$ for $m=4$}

By \ref{BasicStrucureOfRank2Links}, $\lk(p)$ is the projective plane over the quaternions.
The subalgebras $\fp_q,\fp_d,\fp_{d,q}$ form the following diagram, with the 'obvious' inclusions.
We decorate the arrows by the kernels of actions on the respective spheres.
\begin{diagram}
\fsp(1)\oplus\fsp(2)&\lTo_{\scriptstyle\fsp(1)\oplus0\oplus0} &\fsp(1)\oplus\fsp(1)
\oplus\fsp(1)&\rTo_{\scriptstyle0\oplus0\oplus\fsp(1)}& \fsp(2)\oplus\fsp(1)\\
\dTo&&&\ldTo(4,2)\\
\fsp(3)
\end{diagram}
According to the Main Theorem in \cite{GKK2} there are the following possibilities for $\fq$,
with $n\in\{1,5\}\cup(3+4\NN)$.
\begin{Num}[\boldmath$n=4\ell+3$, for $\ell\geq 2$.]
\label{4l+3}
There is one possibility for $\fq$, which is as follows.
\begin{diagram}
 \fsp(1)\oplus\fsp(1)\oplus\fsp(\ell+1) &\rTo &&&&& 
\fsp(2)\oplus\fsp(\ell+2)\\
 \uTo_{\scriptstyle\fsp(1)\oplus0\oplus0} &&&&& \ruTo(4,2) \\
\fsp(1)\oplus\fsp(1)\oplus\fsp(\ell)&\rTo_{\scriptstyle0\oplus0\oplus\fsp(\ell)}& \fsp(2)\oplus\fsp(\ell)
\end{diagram}
Thus $\fa=\fsp(\ell)$ and $\fb=\fsp(1)$. Up to inner automorphisms, there is a unique identification
between $\fp_{d,q}\oplus\fa\rTo\fp_q\oplus\fa$ and $\fq_{p,d}\oplus\fb\rTo\fq_{p}\oplus\fb$
which is compatible with the action on $\SS^4$.
From \ref{StructureOfD} we see that $G_d$ induces
the group $\SO(5)\times(\Sp(\ell+1)\cdot\Sp(1))$ on $\SS^4\times\SS^{4\ell+3}$.
This determines $\fg_d\cong\fsp(2)\oplus\fsp(\ell+1)\oplus\fsp(1)$ 
and the remaining homomorphisms in $\Lie(\cG)$.
By \ref{LieSuffices}, $\cG$ is determined by $\Lie(\cG)$.
\qed
\end{Num}
\begin{Num}[\boldmath$n=7$.]
There is one possibility for $\fq$, which is as follows.
\begin{diagram}
 \fsp(1)\oplus\fsp(1)\oplus\fsp(2) &\rTo &&&&& 
\fsp(2)\oplus\fsp(3)\\
 \uTo_{\scriptstyle\fsp(1)\oplus0\oplus0} &&&&& \ruTo(4,2) \\
\fsp(1)\oplus\fsp(1)\oplus\fsp(1)&\rTo_{\scriptstyle0\oplus0\oplus\fsp(1)}& \fsp(2)\oplus\fsp(1)
\end{diagram}
From this and \ref{StructureOfD} 
we see that $G_d$ induces the group
$\SO(5)\times(\Sp(2)\cdot\Sp(1))$ in its standard action on $\SS^4\times\SS^7$.
In particular, $\fg_{p,d,q}$ contains 
$\fsp(1)\oplus\fsp(1)\oplus\fsp(1)\oplus\fsp(1)$.
From this we see that $\fa=\fsp(1)$ and $\fb=\fsp(1)$. The isomorphism
$\iota$ is unique up to inner automorphisms.
We end up with a unique diagram $\Lie(\cG)$ as in \ref{4l+3}, 
with $\ell=1$.
This diagram determines $\cG$ by \ref{LieSuffices}.
\qed
\end{Num}
\begin{Num}[\boldmath$n=3$.] 
There is one possibility for $\fq$, which is as follows.
\begin{diagram}
 \fsp(1)\oplus\fsp(1)\oplus\fsp(1) &\rTo &&&&& 
\fsp(2)\oplus\fsp(2)\\
 \uTo_{\scriptstyle\fsp(1)\oplus0}   &&&&& \ruTo(4,2) \\
\fsp(1)\oplus\fsp(1)&\rTo_{\scriptstyle0\oplus0}  & \fsp(2)
\end{diagram}
From this and \ref{StructureOfD} we see that $G_d$ induces the group
$\SO(5)\times\SO(4)$ on $\SS^4\times\SS^3$. 
Thus we have $\fb=\fsp(1)$ and $\fa=0$. The isomorphism $\iota$ is unique
up to inner automorphisms and $\Lie(\cG)$ is uniquely determined.
This determines $\cG$ by \ref{LieSuffices}.
\qed
\end{Num}
\begin{Num}[\boldmath$n=5$.]
There is one generalized quadrangle, but two transitive
connected groups, of type $\fsu(5)$ and $\fsu(5)\oplus\RR$,
respectively.
By \ref{CommutatorIsTransitive} and \ref{StabilizersAreSemisimple},
the group $G_q$ is semisimple, hence the ideal $\fq$ is also 
semisimple. The possibility for $\fq$ is thus as follows.
\begin{diagram}
 \fsu(3)\oplus\fsu(2)&\rTo &&&&& 
\fsu(5)\\
 \uTo_{\scriptstyle0\oplus\fsu(2)} &&&&& \ruTo(4,2) \\
\fsu(2)\oplus\fsu(2)&\rTo_{\scriptstyle0\oplus0}& \fso(5)
\end{diagram}
From this and \ref{StructureOfD} we see that the group induced
by $G_d$ on $\SS^4\times\SS^5$ is $\SO(5)\times\SU(3)$.
It follows that $\fb=\fsp(1)$ and $\fa=0$. The isomorphism
$\iota$ is unique up to inner automorphisms and the diagram
$\Lie(\cG)$ is uniquely determined. This determines
$\cG$ by \ref{LieSuffices}.
\qed
\end{Num}
\begin{Num}[\boldmath$n=1$.] 
There is a unique possibility for $\fq$, which is as follows.
\begin{diagram}
\RR\oplus \fso(4)&\rTo &&&&& \RR\oplus\fso(6)\\
 \uTo_{\scriptstyle0} &&&&& \ruTo(4,2) \\
\fso(4)&\rTo_{\scriptstyle0}& \fso(5)
\end{diagram}
By \ref{IsotropyConnectedLemma}, the group $G_d$ is connected.
From this and \ref{StructureOfD} we see that the group induced
by $G_d$ on $\SS^4\times\SS^1$ is $\SO(5)\times\SO(2)$.
Also, we have $\fb=\fsp(1)$ and $\fa=0$. 
The isomorphism $\iota$ is unique up to
inner automorphisms and the diagram $\Lie(\cG)$ is uniquely
determined by this. We have $A=1$ by \ref{AIsTrivial}, hence
$G_p=\mathrm{PSp}(3)$. From 
\ref{LemmaForm>n=1} we see that $\cG$ is uniquely determined.
\qed
\end{Num}
These are all possibilities for $m=4$. In each case, there exists a
building $\Delta$ corresponding to $\cG$. The possibilities for
$G$ are given by \cite[Ch.~X Table V]{He} and \cite{EH}. 
They are as follows.
\begin{Num}
If $n=4\ell+3$, with $\ell\geq 0$, then $\Delta$ is the polar space associated to
the quaternionic $(1,[a\mapstoo\bar{a}])$-hermitian form 
\[
h=(-f_3)\oplus f_{3+\ell}
\]
on $\HH^{3+(3+\ell)}$. In this case $G=\bigl(\Sp(3)\times\Sp(3+\ell)\bigr)/\bra{(-1,-1)}$.

If $n=1,5$, then $\Delta$
is the polar space associated to
the (unique) quaternionic $(-1,[a\mapstoo\bar{a}])$-hermitian form 
on $\HH^6$ or $\HH^7$.
Either $G=\mathrm{U}(6)/\bra{-1}$, with $n=1$, or
$G=\SU(7),\mathrm{U}(7)/\bra{-1}$, with $n=5$.
\end{Num}

\subsection{\boldmath The classification of $\cG$ for $m=2$}

By \ref{BasicStrucureOfRank2Links}, $\lk(p)$ is the projective plane over $\CC$.
The subalgebras $\fp_q,\fp_d,\fp_{d,q}$ form the following diagram, with the 'obvious' inclusions.
We decorate the arrows by the kernels of actions on the respective spheres.
\begin{diagram}
\RR\oplus\fsu(2)&\lTo_{\scriptstyle\RR\oplus0} &\RR\oplus\RR
&\rTo_{\scriptstyle0\oplus\RR}& \fsu(2)\oplus\RR\\
\dTo&&&\ldTo(4,2)\\
\fsu(3)
\end{diagram}
This case $m=2$ is more complicated since we have to deal with reductive Lie
algebras, where the complement of an ideal is not necessarily unique.
We fix some more notation. We identify the Lie algebra
$\fsu(3)$ with the algebra of
complex traceless skew-hermitian $3\times 3$ matrices, and the upper line
in the diagram above with the following inclusions in $\fsu(3)$.
\begin{diagram}
\fp_d&\lTo&\fp_{d,q}&\rTo&\fp_q\\
\left(\begin{smallmatrix}
* &   &   \\ 
  & * & * \\
  & * & *
\end{smallmatrix}\right)
&\lTo&
\left(\begin{smallmatrix}
* &   &   \\ 
  & * &   \\
  &   & *
\end{smallmatrix}\right)
&\rTo&
\left(\begin{smallmatrix}
* & * &   \\ 
* & * &   \\
  &   & *
\end{smallmatrix}\right)
\end{diagram}
We have the following four $1$-dimensional subalgebras of $\fp_{d,q}\cong\RR^2$. Each pair of
them spans $\fp_{d,q}$, and we use them below.
\begin{align*}
\fz_d&=\Cen(\fp_d)=\textstyle\left\{\left.\left(\begin{smallmatrix}  -2s\bi&&\\&s\bi\\&&s\bi  \end{smallmatrix}\right)\right| s\in\RR\right\}
  & \ft_d&=\fp_{d,q}\cap[\fp_d,\fp_d]=\left\{\left.\left(\begin{smallmatrix}  0&&\\&s\bi\\&&-s\bi  \end{smallmatrix}\right)\right| s\in\RR\right\} \\
\fz_q&=\Cen(\fp_q)=\left\{\left.\left(\begin{smallmatrix}  s\bi&&\\&s\bi\\&&-2s\bi  \end{smallmatrix}\right)\right| s\in\RR\right\} 
  & \ft_q&=\fp_{d,q}\cap[\fp_q,\fp_q]=\left\{\left.\left(\begin{smallmatrix}  s\bi&&\\&-s\bi\\&&0  \end{smallmatrix}\right)\right| s\in\RR\right\}
\end{align*}
According to \cite{GKK1,GKK2}, we have
$n\in\{2\}\cup(2\NN+1)$, and
there are the following possibilities for~$\fq$.
\begin{Num}[\boldmath$n=2\ell+1$ and $\ell\geq 2$.]
\label{2l+1}
By \ref{StabilizersAreSemisimple}, the group $G_p$ is semisimple and
$G_d$ and $G_q$ have $1$-dimensional centers.
The Lie algebra $\fa$ is then also semisimple.
We let $L=G_q/B$ denote the group induced by $G_q$ on $\lk(q)$, with
$\fq\cong\Lie(L)$.
From the Main Theorem
in \cite{GKK2} we see that there are two possibilities for $L$, both
acting on the same generalized quadrangle. These actions
can be understood from the two orbit equivalent polar representations of
$\SU(2)\times\SU(\ell+2)$ and $\mathrm U(2)\times\SU(\ell+2)$
on $\CC^{2\times(\ell+2)}$, as described in \ref{Veronesean}.
The  semisimple commutator group $[L,L]$ acts transitively on 
$\lk(q)$. We denote its Lie algebra by $\fq'=[\fq,\fq]$. The diagram for $\fq'$ looks
as follows.
\begin{diagram}
\RR\oplus \fsu(\ell+1)&\rTo &&&&& \fsu(2)\oplus\fsu(\ell+2)\\
 \uTo_{\scriptstyle0\oplus0} &&&&& \ruTo(4,2) \\
\RR\oplus\fsu(\ell)&\rTo_{\scriptstyle0\oplus\fsu(\ell)}& \fsu(2)\oplus\fsu(\ell)
\end{diagram}
If $\fq$ is not semisimple (and we will see shortly that this is indeed the case),
then $\fq_\alpha=\fq'_\alpha\oplus\RR$ 
in the diagram above. We note also that $\fsu(\ell)\subseteq\fa$.
From the polar representation we see that the group induced by $[L,L]_d$ on 
$\lk(\{d,q\})\cong\SS^{2\ell+1}$ is $\mathrm{U}(\ell+1)$.
From this and \ref{StructureOfD} we see that the group induced by $G_d$ on 
$\lk(d)$ is $\SO(3)\times\mathrm{U}(\ell+1)$, in its natural action on
$\SS^2\times\SS^{2\ell+1}$.
Now we determine the isomorphism $\iota$
\begin{diagram}[size=2em]
\fq_{p,d}\oplus\fb&&&\rTo&&&\fq_p\oplus\fb\\
&&& \uTo^\cong_\iota \\
\ft_d\oplus\ft_q\oplus\fa&&&\rTo&&&\fz_q\oplus[\fp_q,\fp_q]\oplus\fa .
\end{diagram}
Since $\fg_p$ is semisimple, we have
$\fa=\fsu(\ell)$. The pair $\ft_q\subseteq [\fp_q,\fp_q]\cong\fsu(2)$ is identified
with the pair $\RR\subseteq \fsu(2)\subseteq\fq'_p$, and $\iota$ is unique
up to inner automorphisms on this part. The group corresponding to the 
algebra $\ft_d$ acts trivially on $|\lk(\{d,q\})|_K\cong\SS^n$,
because we have a product action
on $\lk(d)$. It acts, however, nontrivially on $|\lk(\{p,q\})|_K\cong\SS^m$. 
There is a unique
homomorphism $\ft_d\rTo\fq$ corresponding to such an action. It follows
that $\fq=\fq'\oplus\RR$ is not semisimple, that $\fb=0$, and now we have determined
the isomorphism $\iota:\fp_d\oplus\fa\rTo\fq_p$.  The structure of $\lk(d)$
was already determined above. 
Thus $\Lie(\cG)$ is uniquely determined, and so is $\cG$ by \ref{LieSuffices}.
\qed
\end{Num}
Now we get to the interesting case where the exceptional geometry occurs.
\begin{Num}[\boldmath$n=3$.]
\label{2*1+1}
By \ref{StabilizersAreSemisimple}, the group
$G_p$ is semisimple, and so is $\fa$.
From \ref{ConnDiagram}
and \ref{pi_1_vs_abelianization} we see again that $\fg_q$ has 
a $1$-dimensional center.
We use the same notation as in \ref{2l+1}. 
The diagram for $\fq'$ is as follows.
\begin{diagram}
\RR\oplus \fsu(2)&\rTo &&&&& \fsu(2)\oplus\fsu(3)\\
 \uTo &&&&& \ruTo(4,2) \\
\RR&\rTo& \fsu(2)
\end{diagram}
If $\fq$ is not semisimple, then we have again
$\fq_\alpha=\fq'_\alpha\oplus\RR$ in the diagram above.
In either case, $\fa=0$ (because it is semisimple), hence $\fg_{p,d,q}\cong\RR^2$.
Let $L=G_q/B$. The group induced by $[L,L]_d$ on $\lk(\{d,q\})\cong\SS^3$ is
$\mathrm U(2)$. The isomorphism $\iota$ identifies the pair $\ft_q{\rTo} [\fp_q,\fp_q]$ with
the pair $\RR\rTo\fsu(2)$ in the diagram above, uniquely up to inner automorphisms.
So far, everything is completely analogous to \ref{2l+1}.

By \ref{StructureOfD}, we have two possibilities for the
group induced on $\lk(d)$ by $G_d$. It is either the product action of $\SO(3)\times\mathrm U(2)$
on $\SS^2\times\SS^3$ or the exceptional action of $\SO(4)$ on $\SS^2\times\SS^3$
described after \ref{StructureOfD}.

(1) Assume that we are in the case of the product action. Then $\ft_d$ acts 
trivially on $\lk(\{d,q\})$,
but nontrivially on $\lk(\{p,q\})$. As in the case $\ell\geq 1$ before,
there is a unique homomorphism $\ft_d\rTo\fq$ corresponding to this action,
and we find that $\fq=\fq'\oplus\RR$ is not semisimple.
Thus $\iota$ is uniquely determined on $\fp_d=[\fp_q,\fp_q]+\ft_d$, and $\fb=0$.
This determines also $\fg_d$ and thus $\Lie(\cG)$. By \ref{LieSuffices},
the complex $\cG$ is uniquely determined. This case corresponds
to the building.

(2) Suppose that $G_d$ induces $\SO(4)$ in the exceptional action on
$\SS^2\times\SS^3$. 
Let $C$ denote the kernel of the action of $G_d$ on $\lk(d)$, with Lie algebra
$\mathfrak c$. 
We have $\fg_d=\mathfrak{c}\oplus\mathfrak{d}$, and
$\mathfrak d\cong\fso(4)$.
Since $\ft_d\subseteq[\fp_d,\fp_d]$, we see that 
$\mathfrak d_{p,q}=\ft_d$. The group $C\subseteq G_{p,d,q}$,
on the other hand, acts trivially on $\lk(\{p,d\})$. There is a unique
subalgebra in $\fp_{d,q}=\fg_{p,d,q}$ with this property, namely $\fz_d$.
Thus $\fz_d=\mathfrak c$ acts trivially on $\lk(\{d,q\})$.
This determines a unique homomorphism $\fz_d\rTo\fq$. Thus
$\iota:\fp_q\rTo\fq_p$ is uniquely determined. Also,
$\fg_d$ is now uniquely determined, hence the same is true for 
$\Lie(\cG)$ and, by \ref{LieSuffices}, for $\cG$.
This case does not correspond to a building, but to the
exceptional polar action of $\mathrm{PSU}(3)\times\SU(3)$ on the Cayley plane,
as described in Section~\ref{TheExceptionalC3Geometry}.
\qed
\end{Num}

\begin{Num}[\boldmath$n=2$.]
By \ref{CommutatorIsTransitive}, the group $G$ is semisimple.
By \cite{GKK2},
there are two non-isomorphic possibilities for $\fq\cong\fsp(2)\cong\fso(5)$, both of which
are given by the following diagram, with different homomorphisms.
One arrow corresponds to the natural inclusion
$\fu(2)\subseteq\fso(5)$ (or $\fsp(1)\oplus\fu(1)\subseteq\fsp(2)$), 
the other to the natural inclusion
$\fso(2)\oplus\fso(3)\subseteq\fso(5)$
(or $\fu(2)\subseteq\fsp(2)$).
\begin{diagram}
\RR\oplus \fsu(2)&\rTo &&&&& \fsp(2)\\
 \uTo_{\scriptstyle\RR\oplus0} &&&&& \ruTo(4,2) \\
\RR\oplus\RR&\rTo_{\scriptstyle0\oplus\RR}& \fsu(2)\oplus\RR
\end{diagram}
From this diagram we see that either $\fa=0=\fb$ or $\fa\cong\RR\cong\fb$.
In particular, $2\leq \dim(G_{p,d,q})\leq 3$, and thus
$\dim(G)\leq 21$. Since $\fp$ is simple,
there exists a simple factor $\fh$ of $\fg$ such that
the canonical projection $\pr_\fh:\fg\rTo\fh$ is injective on
$\fp$. Since $\fq$ is also simple and $\fp\cap\fq\neq 0$,
the map $\pr_\fh$ is also injective on $\fq$. Thus $\fh$ is a
simple compact Lie algebra which contains copies of $\fsu(3)$
and $\fso(5)$, and with $\dim\fh\leq 21$. From the list
of low-dimensional simple compact Lie algebras \cite[94.33]{BlauesBuch}
and the low-dimensional representations of $\fsu(3)$ and $\fso(5)$,
see eg.~\cite[95.10]{BlauesBuch} and \cite[Ch.~4]{KramerHabil}, 
we see readily that $\fh\in\{\fsu(4),\fso(7),\fsp(3)\}$.
We consider these three cases separately.

\medskip\noindent
{\bf The case \boldmath$\fh=\fsu(4)$ is not possible}
Suppose to the contrary that $\fh\cong\fsu(4)$.
We consider the natural representation $\CC^4$ of $\fsu(4)$.
All copies of $\pr_\fh(\fp)\cong\fsu(3)$ in $\fsu(4)$ are conjugate and fix 
in this $4$-dimensional representation a $1$-dimensional complex subspace
pointwise.
Similarly, all copies of $\pr_\fh(\fq)\cong\fsp(2)$ in $\fsu(4)$ are conjugate,
with trivial $\fsu(4)$-centralizers. In particular, $\pr_\fh(\fb)=0$ and 
therefore $\pr_\fh(\fg_{p,q})\cong\fsu(2)\oplus\RR$.
Since $\fp_q\cong\fu(2)$, we see that 
$\pr_\fh(\fg_{p,q})=\pr_\fh(\fp_q)\subseteq\pr_\fh(\fp)$.
Thus $\pr_\fh(\fg_{p,q})$ fixes a $1$-dimensional subspace in $\CC^4$
pointwise. On the other hand, neither the subalgebra 
$\fu(2)\subseteq\fsp(2)\subseteq\fsu(4)$
nor the subalgebra $\fsp(1)\oplus\fu(1)\subseteq \fsp(2)\subseteq\fsu(4)$ 
fix a $1$-dimensional subspace in $\CC^4$ pointwise.
Therefore this case is not possible.

\medskip\noindent
{\bf The case \boldmath$\fh=\fso(7)$ is possible in a unique way}
We consider the standard representation $\RR^7$ of $\fso(7)$.
Since $\dim(\fso(7))=21$, we have $\fh=\fg$ and $\fa\cong\RR\cong\fb$.
The inclusion $\fg_p\cong\fu(3)\subseteq\fso(7)$
is unique up to conjugation and fixes a unique $1$-dimensional real
subspace pointwise. This determines also how $\fsu(2)\cong[\fg_{p,q},\fg_{p,q}]$
is embedded in $\fso(7)$, namely as a conjugate of its standard real 
$4$-dimensional representation.

The inclusion 
$\fg_q\cong\fso(2)\oplus\fso(5)\subseteq\fso(7)$ is also unique up
to conjugation. We fix once and for all
the standard inclusion of this algebra corresponding to the decomposition
$\RR^7=\RR^2\oplus\RR^5$, and we identify $[\fg_{p,q},\fg_{p,q}]$ 
with $\fsu(2)\subseteq\fso(5)$ acting on $\CC^2\oplus\RR=\RR^5$. 
Then $\fg_{p,q}$ has a unique real $1$-dimensional fixed space in $\RR^7$,
and this determines the subalgebra $\fg_p\cong\fu(3)$
uniquely. This shows that there is at most one possibility for 
$\Lie(\cG)$, and by \ref{LieSuffices} also for $\cG$.

\medskip\noindent
{\bf The case \boldmath$\fh=\fsp(3)$ is possible in a unique way}
We consider the standard representation $\HH^3$ of $\fsp(3)$.
Since $\dim(\fsp(3))=21$, we have $\fh=\fg$ and $\fa\cong\RR\cong\fb$.
The inclusion $\fg_p\cong\fu(3)\subseteq\fsp(3)$ is unique up to
conjugation. It corresponds to the extension of scalars given
by $\HH^3=\CC^3\otimes_\CC\HH$.
We identify $\fg_{p,q}$ with 
$\fu(1)\oplus\fu(2)$ in the standard inclusion
coming from the splitting $\HH\oplus\HH^2=(\CC\oplus\CC^2)\otimes_\CC\HH$.
The inclusion of
$\fg_q\cong\fu(1)\oplus\fsp(2)\subseteq\fsp(3)$ is also unique up to
conjugation. From the splitting of $\HH^3$ we see that there is
a unique conjugate of $\fg_q$ containing $\fg_{p,q}$.
Thus, there is at most one possibility for $\Lie(\cG)$, 
and by \ref{LieSuffices} also for $\cG$.

\medskip
Thus there are precisely two possibilities for $\cG$ with $m=n=2$.
Both are realized by polar spaces over the complex numbers, one
corresponding to the $5$-dimensional nondegenerate quadratic form
over $\CC$, and the other
to the $6$-dimensional symplectic form over $\CC$.
\qed
\end{Num}
\begin{Num}[\boldmath$n=1$.]
There is a unique possibility for $\fq$, which is as follows.
\begin{diagram}
\RR\oplus \RR&\rTo &&&&& \RR\oplus\fso(4)\\
 \uTo
&&&&& \ruTo(4,2) \\
\RR&\rTo
& \fso(3)
\end{diagram}
The groups $G_q$, $G_{d,q}$ and $G_d$ are connected by
\ref{IsotropyConnectedLemma}, and $A=1$.
By \ref{StructureOfD}, the group induced by
$G_d$ on $D$ is $\SO(3)\times\SO(2)$, in its
natural action on $\SS^2\times\SS^1$.
In particular, $G_{p,d}$ induces the group $\SO(3)$ on
$\lk(\{p,d\})$, and not the group $\mathrm O(3)$.
Therefore $G_p=\mathrm{PSU}(3)$, and all groups in
$\cG$ are connected. We now apply Lemma \ref{LemmaForm>n=1} and
conclude that $\cG$ is uniquely determined.
\qed
\end{Num}
These are all possibilities for $m=2$. 
The corresponding universal compact geometries are as follows.
\begin{Num}
If $n=2\ell+1$, with $\ell\geq 0$ and $\ell\neq 1$, then $\Delta$ is the polar space associated to
the complex $(1,[a\mapstoo\bar{a}])$-hermitian form 
\[
h=(-f_3)\oplus f_{3+\ell}
\]
on $\CC^{3+(3+\ell)}$.

If $n=3$, then either $\Delta$ is the polar space associated to
the complex $(1,[a\mapstoo\bar{a}])$-hermitian form 
\[
h=(-f_3)\oplus f_{4}
\]
on $\CC^{3+4}$, or the exceptional geometry from Section~\ref{TheExceptionalC3Geometry}.

If $n=2$, then $\Delta$ is either the polar space associated to the
complex symplectic form on $\CC^6$ or the polar space associated to the complex quadratic form on
$\CC^7$.
\end{Num}
The compact connected chamber-transitive groups $G$ on the universal geometry
$\Delta$ are as follows.
In the hermitian case we have $G=\SU(3)\cdot\mathrm{U}(3+\ell)$, or
$G=\SU(3)\cdot\SU(3+\ell)$ for $\ell\geq 1$.
In the symplectic case, the group is $G=\Sp(3)/\bra{-1}$, and 
in the orthogonal case it is $G=\SO(7)$. This follows from
\cite[Ch.~X Table V]{He} and \cite{EH}. 
In the case of the exceptional $\mathsf C_3$ geometry, $G=\mathrm{PSU}(3)\times\SU(3)$.

\subsection{\boldmath The classification of $\cG$ for $m=1$}

By \ref{BasicStrucureOfRank2Links}, the link $\lk(p)$ is the projective plane over
$\RR$. This will again be the starting point for our classification. The
Lie algebra $\fp$ is isomorphic to $\fso(3)$, and $G_p$ induces the group 
\[
K=G_p/A\cong\SO(3)
\]
on $\lk(p)$. We have $K_d\cong\mathrm O(2)\cong K_q$ and
$K_{d,q}\cong\ZZ/2\oplus\ZZ/2$. In particular, the groups $G_{p,d}$,
$G_{p,d,q}$ and $G_{p,q}$ are not connected.
We put 
\[
M=G_d/C\qquad\text{and}\qquad L=G_q/B.
\]
By \ref{StructureOfD} we have a product action of $M^\circ$ 
on $\SS^1\times\SS^n$.
The group $M$ is not connected, because $K_d$ induces the group
$\mathrm{O}(2)$ on the $1$-sphere $\lk(\{p,d\})$.
We note also that both $B$ and $C$ inject into $\mathrm{O}(1)\cong\ZZ/2$,
whence $\fb=\fc=0$.
\begin{Num}\textbf{\boldmath The structure of $\lk(q)$ and of $L=G_q/B$\ }
\label{StructureOfQForm=1}
If $n\geq 2$, then the generalized quadrangle $\lk(q)$ belongs by \cite{GKK2} to the
symmetric bilinear form $(-f_2)\oplus f_{n+2}$ on $\RR^{2+(n+2)}$. The action of 
\[
L=G_q/B
\]
is given by a polar representation of
$L\subseteq\mathrm{O}(2)\cdot\mathrm{O}(n+2)$ 
on $\RR^{2\times(n+2)}$
which is orbit equivalent to
the polar representation of $\mathrm{O}(2)\cdot\mathrm{O}(n+2)$ 
described in \ref{Veronesean}.
The dot $\cdot$ indicates that the element $(-1,-1)=(-\id_{\RR^2},-\id_{\RR^{n+2}})$
acts trivially.
By \cite{GKK2}, the identity component of $L$ is either
$\SO(2)\cdot\SO(n+2)$, or $\SO(2)\times\mathrm G_2$ for $n=5$,
or $\SO(2)\cdot\Spin(7)$ for $n=6$. These connected groups induce $\SO(2)$
on $\lk(\{p,q\})$, as one can easily check in the polar representation.
Since we know that $G_{p,q}$ induces the group $\mathrm{O}(2)$ on $\lk(\{p,q\})$, 
we see that $L$ is not connected.
We have $\fq_{p,d}=\fso(n)$, or, for the two exceptional actions,
$\fq_{p,d}=\fsu(2)$ or $\fq_{p,d}=\fsu(3)$. In any case, the Lie algebra
$\fa=\fq_{p,d}=\fg_{p,d,q}$ is semisimple for $n\geq 3$.
We will see that the group $\SO(2)\cdot\Spin(7)$ cannot occur in this setting.

If $n=1$, then there are two possibilities. Either $\lk(q)$ is the
generalized quadrangle of the symmetric bilinear form $(-f_2)\oplus f_3$
on $\RR^{2+3}$ and $L^\circ=\SO(2)\times\SO(3)$, or $\lk(d)$ is the generalized
quadrangle associated to the standard symplectic form on $\RR^4$.
The group $L^\circ$ is then $\mathrm{U}(2)/\{\pm1\}\cong\SO(2)\times\SO(3)$.
These two generalized quadrangles are dual to each other.
\end{Num}
\begin{Lem}
\label{StructureOfG_p}
We have $G_p=\SO(3)\times A$. If $n\geq 2$, then $A$ is connected.

\proof
Let $P$ be a compact connected supplement of $A^\circ$ in $(G_p)^\circ$, 
such that $(G_p)^\circ=P\cdot A^\circ$.
We claim that $P=\SU(2)$ is not possible. Assume to the contrary
that $P=\SU(2)$.
Suppose first that $n\geq 3$. We put $Z=(P_d)^\circ$. This circle 
group contains the unique nontrivial central element $z$ of $P$.
Since $\fg_{p,d,q}=\fa$ is semisimple by \ref{StructureOfQForm=1},
we have $Z=\Cen((G_{p,d})^\circ)^\circ$.
Since we have a product action of $M^\circ$ on $\lk(d)$ and since
$Z$ is connected, 
$Z$ acts trivially on $\lk(d,q)$ under the equivariant projection
$\SS^1\times\SS^n\rTo\SS^n$. In particular, $z$ acts trivially
on $\mathcal E_1(\{p,d,q\})$. This is a contradiction to 
\ref{EffectiveOnE_1}, 
hence $P=\SO(3)$ is centerless.
Suppose now that $n\leq 2$ and that $P=\SU(2)$. Then $G_{p,d,q}$ contains the 
quaternion group $Q=\{\pm1,\pm\bi,\pm\bj,\pm\bk\}$.
On the other hand, $G_{p,d,q}$ embeds into 
$\mathrm O(1)\times\mathrm O(1)\times\mathrm O(n)$. 
But this is impossible: every
$1$- or $2$-dimensional real representation of $Q$ annihilates
the element $-1$.

Thus $P=\SO(3)$ is centerless simple.
It follows that $A\cap P=1$, hence $G_p=P\ltimes A$ is a semidirect product,
and $P$ centralizes $A^\circ$. If $n\geq 2$, then $G_p$ is connected 
by \ref{IsotropyConnectedLemma},
whence $A^\circ=A$. If $n=1$, then $A$ is discrete and therefore
centralized by $P$.
\qed
\end{Lem}
As a consequence of the proof, we note that
\[
G_{p,d}\cong G_{p,q}\cong\mathrm{O}(2)\times A
\]
and that this product splitting is canonical, for all $n$.
\begin{Cor}
\label{CorForPi0AndKernels}
If $n\geq 2$, then $C\subseteq(G_d)^\circ$ and $B\subseteq(G_q)^\circ$ and
\[
\pi_0(M)\cong\pi_0(G_d)\cong\pi_0(G_{p,d})\cong\pi_0(G_{p,q})\cong\pi_0(G_q)
\cong\pi_0(L)\cong\ZZ/2.
\]

\proof
From \ref{StructureOfG_p} we know that 
$\pi_0(G_{p,d})\cong\pi_0(G_{p,q})\cong\ZZ/2$. By 
\ref{ConnDiagram}, we also have $\pi_0(G_{d})\cong\pi_0(G_{q})\cong\ZZ/2$.
The groups $M$ and $L$ are not connected, but
they cannot have more components than $G_d$ and $G_q$ have,
hence $\pi_0(M)\cong\pi_0(L)\cong\ZZ/2$.
We have $C\subseteq G_d$, and if $C\not\subseteq(G_d)^\circ$, then $M$ would be connected.
Similarly, $B$ has to be contained in $(G_q)^\circ$.
\qed
\end{Cor}
\begin{Cor}
\label{ExcludeSpin7}
The case $n=6$ with $\fq=\RR\oplus\fso(7)$ cannot occur.

\proof
Consider the $8$-dimensional real irreducible
representation of $\Spin(7)$. 
This representation is of real type, see \cite[p.~625]{BlauesBuch}.
Since the nontrivial center of $\Spin(7)$ acts faithfully on $\RR^8$,
we have $-\id_{\RR^8}\in\Spin(7)$.

Assume now to the contrary that $\fq=\RR\oplus\fso(7)$, with $n=6$.
Because $\Spin(7)$ is self-normalizing in $\mathrm O(8)$
and $L$ is not connected, we have necessarily
\[
L=(\mathrm O(2)\times\Spin(7))/\bra{(-1,-1)}.
\]
Let $\SU^{\bar\ }(3)$ denote the group generated by
$\SU(3)$ and by complex conjugation on $\CC^3$. 
From the polar representation on $\RR^{2\times 8}$ we 
see that
\[
L_p=\mathrm S(\mathrm O(2)\times\SU^{\bar\ }(3))/\bra{(-1,-1)}.
\]
But $L_p=\mathrm S(\mathrm O(2)\times\SU^{\bar\ }(3))/\bra{(-1,-1)}$ cannot be 
written as a quotient of $G_{p,q}=\mathrm O(2)\times\SU(3)$.
The reason for this is that the adjoint representation of 
$\mathrm S(\mathrm O(2)\times\SU^{\bar\ }(3))/\bra{(-1,-1)}$ 
splits off a module $\fsu(3)$ with $\SU^{\bar\ }(3)$ acting faithfully
on it, which is not the case for $\mathrm O(2)\times\SU(3)$.
This is a contradiction to \ref{StructureOfG_p}.
\qed
\end{Cor}
\begin{Lem}
\label{LemmaA}
Corresponding to each diagram
\begin{diagram}[height=1.5em]
\fq_d&\rTo&\fq\\
\uTo&&\uTo\\
\fq_{p,d}&\rTo&\fq_p
\end{diagram}
as in \ref{StructureOfQForm=1}, there is, up to isomorphism,
at most one possibility for the diagram $\Lie(\cG)$.

\proof
If $n\neq 2$, then $\fa$ is semisimple and 
$\fg_{p,d}$ and $\fg_{p,q}$ have $1$-dimensional centers, which must correspond to the circle
groups acting on the $1$-spheres $\lk(\{p,d\})$ and $\lk(\{p,q\})$. 
If $n=2$, then $\fa\cong\RR$ is not semisimple.
But since we have $G_{p,q}=\SO(2)\times\mathrm{O}(2)$ by \ref{StructureOfG_p},
the Lie algebra of the circle
group $(K_q)^\circ$ is distinguished in $\fg_{p,q}$ by the fact that 
$\mathrm{Ad}(G_{p,q})$ acts
nontrivially on it. The same applies to $\fg_{p,d}$.
Thus, $\Lie(K_d)$ and $\Lie(K_q)$ are in any case
distinguished subalgebras, and the possible isomorphisms
\begin{diagram}[height=1.5em]
\fd_{p}&&\lTo&&\fd_{p,q}&& \fq_{p,d}&&\rTo&&\fq_p\\
&&\uTo^\cong&&&&&& \uTo^\cong \\
\Lie(K_d)\oplus\fa&&\lTo&&\fa &&\fa&&\rTo&&\Lie(K_q)\oplus\fa .
\end{diagram}
are parametrized by nonzero reals.
However, all choices of these parameters lead to isomorphic diagrams for the Lie groups.
All other identification maps in $\Lie(\cG)$
are also unique up to automorphisms.
\qed
\end{Lem}
Recall from \ref{LemmaForn>m=1} that $\cG^\circ$ denotes the subdiagram of $\cG$
consisting of the identity components of the stabilizers.
\begin{Cor}
\label{UniqueDiagramForIdentityComponents}
Suppose that $n\geq 2$. Corresponding to each diagram
\begin{diagram}[height=1.5em]
\fq_d&\rTo&\fq\\
\uTo&&\uTo\\
\fq_{p,d}&\rTo&\fq_p
\end{diagram}
as in \ref{StructureOfQForm=1}, there is, up to isomorphism,
at most one possibility for the diagram $\cG^\circ$.

\proof
This follows from \ref{StructureOfG_p}, \ref{LemmaA} and \ref{LemmaForn>m=1}.
\qed
\end{Cor}
We identify $K_d\lTo K_{d,q}\rTo K_q$ with the matrix groups
\begin{diagram}
\left(\begin{smallmatrix}
* &   &   \\ 
  & * & * \\
  & * & *
\end{smallmatrix}\right)
&\lTo&
\left(\begin{smallmatrix}
* &   &   \\ 
  & * &   \\
  &   & *
\end{smallmatrix}\right)
&\rTo&
\left(\begin{smallmatrix}
* & * &   \\ 
* & * &   \\
  &   & *
\end{smallmatrix}\right)
\end{diagram}
in $\SO(3)$ and we put
\[\textstyle
u=\left(
\begin{smallmatrix}
1 &   &   \\ 
  & -1 &   \\
  &   & -1
\end{smallmatrix}\right)\qquad\text{ and }\qquad
v=\left(
\begin{smallmatrix}
-1 &   &   \\ 
  & -1 &   \\
  &   &\phantom{-} 1
\end{smallmatrix}\right) .
\]
\begin{Prop}
Suppose that $n\geq 2$. Corresponding to each diagram
\begin{diagram}[height=1.5em]
\fq_d&\rTo&\fq\\
\uTo&&\uTo\\
\fq_{p,d}&\rTo&\fq_p
\end{diagram}
as in \ref{StructureOfQForm=1}, there is, up to isomorphism,
at most one possibility for the diagram $\cG$.

\proof
There is a unique possibility for the diagram $\cG^\circ$ by
\ref{UniqueDiagramForIdentityComponents}. The group $H=(G_q)^\circ$
induces the group $L^\circ$ on $\lk(q)$ and acts transitively on the 
chambers of this link. From \ref{ConnDiagram}
we see that $H_p=(G_{p,q})^\circ$, whence
$H_{p,d}=H_p\cap G_{p,d,q}$ and
$H_d=H_{p,d}(G_{d,q})^\circ$.
Thus the action of $(G_q)^\circ$ on $\lk(q)$ is uniquely determined by $\cG^\circ$,
and so is the kernel $B$. Similarly, the action of $(G_d)^\circ$ on $\lk(d)$
and the kernel $C$ are uniquely determined.

We note that, by \ref{ConnDiagram}, the group $G_q$ is generated by
$H\cup\{u \}$. We know that $u$ acts on the $1$-sphere $\lk(\{p,q\})$
as a reflection. On the other hand, $u$ is contained in $(G_d)^\circ$
and therefore we know in particular how it acts on $\lk(\{d,q\})$.
Thus, the image of $u$ in $L$ is uniquely determined, and hence $L$
is uniquely determined. 

We put $B=\{1,z\}$. If $z=1$, then $L=G_q$ is uniquely determined, and so 
is the image of $u$ in $L$. If $z\neq 1$, then
we apply \ref{HomologyLemma} to the problem 
\begin{diagram}[height=1.5em]
(G_q)^\circ & \rDotsto& G_q\\
\dTo &&\dDotsto\\
L^\circ&\rTo& L .
\end{diagram}
By \ref{HomologyLemma} and the remarks following it,
there are two possibilities for the multiplication on $G_q$.
We know that $u^2=1$. One of the two possible multiplications
would give us $u*u=(uz)*(uz)=z$, which is wrong. So the correct multiplication
on $G_q$ is uniquely determined. There are two possible 
targets for $u$ in $G_q$, differing by $z$,
which act in the same way on $\lk(q)$. The element $u$ acts trivially
on $\lk(\{p,d\})$, while the product $uz$ acts 
as a reflection on $\lk(\{p,d\})$, hence
we know also the correct image of $u$ in $G_q$. This determines the 
map $K_q\times A\rTo G_q$ uniquely.

A completely analogous discussion shows that there is a unique possibility
for $G_d$ and the map $K_d\times A\rTo G_d$. In particular, there is
a unique possibility for $G_d\lTo K_{d,q}\times A\rTo G_q$.
The diagram $\cG$ is now uniquely determined.
\qed
\end{Prop}
It remains to consider the case $n=1$. We have seen in \ref{LemmaA}
that there are precisely two possibilities for $\Lie(\cG)$. One is realized 
in the polar space associated to the symmetric bilinear form 
$(f_{-3})\oplus f_4$ on $\RR^{3+4}$ and the associated polar representation
of $\SO(3)\times\SO(4)$ on $\RR^{3\times 4}$. 
It is analogous to the case $n>1$ considered
before, and we call this the \emph{orthogonal situation}. 

The other possibility is associated to the polar space corresponding to the standard
symplectic form on $\RR^6$. The associated polar representation is
$\mathrm U(3)/\{\pm1\}$ acting on the space of complex symmetric $3\times 3$-matrices,
via $(g,X)\mapstoo gXg^T$. In this action the group $G_p$ is the stabilizer of
{\tiny$\left(\begin{smallmatrix} 1&&\\&1&\\&&1    \end{smallmatrix}\right)$},
the group $G_q$ is the stabilizer of 
{\tiny$\left(\begin{smallmatrix} 1&&\\&0&\\&&0    \end{smallmatrix}\right)$},
and the group $G_d$ is the stabilizer of
{\tiny$\left(\begin{smallmatrix} 1&&\\&1&\\&&0    \end{smallmatrix}\right)$}.
We call this the \emph{symplectic situation}.

The two generalized quadrangles that may appear as the link at $q$ 
are dual to each other (isomorphic under a not type-preserving 
simplicial isomorphism).
The connected component of $L=G_q/B$ is 
\[
L^\circ=\SO(2)\times\SO(3).
\]
In the orthogonal situation, $(L^\circ)_p\cong\SO(2)$ acts with a two-element
kernel on $\lk(\{p,q\})$, while $(L^\circ)_d\cong\mathrm O(2)$ acts faithfully
on $\lk(\{d,q\})$. The $L^\circ$-stabilizer of $\gamma=\{p,d,q\}$
acts trivially on $\lk(\{p,q\})$, and as a reflection on $\lk(\{d,q\})$. 

In the symplectic situation, it is the other way around.

We note that in both cases $u$ becomes trivial in 
$\pi_0(G_{p,d})$, $\pi_0(G_p)$ and $\pi_0(G_d)$,
and that $v$ becomes trivial in $\pi_0(G_{p,q})$, $\pi_0(G_p)$ and $\pi_0(G_q)$.
The element $v$ is not trivial in $\pi_0(G_d)$, because its action 
on $\SS^1\times\SS^1$ is not orientation preserving.
Since we have a product action of $M^\circ$ on $\lk(d)$, the element $u$ acts
trivially on $\lk(d)$, and in particular $C=\{1,u\}\subseteq (G_d)^\circ$.
\begin{Num}[\boldmath In the symplectic situation $\cG$ is unique]
The circle group $(K_q)^\circ$ acts with kernel $\{1,v\}$ on $\lk(\{p,q\})$.
The group $(L_p)^\circ$, which must be its image, acts faithfully on $\lk(\{p,q\})$.
Therefore we have $B=\{1,v\}\subseteq (G_q)^\circ$. 

We claim that $A=1$. Suppose to the contrary that $1\neq a\in A$.
Then $a$ is nontrivial in $\pi_0(G_p)$ by \ref{StructureOfG_p}.
It is nontrivial in $\pi_0(G_d)$, since its action on $\SS^1\times\SS^1$ is not
orientation preserving. 
Also, it acts differently than the $L^\circ$-stabilizer
of $\gamma$, hence $a$ is nontrivial in $\pi_0(G_q)$. It follows now easily that 
$\cG$ admits a simple homomorphism onto $\ZZ/2$, hence $G$ is by \ref{UniversalProperty}
and the remark preceding \ref{IsotropyConnectedLemma} not connected, a contradiction.

Therefore $A=1$ and $\pi_0(G_\gamma)\cong\ZZ/2\oplus\ZZ/2$ has $u,v$ as a 
$\ZZ/2$-basis.
From the action of $\bra{u,v}=G_\gamma$ on $\lk(\{p,q\})\cup\lk(\{d,q\})$ we see that 
$G_q/B=L=\SO(3)\times\SO(2)$ is connected. Since $B\subseteq (G_q)^\circ$,
the group $G_q$ is also connected. 
Since $\fq_p$ is contained in the simple part $\fq'=[\fq,\fq]\cong\fso(3)$
and since the corresponding connected circle group contains the 
nontrivial kernel $B$, we have $[G_q,G_q]\cong\SU(2)$ and we may identify
$(K_q)^\circ$ with the subgroup $\SO(2)\subseteq\SU(2)$.
The group $\SU(2)\times\SO(2)$ contains no involution $u$ that
normalizes $\SO(2)$ and acts by inversion. Thus
$G_q=\mathrm U(2)$ and $G_{p,q}=\mathrm O(2)\subseteq\mathrm U(2)$, embedded in the standard
way as the group of elements fixed by complex conjugation.
The group $(G_{p,d})^\circ$ is determined by its Lie algebra, and
$G_{p,d}=G_{p,d,q}(G_{p,d})^\circ$.
This group can be identified with
$\left(\begin{smallmatrix}\mathrm O(1)\\&\mathrm U(1)\end{smallmatrix}\right)\subseteq\mathrm U(2)$.
The image of $v$ in $G_q$, 
being in the kernel $B$, is $\left(\begin{smallmatrix}\mathrm -1\\&-1\end{smallmatrix}\right)$.
The action of $(G_{p,d})^\circ$ on $\lk(\{d,q\})$ has in its kernel 
the element $\left(\begin{smallmatrix}\mathrm 1\\&-1\end{smallmatrix}\right)$.
Since $(G_{p,d})^\circ$ acts trivially on $\lk(\{p,d\})$, we conclude that this matrix is the
image of $u$ (and not of $uv$, which acts nontrivially on $\lk(\{d,q\})$) in $G_q$. 
Thus we know $G_q$ and its stabilizers, and the homomorphism $K_q\rTo G_q$.
It remains to determine $G_d$. But $G_d$ is the quotient of $G_{p,d}\times G_{d,q}$,
where we identify the respective images of $u$, $v$ and $uv$.
The diagram $\cG$ is now completely determined.
\qed
\end{Num}
\begin{Num}[\boldmath In the orthogonal situation $\cG$ is unique]
The circle group $(K_q)^\circ$ acts with kernel $\{1,v\}$ on $\lk(\{p,q\})$.
The group $(L_p)^\circ$, which must be its image, acts also with a $2$-element kernel 
on $\lk(\{p,q\})$. 
The element $v$ acts therefore as a reflection on $\lk(\{d,q\})$ and on $\lk(\{p,d\})$.

Let $Q\subseteq G_q$ denote the connected normal subgroup with Lie algebra $\fso(3)$.
We claim that $Q=\SO(3)$. Suppose to the contrary that $Q=\SU(2)$, with center $\{1,z\}$.
The circle group $(G_{d,q})^\circ\subseteq Q$ contains $z$. Since we have a product action on
$\lk(d)$, the element $z$ acts trivially on $\lk(d)$, and it acts trivially on $\lk(q)$.
This contradicts \ref{EffectiveOnE_1}. Thus $Q=\SO(3)$ has trivial center. It follows that 
$(G_q)^\circ=\SO(3)\times\SO(2)$. We have $G_q/G_{p,q}\cong\RRP^3$.
From the exact sequence
\[
 1\rTo
 \underbrace{\pi_1(G_{p,q})}_\ZZ\rTo 
 \underbrace{\pi_1(G_q)}_{\ZZ\oplus\ZZ/2}\rTo
 \underbrace{\pi_1(\RRP^3)}_{\ZZ/2}\rTo\pi_0(G_{p,q})\rTo\pi_0(G_q)\rTo 1
\]
we have an isomorphism $\pi_0(G_{p,q})\cong\pi_0(G_q)$. 
If $A\neq 1$, then we see from the action on $\mathcal E_1(\{p,d,q\})$ that
$G_{p,d,q}$ has $8$ elements, and from the action of $G_{p,q}$ on 
$\lk(\{p,q\})$ that $\pi_0(G_{p,q})\cong\ZZ/2\oplus\ZZ/2$. As in the symplectic
case, we conclude that $\pi_0(\cG)$ has a nontrivial simple homomorphism to $\ZZ/2$,
which is impossible. Therefore $A=1$ and $G_{p,d,q}=\bra{u,v}$ has $4$ elements.
From the action of $G_{p,d,q}$ on $\mathcal E_1(\{p,d,q\})$ we see that $B=1$.
Also, we know $(K_q)^\circ\rTo G_q$. Since $B=1$, there is a unique target for
$u$ in $G_q$. This determines $K_q\rTo G_q$ completely. Also, $G_{d,q}$ is now 
determined by its Lie algebra and by $K_{d,q}\rTo G_q$.

Finally, $(G_{p,d})^\circ\cap(G_{d,q})^\circ=1$, since $u$ is not contained in
$G_{d,q}$, hence $(G_d)^\circ=(G_{p,d})^\circ\times(G_{d,q})^\circ$.
Also, we have $G_d\subseteq G_{p,d}\times G_{d,q}$ and we know the image of
$v$ in the first factor. The image in the second factor is uniquely determined
by the action of $v$ on $\lk(\{d,q\})$, a reflection, since $G_{d,q}$
acts faithfully on $\lk(\{d,q\})$. Thus, the target of $v$ in $G_d$
is uniquely determined, and this determines the remaining homomorphisms in $\cG$.
The diagram $\cG$ is now completely determined.
\qed
\end{Num}
{\em Proof of Theorem \ref{C3ClassificationTheorem}.}
In the previous sections we have determined, up to isomorphism,
all possibilities for a simple complex of compact groups $\cG$
arising from a homogeneous compact geometry in $\mathbf{HCG}(\mathsf C_3)$,
with $G$ being minimal. Each example that we found is either realized
by a rank $3$ polar space (a building), or by the exceptional
$\mathsf C_3$ geometry. 
\qed

\section{Consequences and applications}

In this last section we show first that in homogeneous compact 
geometries of higher rank, no
exceptions occur.
\begin{Lem}
\label{F4Lemma}
Suppose that $(G,\Delta)$ is a homogeneous compact geometry 
in $\mathbf{HCG}(\mathsf F_4)$.
Then $\Delta$ is continuously and equivariantly 
$2$-covered by a compact connected
Moufang building of type $\mathsf F_4$.

\proof
In the exceptional geometry of type $\mathsf C_3$ from Section~\ref{TheExceptionalC3Geometry},
the panels have dimensions $2$ and $3$. In the geometry $\Delta$, however, all panels
belong to compact Moufang planes and have therefore dimensions
$1$, $2$, $4$, or $8$ by \ref{BasicStrucureOfRank2Links}. 
By \ref{C3ClassificationTheorem},
every link of type $\mathsf C_3$ in $\Delta$
is covered by a $\mathsf C_3$ building. By Tits' Theorem
\ref{TitsCoveringTheorem}, there exists a building $\widetilde\Delta$
and a $2$-covering $\rho:\widetilde\Delta\rTo\Delta$.
By \ref{BuildingCoveringsAreContinuous}, 
the building $\tilde\Delta$ can be topologized in such a way that
$\rho$ becomes equivariant and continuous, 
and $\widetilde\Delta$ is the compact Moufang
building associated to a simple noncompact Lie group.
\qed
\end{Lem}
\begin{Lem}
\label{C4Lemma}
Suppose that $(G,\Delta)$ is a homogeneous compact geometry 
in $\mathbf{HCG}(\mathsf C_4)$.
Then $\Delta$ is continuously and equivariantly $2$-covered by a compact connected
Moufang building of type $\mathsf C_4$.

\proof
We label the vertex types as follows.
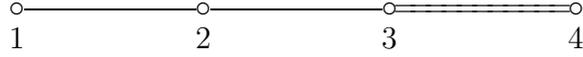
\begin{diagram}[height=1em,abut]
\circ & \rLine & \circ & \rLine & \circ & \rEq & \circ\\
1&&2&&3&&4
\end{diagram}
Let $\gamma=\{v_1,v_2,v_3,v_4\}$ be a chamber, where $v_i$ has type $i$.
We have to show that $\lk(v_1)$ cannot be the exceptional 
$\mathsf C_3$ geometry from Section~\ref{TheExceptionalC3Geometry}. 
Assume to the contrary that this is the case.
For $\emptyset\neq\alpha\subseteq\gamma$, we put $\fg_\alpha=\Lie(G_\alpha)$
and we let $\fn_\alpha\unlhd\fg_\alpha$ denote the Lie algebra of the
kernel of the action on $\lk(\alpha)$. Finally, we put 
$\fh_\alpha=\fg_\alpha/\fn_\alpha$. This is the Lie algebra of
the group induced by $G_\alpha$ on $\lk(\alpha)$.

We have $\fh_{v_1}=\fsu(3)\oplus\fsu(3)$ and,
by \ref{2*1+1}, we have $\fh_{v_1,v_3}=\fso(4)$, corresponding to the
exceptional action of $\SO(4)$ on $\SS^2\times\SS^3$.

Now we consider $\fh_{v_3}\subseteq\fsu(3)\oplus\fso(4)$.
The projection $\pr_1:\fh_{v_3}\rTo\fsu(3)$ to the first factor is onto, 
since  $\mathrm{PSU}(3)$ has no chamber-transitive closed subgroups.
Let $\fh_2$ denote the kernel of this projection, and let
$\fh_1\cong\fsu(3)$ be a supplement of the kernel,
$\fh_{v_3}=\fh_1\oplus\fh_2$. Every homomorphism from
$\fsu(3)$ to $\fso(4)$ is trivial, hence $\fh_1$ is 
the kernel of the projection $\pr_2:\fh_{v_3}\rTo\fso(4)$. 
The Lie algebra $\fh_{v_3}$ splits therefore in its action on $\lk(v_3)$
as a direct sum.
It follows that the Lie algebra $\fh_{v_1,v_3}$ splits also in its 
action on $\lk(\{v_1,v_3\})$.
We have reached a contradiction.

As in the previous lemma, we conclude from
\ref{TitsCoveringTheorem} and \ref{BuildingCoveringsAreContinuous} 
that there exists a compact building 
$\widetilde\Delta$ associated to a simple noncompact Lie group and a
continuous equivariant covering $\rho:\widetilde\Delta\rTo\Delta$.
\qed
\end{Lem}
The following two theorems summarize the main results of our classification.
\begin{Thm}
\label{MainTheoremInTheIrreducibleCase}
Let $M$ be a spherical irreducible Coxeter matrix of rank at least $4$
and suppose that $(G,\Delta)$ is a homogeneous compact geometry 
in $\mathbf{HCG}(M)$. Then there exists a compact connected spherical
building $\widetilde\Delta$ and a continuous $2$-covering
$\rho:\widetilde\Delta\rTo\Delta$.

\proof
By the previous two lemmata and by induction we see that the link of every
vertex is $2$-covered by a building. The claim follows now as in the 
proof of \ref{F4Lemma}.
\qed
\end{Thm}
The next theorem is an immediate consequence of 
\ref{JoinDecompositionLemma}, \ref{MainTheoremInTheIrreducibleCase} and
\ref{C3ClassificationTheorem}. It contains the Theorem~A of the 
introduction as a special case.
\begin{Thm}
\label{MainTheoremInTheReducibleCase}
Let $M$ be a Coxeter matrix of spherical type and let
$(G,\Delta)$ be a homogeneous compact geometry in $\mathbf{HCG}(M)$.
Suppose that the Coxeter diagram of $M$ has no isolated nodes.
Then there exists a homogeneous compact geometry $(K,\widetilde\Delta)$
in $\mathbf{HCG}(M)$ which is a join of buildings associated to simple noncompact Lie groups
and geometries of type $\mathsf C_3$ which are isomorphic to the exceptional
geometry in Section~\ref{TheExceptionalC3Geometry}, and a continuous $2$-covering
$\widetilde\Delta\rTo\Delta$, which is equivariant with respect to a
compact connected Lie group $K$ acting transitively on the chambers of
$\widetilde\Delta$.

\proof
We decompose $\Delta$ as a join $\Delta=\Delta_1*\Delta_2*\cdots*\Delta_m$
of irreducible factors, and we let $H_i$ denote the group induced by $G$ on 
$\Delta_i$. Now we apply \ref{C3ClassificationTheorem} and
\ref{MainTheoremInTheIrreducibleCase} to the homogeneous compact geometries
$(H_i,\Delta_i)$. We obtain equivariant $2$-coverings 
$\widetilde\Delta_i\rTo\Delta_i$, where $\widetilde\Delta_i$ is either a compact
building or the exceptional $\mathsf C_3$ geometry from Section~\ref{TheExceptionalC3Geometry}.
Taking the join of these $2$-coverings, we obtain the result that we claimed.
We note that the group induced by $K$ on $\Delta$ may be strictly larger than
the group $G$ we started with.
\qed
\end{Thm}
Recall from Section \ref{Section3} 
that an isometric group action $G\times X\rTo X$ on a complete Riemannian manifold $X$ is
called \emph{polar} if it admits a \emph{section} $\Sigma\subseteq X$, 
i.e. a complete totally geodesic submanifold
that intersects every orbit perpendicularly. A polar action is called
\emph{hyperpolar} if the section $\Sigma$ is flat.
One motivation for the present work is a recent result by the
second author \cite{Lytchak}.
Relying on the classification of buildings by Tits and Burns-Spatzier
\cite{BS}, the main theorem of \cite{Lytchak} contains
a classification of polar foliations on symmetric spaces of compact type
under the assumption that the irreducible parts of the foliation have codimension
at least $3$.
We refer to \cite{Lytchak} for the history of the subject and an extensive list of
literature about this problem. Our main result now
covers the remaining cases of codimension $2$,
provided that the polar foliation arises
from a polar action. The result is as follows.
\begin{Thm} \label{cons}
Suppose that $G\times X\rTo X$ is a polar action of a compact connected Lie group $G$
on a symmetric space $X$ of compact type. 
Then, possibly after replacing $G$ by a larger orbit equivalent
group, we have splittings $G=G_1 \times\cdots\times G_m$ and $X=X_1  \times\cdots\times X_m$, 
such that the action of $G_i$ on $X_i$ is either trivial or hyperpolar or the
space $X_i$ has rank $1$, for $i=1,\ldots m$.  
\end{Thm}
We indicate briefly the connection between our main theorem and \ref{cons} and 
refer the interested reader to \cite{Lytchak}.
Given a polar action on a symmetric space of compact type, the de Rham decomposition of 
a section of the action
gives rise to an equivariant product decomposition of the whole space. Removing the 
hyperpolar pieces (corresponding to the flat factor of the section) and the pieces that
do not admit reflection groups (corresponding to trivial actions), 
one is left with polar actions whose sections have constant positive curvature.
Moreover, one can split off another factor with a trivial
action of our group $G$, unless any two points of the manifold can be connected by a
sequence of points $p_1,\ldots,p_r$, where each consecutive pair $p_i,p_{i+1}$
is contained in a section.  
All these decomposition results rely heavily on results by Wilking \cite{Wilk}.
Then it remains to show that in this case  the symmetric space has  rank $1$. 

In order to do this, one observes that each section is a sphere or a projective space,
and that  the quotient space of the action
is isometric to the quotient of the universal covering of a section by a finite Coxeter
group.  If the Coxeter group is reducible,
the decomposition of the Coxeter group implies the existence of very special ``polar'' 
submanifolds of our symmetric space and from this
one deduces that the rank must be $1$.  
In the irreducible case, the quotient is a Coxeter simplex and each section is
decomposed by such Coxeter simplices.  
Taking all these simplices from all sections together one finds
a huge polyhedral complex. This polyhedral complex turns out 
to be a geometry of spherical type, each link of which is a spherical building
(defined by the corresponding slice representations). 

Thus we have found a homogeneous compact geometry of spherical type.
If the geometry is covered by a building (which is always the case if the geometry
is not of type $\mathsf C_3$) then
the covering complex is a Moufang building $\Delta$ belonging to a simple noncompact Lie group
and the manifold we started with is the base of a principal bundle
with total space homeomorphic to a sphere (the geometric realization
$|\Delta|_K$ in the coarse topology).  
Thus $X$ turns out to be of rank $1$ in this case.  
The  $\mathsf C_3$ case cannot be handled in this way,
since the Cayley plane $\OOP^2$ is not the quotient of a free action
of a compact group on a sphere. Indeed, the
exceptional geometry described in Section~\ref{TheExceptionalC3Geometry} 
cannot arise in this way. However, \ref{C3ClassificationTheorem}
says that the Podest\`a-Thorbergsson example is the only exception. 
\qed

\smallskip
Using the methods of \cite{Kol1}, the above  extension of \cite{Lytchak} was obtained
in \cite{KolLyt} under the additional assumption that $X$ is irreducible.  
Kollross has  announced an independent  extension of \cite{KolLyt} to the reducible case, 
with the methods of \cite{Kol1} (unpublished).

Key ideas and methods of \cite{Lytchak} were discovered and used independently by
Fang-Grove-Thorbergsson in their classification of polar actions in
positive curvature \cite{FGTh}. In particular, the classification of compact
Moufang buildings by Tits and Burns-Spatzier, and Tits' local approach to
buildings play crucial roles. With a few exceptions, the classification in
cohomogeneity two involving $\mathsf C_3$ geometries is based on an axiomatic
characterization by Tits. Of course our main result here can also be used
for this purpose.

{\bigskip
\raggedright
Linus Kramer\\
Mathematisches Institut, 
Universit\"at M\"unster,
Einsteinstr. 62,
48149 M\"unster,
Germany\\
{\makeatletter
e-mail: {\tt linus.kramer{@}uni-muenster.de}}\\
\smallskip
Alexander Lytchak\\
Mathematisches Institut,
Universit\"at K\"oln,
Weyertal 86-90,
50931 K\"oln, Germany\\
{\makeatletter
e-mail: {\tt alytchak@math.uni-koeln.de}}\\}

\end{document}